\documentclass[12pt,leqno]{article}
\usepackage[mathscr]{eucal}
\usepackage{amsmath,amssymb,latexsym,theorem,bbm,amsfonts}
\usepackage{verbatim}
\usepackage{color}
\usepackage{appendix}
\usepackage{array,epsfig,theorem,bbm,graphics,booktabs,url}

\newcommand{\CC}{\mathbb{C}}

\newcommand{\NN}{\mathbb{N}}
\newcommand{\QQ}{\mathbb{Q}}
\newcommand{\RR}{\mathbb{R}}

\newcommand{\SSS}{\mathbb{S}}

\newcommand{\ZZ}{\mathbb{Z}}

\newcommand{\bc}{{\boldsymbol{c}}}

\newcommand{\tg}{\widetilde{g}}

\newcommand{\tI}{\widetilde{I}}

\newcommand{\bv}{{\boldsymbol{v}}}

\newcommand{\bx}{{\boldsymbol{x}}}

\newcommand{\by}{{\boldsymbol{y}}}
\newcommand{\bY}{{\boldsymbol{Y}}}

\newcommand{\bzero}{{\boldsymbol{0}}}
\newcommand{\bone}{{\boldsymbol{1}}}

\newcommand{\cB}{{\mathcal B}}

\newcommand{\cD}{{\mathcal D}}

\newcommand{\cF}{{\mathcal F}}

\newcommand{\cM}{{\mathcal M}}

\newcommand{\cS}{{\mathcal S}}

\newcommand{\tX}{\widetilde{X}}

\newcommand{\cZ}{{\mathcal Z}}

\newcommand{\dd}{\mathrm{d}}
\newcommand{\ee}{\mathrm{e}}
\newcommand{\ff}{\mathrm{f}}
\newcommand{\ii}{\mathrm{i}}

\newcommand{\EE}{\operatorname{\mathbb{E}}}

\newcommand{\PP}{\operatorname{\mathbb{P}}}

\newcommand{\OO}{\operatorname{O}}

\newcommand{\sign}{\operatorname{sign}}

\renewcommand{\Re}{\operatorname{Re}}
\renewcommand{\Im}{\operatorname{Im}}

\newcommand{\talpha}{\widetilde{\alpha}}

\newcommand{\vare}{\varepsilon}

\renewcommand{\mid}{\,|\,}

\renewcommand{\leq}{\leqslant}
\renewcommand{\geq}{\geqslant}

\newcommand{\distr}{\stackrel{\cD}{\longrightarrow}}
\newcommand{\distrf}{\stackrel{\cD_\ff}{\longrightarrow}}
\newcommand{\distre}{\stackrel{\cD}{=}}
\newcommand{\distrw}{\stackrel{w}{\longrightarrow}}
\newcommand{\distrv}{\stackrel{v}{\longrightarrow}}

\newcommand{\bbone}{\mathbbm{1}}

\newcommand{\nt}{{\lfloor nt\rfloor}}

\newcommand{\proofend}{\hfill\mbox{$\Box$}}

\numberwithin{equation}{section}

\theoremstyle{change} \theorembodyfont{\em}
\newtheorem{Lem}{Lemma.}[section]
\newtheorem{Thm}[Lem]{Theorem.}
\newtheorem{Pro}[Lem]{Proposition.}
\newtheorem{Cor}[Lem]{Corollary.}
\newtheorem{Def}[Lem]{Definition.}

\theorembodyfont{\rm}
\newtheorem{Rem}[Lem]{Remark.}

\begin{document}

\begin{center}
 {\bfseries\Large Convergence of partial sum processes to stable processes \\[2mm] with application for aggregation of branching processes}

\vskip0.5cm

 {\sc\large
  M\'aty\'as $\text{Barczy}^{*,\diamond}$,
  Fanni K. $\text{Ned\'enyi}^{*}$,
  Gyula $\text{Pap}^{**}$}

\end{center}

\vskip0.2cm

\noindent
 * MTA-SZTE Analysis and Stochastics Research Group,
   Bolyai Institute, University of Szeged,
   Aradi v\'ertan\'uk tere 1, H--6720 Szeged, Hungary.

\noindent
 ** Bolyai Institute, University of Szeged,
     Aradi v\'ertan\'uk tere 1, H--6720 Szeged, Hungary.

\noindent e--mails: barczy@math.u-szeged.hu (M. Barczy),
                    nfanni@math.u-szeged.hu (F. K. Ned\'enyi).

\noindent $\diamond$ Corresponding author.

\vskip0.2cm


\renewcommand{\thefootnote}{}
\vspace*{0.2cm}
\footnote{M\'aty\'as Barczy is supported by the J\'anos Bolyai Research Scholarship of the Hungarian Academy of Sciences.
Fanni K. Ned\'enyi is supported by the UNKP-18-3 New National
Excellence Program of the Ministry of Human Capacities.
Gyula Pap was supported by grant NKFIH-1279-2/2020 of the Ministry for Innovation and Technology, Hungary.}

\vspace*{-10mm}

\begin{abstract}
We provide a generalization of Theorem 1 in Bartkiewicz et al.\ \cite{BarJakMikWin} in the sense that we give sufficient conditions for weak convergence of finite dimensional
 distributions of the partial sum processes of a strongly stationary sequence to the corresponding finite dimensional distributions of a
 non-Gaussian stable process instead of weak convergence of the partial sums themselves to a  non-Gaussian stable distribution.
As an application, we describe the asymptotic behaviour of finite dimensional distributions of aggregation of independent copies of a strongly stationary subcritical
 Galton--Watson branching process with regularly varying immigration having index in \ $(0, 1) \cup (1, 4/3)$ \
 in a so-called iterated case, namely when first taking the limit as the time scale and then the number of copies tend to infinity.
\end{abstract}

\section{Introduction}
\label{intro}

Studying asymptotic behaviour of appropriately centered and scaled partial sum processes
 \ $\bigl(\sum_{k=1}^\nt Y_k\bigr)_{t\in[0,\infty)}$, \ $n = 1, 2, \ldots$, \ of a strongly stationary sequence
 \ $(Y_k)_{k\in\{0,1,\ldots\}}$ \ to a non-Gaussian stable process has a long history, see,
 e.g., Resnick \cite{Res}, Beran et al. \cite{BerFenGhoKul} and Kulik and Soulier \cite[Chapter 15]{KulSou}.
In this paper we provide a generalization of Theorem 1 in Bartkiewicz et al.\ \cite{BarJakMikWin} in the sense that
 we give sufficient conditions for weak convergence of finite dimensional distributions of the partial sum processes
 of a strongly stationary sequence to a non-Gaussian stable process instead of weak convergence of the partial sums themselves
 to a non-Gaussian stable distribution, see Theorem \ref{BJMW1}.
These sufficient conditions in question include mixing-type and anti-clustering type conditions (see (ii) and (iii) of Theorem \ref{BJMW1}.
 For strongly mixing, strongly stationary jointly regularly varying sequences we give some well-applicable sufficient conditions under which the
 mixing-type and anti-clustering type conditions in question hold, see Lemma \ref{BJMW1} and \ref{Lemma3}, respectively.
Theorem 1 in Bartkiewicz et al.\ \cite{BarJakMikWin} has been applied for stochastic recurrence equations,
 ARCH(1) and GARCH(1,1) processes and stable stationary sequences, see Bartkiewicz et al.\ \cite[Section 4]{BarJakMikWin}.
As an application of our Theorem \ref{BJMW1}, we describe the asymptotic behaviour of finite dimensional distributions of aggregation of independent copies
 of a strongly stationary subcritical Galton--Watson branching process \ $(X_k)_{k\in\{0,1,\ldots\}}$ \
 with regularly varying immigration having index in \ $(0, 1) \cup (1, 4/3)$ \
 in a so-called iterated case, namely when first taking the limit as the time scale and then the number of copies tend to infinity.
For this application, our generalization Theorem \ref{BJMW1} of Theorem 1 in Bartkiewicz et al.\ \cite{BarJakMikWin} is needed,
 since we consider weak convergence of finite dimensional distributions of the aggregated process in question, and not only those
 of the one-dimensional distributions.
In what follows we will simultaneously review the existing literature and present our results
 on the asymptotic behaviour of a strongly stationary sequence towards a non-Gaussian stable process
 and that on the aggregation of Galton--Watson branching processes with immigration.
In this way, we also motivate our generalization of Theorem 1 in Bartkiewicz et al.\ \cite{BarJakMikWin}.

Modelling the growth of a population in time, one may naturally be interested in the number of individuals
 who were born or joined as immigrants up to a given time point.
This quantity corresponds to the appropriate partial sum of the Galton--Watson branching process with immigration that is modelling the growth of  the population.
While keeping track of the growth of a single community is a meaningful task on its own, we are interested in doing so having multiple independent communities.
(For example one may consider a big area where many hives coexist without any interference.)
Therefore, it is important to study the total number of individuals who were born or joined as immigrants up to a given time point \ $n$ \ in some of the given \ $N$ \ populations.
We can do so by describing the asymptotic behaviour of this double sum as \ $n$ \ and \ $N$ \
 tend to \ $\infty$ \ in an iterated way or simultaneously.
In the iterated case there are two subcases, namely, first taking the limit \ $n\to\infty$ \ and then \ $N\to\infty$, \
 and vica versa.
In the literature the above mentioned procedure is called an iterated and simultaneous aggregation, respectively,
 see, e.g., Pilipauskait\.{e} and Surgailis \cite{PilSur}.
For a review of the literature on aggregation of time series and branching processes, see the Introduction of Barczy et al. \cite{BarNedPap3}.
Here we mention only our paper Barczy et al. \cite{BarNedPap2}, where we investigated the limit behaviour of the same aggregation scheme for a stationary multitype Galton--Watson branching process with immigration under third order moment conditions
 on the offspring and immigration distributions in the iterated and simultaneous cases as well.

First, we examine the asymptotic behaviour of the finite dimensional distributions of the temporal aggregates
 \ $\bigl(\sum_{k=1}^\nt X_k\bigr)_{t\in[0,\infty)}$ \ as \ $n\to\infty$, \ where \ $(X_k)_{k\in\{0,1,2,\ldots\}}$ \
  is a strongly stationary subcritical Galton-–Watson branching process with regularly varying immigration.
Under the assumptions that the index \ $\alpha$ \ of the regularly varying immigration is in \ $(0,2)$, \ and
 in case of \ $\alpha\in[1,2)$ \ additionally assuming that the offspring distribution admits a finite second moment,
 the strongly stationary process \ $(X_k)_{k\in\{0,1,\ldots\}}$ \ is jointly regularly varying, see Basrak et al. \cite[Theorem 2.1.1]{BasKulPal}
 (also Theorem \ref{Xtail}) and Basrak and Segers \cite[Theorem 2.1]{BasSeg}.
As we mentioned, in the literature one can find several results, mainly based on point processes and characteristic functions, for deriving convergence of appropriately centered and scaled partial sum processes \ $\bigl(\sum_{k=1}^\nt Y_k\bigr)_{t\in[0,\infty)}$, \ $n = 1, 2, \ldots$, \ of a strongly stationary sequence \ $(Y_k)_{k\in\{0,1,\ldots\}}$ \ to a non-Gaussian stable process, which have been successfully applied for time series, especially
 for autoregressive moving average sequences.
Much less is known about application for branching processes.
Basrak et al. \cite[Section 3.2]{BasKulPal} described the limit behaviour of appropriately normalized partial sums of \ $(X_k)_{k\in\{0,1,\ldots\}}$ \ in case of \ $\alpha \in (0, 1) \cup (1, 2)$.
\ In their work, they identified the limit distribution as an \ $\alpha$-stable distribution, however they did not present the characteristic function of this limit law.
Furthermore, in case of \ $\alpha \in (1, 2)$, \ they wrote that they need an additional technical condition which is not formulated in their paper, they refer to the vanishing small value condition (3.2) in Davis and Hsing \cite{DavHsi}.
We remove this additional technical condition in case of \ $\alpha\in(1,\frac{4}{3})$.
 \ Further, instead of weak convergence of appropriately normalized partial sums we can prove weak convergence of appropriately normalized
 finite dimensional distributions of the partial sum processes providing an explicit characteristic function of the
 limit distributions as well, see Theorem \ref{aggr_time}.

Roitershtein and Zhong \cite[Theorems 2.11 and 2.12]{RoiZho} described the asymptotic behaviour of partial sums of a strongly stationary first order random coefficient integer-valued autoregressive process (abbreviated as RCINAR(1)) with regularly varying immigration having index \ $\alpha \in (0, 2]$.
\ Note that a first order integer-valued autoregressive process is a Galton--Watson branching process with immigration admitting Bernoulli offspring distribution.
Roitershtein and Zhong \cite{RoiZho} showed that appropriately scaled and centered partial sums of the RCINAR(1) process in question converge in distribution to an \ $\alpha$-stable law.
Surprisingly, they can also handle the case \ $\alpha = 1$ \ which is usually excluded in papers on similar studies for time series (see, e.g.,
 Bartkiewicz et al. \cite[Propositions 3 and 4]{BarJakMikWin}).
We note that in case of \ $\alpha \in [1, 2]$, \ Roitershtein and Zhong \cite{RoiZho} do not provide any proof, they refer to a ''standard technique'' due to Kesten et al. \cite{KesKozSpi}.

To describe the asymptotic behaviour of finite dimensional distributions of \ $\bigl(\sum_{k=1}^\nt X_k\bigr)_{t\in[0,\infty)}$ \ as \ $n\to\infty$, \ in principle, we could use the results of Tyran-Kami\'nska \cite[Theorem 1.1]{Tyr}, Basrak et al. \cite[Theorem 3.4]{BasKriSeg} or Cattiaux and Manou-Abi \cite[Theorem 3.1]{CatMan}, but we were not able to check some of their vanishing small value type conditions in case of \ $\alpha \in [1, 2)$.
To be more detailed, Tyran-Kami\'nska \cite[Theorem 1.1]{Tyr} established a functional limit theorem for an appropriately centered and normalized partial sum process of
 a strongly stationary and strongly mixing stochastic process \ $(Y_k)_{k\in \{1,2,\ldots\}}$ \ such that \ $Y_1$ \ is regularly varying with index \ $\alpha\in (0,2)$.
\ In case of \ $\alpha\in[1,2)$ \ a vanishing small value type condition is needed (see (1.6) in Tyran-Kami\'nska \cite[Theorem 1.1]{Tyr}), which
 is hard to check in general.
In Lemma 4.8 and Remark 4.9 in \cite{Tyr}, sufficient conditions are given for the condition (1.6) in \cite{Tyr} in terms of the maximal correlation coefficients
 of \ $(Y_k)_{k\in\{1,2,\ldots\}}$.
Unfortunately, we were not able to check this condition (1.6) in case of \ $\alpha\in[1,2)$, \ nor to use the sufficient conditions given
 in Lemma 4.8 and Remark 4.9 in \cite{Tyr}.
We also note that in Theorem 1.1 in Tyran-Kami\'nska \cite{Tyr}, the characteristic function of the limit $\alpha$-stable law remains somewhat hidden.

Basrak et al. \cite[Theorem 3.4]{BasKriSeg} formulated another functional limit theorem for an appropriately centered
 and scaled partial sum process of a strongly stationary, jointly regularly varying stochastic process under some conditions.
In case of \ $\alpha\in[1,2)$, \ similarly to Tyran-Kami\'nska \cite[Theorem 1.1]{Tyr}, they also suppose
 a vanishing small value condition, which we could not check for a subcritical strongly stationary Galton--Watson branching process with regularly varying immigration.
In Section 3.5.2 in Basrak et al. \cite{BasKriSeg}, they discuss different centralizing constants such as truncated mean, the mean itself
 or zero centralizing constant.
In our paper we also present limit theorems with these three kinds of centralizing constants separately.

We overcome the above mentioned difficulties
 by deriving a generalization of Theorem 1 in Bartkiewicz et al. \cite{BarJakMikWin} on weak convergence of partial sums of strongly stationary
 jointly regularly varying sequences with index in \ $(0, 2)$ \ in the sense that we give some sufficient conditions under which
 the finite dimensional distributions of the corresponding partial sum processes converge weakly to those of a non-Gaussian stable process,
 see Theorem \ref{BJMW1}.
Our sufficient conditions are the same as those of Theorem 1 in Bartkiewicz et al. \cite{BarJakMikWin} except their (MX) mixing-type condition,
 which we adjusted for handling finite dimensional distributions.
More precisely, the conditions (i), (iii), (iv) and (v) of Theorem \ref{BJMW1} coincide with the conditions 1, 3, 4 and 5 of Theorem 1
 in Bartkiewicz et al.\ \cite{BarJakMikWin}, respectively.
Further, the condition (ii) of Theorem \ref{BJMW1} can be considered as a ''finite dimensional'' counterpart of the condition 2 ((MX) condition)
 of Theorem 1 in Bartkiewicz et al. \cite{BarJakMikWin}.
Note also that condition (ii) of Theorem \ref{BJMW1} in the special case \ $d=1$ \ gives back
 condition (MX) in Bartkiewicz et al.\ \cite{BarJakMikWin}.
In what follows, we discuss some sufficient conditions under which the condition (ii) of Theorem \ref{BJMW1} holds, and we also discuss
 their connections with the existing sufficient conditions for (MX) condition in Bartkiewicz et al.\ \cite{BarJakMikWin}.
The condition (MX) in Bartkiewicz et al.\ \cite{BarJakMikWin} and condition (ii) in Theorem \ref{BJMW1} are formulated
 in terms of characteristic functions of appropriately normalized partial sums, and, in general, it is not easy to check them.
However, Lemma 3 in Bartkiewicz et al.\ \cite{BarJakMikWin} contains a sufficient condition under which a strongly mixing, strongly stationary
sequence satisfies condition (MX), further, on page 360 in Bartkiewicz et al. \cite{BarJakMikWin} it is stated that a strongly mixing, strongly stationary sequence
 with geometric rate automatically satisfies condition (MX).
Unfortunately, we can not fully understand the proof of Lemma 3 in Bartkiewicz et al. \cite{BarJakMikWin},
 namely, we do not understand why the given definition of \ $U_{ji}$ \ is appropriate,
 and we do not understand why the inequality \ $\EE\big(\big \vert \frac{x}{a_n}\sum_{j=1}^{k_n} U_{\delta j}\vert \wedge 2\big)
 \leq \EE\big( \frac{x\delta n}{m a_n}\vert X_1\vert \wedge 2\big)$ \ holds at the second step of the
 estimation of \ $|\varphi_q(x) - \varphi_{nm\delta}(x)|$.
\ We also remark that Bartkiewicz et al. \cite{BarJakMikWin} wrote that their Lemma 3 yields
 that a strongly mixing, strongly stationary sequence with geometric rate automatically satisfies condition (MX) without giving any proof.
Because of these reasons, we decided to formulate a corresponding lemma in an extended form with a detailed proof, see Lemma \ref{Lemma3}.
Next, we discuss some sufficient conditions under which the condition (iii) of Theorem \ref{BJMW1} (which coincides with the (AC) condition
 in Bartkiewicz et al.\ \cite{BarJakMikWin}) holds.
Lemma 2 in Bartkiewicz et al. \cite{BarJakMikWin} gives several sufficient conditions under which the above mentioned (AC) condition holds.
However, we feel that the proof of Lemma 2 in Bartkiewicz et al.\ \cite{BarJakMikWin} contains some gaps, for example,
 we do not know how one can verify the estimates in lines 4-5 and 9 on page 352 in Bartkiewicz et al. \cite{BarJakMikWin}.
Due to this, for the applicability of Theorem \ref{BJMW1}, we decided to use another anti-clustering type condition \eqref{3.9}
 presented in Lemma \ref{Lemma2}, suggested by Thomas Mikosch.
Roughly speaking, the advantage of condition \eqref{3.9} compared to the (AC) condition in Bartkiewicz et al.\ \cite{BarJakMikWin}
 is that it involves only \ $Y_k$, $k\in\{0,1,\ldots\}$, \ instead of the partial sums $Y_1+\cdots+Y_n$, $n\in\{1,2,\ldots\}$,
 of a strongly stationary sequence \ $(Y_k)_{k\in\{0,1,\ldots\}}$, \ and hence it is easier to check in applications.

Based on Theorem 1 in Bartkiewicz et al. \cite{BarJakMikWin}, Mikosch and Wintenberger \cite[Theorem 4.1]{MikWin} gave sufficient conditions for weak convergence of partial sums
 of a strongly stationary, jointly regularly varying sequence which is a function of an irreducible aperiodic Markov chain.
Studying the proof of Theorem 4.1 in Mikosch and  Wintenberger \cite{MikWin}, we realized that in order to derive sufficient conditions for weak convergence of finite dimensional
 distributions in their setup, one needs our above mentioned generalization of Theorem 1 in Bartkiewicz et al. \cite{BarJakMikWin}.
Since this generalization of Theorem 1 in Bartkiewicz et al. \cite{BarJakMikWin} turned out to be enough for us,
 we have not worked out the generalization of Theorem 4.1 in Mikosch and Wintenberger \cite{MikWin}.

As an application of our Theorem \ref{BJMW1}, Lemmas \ref{Lemma3} and \ref{Lemma2}, we derive weak convergence of finite dimensional distributions
 of appropriately centered and scaled partial sum processes \ $\bigl(\sum_{k=1}^\nt X_k\bigr)_{t\in[0,\infty)}$ \ as \ $n \to \infty$ \
  in case of \ $\alpha \in (0, 1) \cup (1, \frac{4}{3})$, \ see Theorem \ref{aggr_time}.
We present limit theorems with three kinds of centralizing constants, namely with the truncated mean, the mean itself or zero centralizing constant.
We point out the fact that in case of \ $\alpha \in (1, \frac{4}{3})$ \ we managed to get rid of the additional vanishing small value condition mentioned
 in Basrak et al. \cite[Section 3.2]{BasKulPal}.
Note that Theorem \ref{BJMW1} can not be applied for \ $(X_n)_{n\in\{0,1,\ldots\}}$ \ in case of \ $\alpha = 1$ \ (see part (i) of Remark \ref{Rem_support}),
 which is the reason for not having any result in case of \ $\alpha=1$.
\ Further, in part (ii) of Remark \ref{Rem_support} we explain why our present technique is not suitable for handling the case \ $\alpha\in[\frac{4}{3},2)$.
\ As a consequence of Theorem \ref{aggr_time}, we describe the asymptotic behaviour of finite dimensional distributions of the aggregated stochastic process
 \ $\bigl(\sum_{j=1}^N \sum_{k=1}^{\lfloor nt \rfloor} X^{(j)}_k\bigr)_{t\in[0,\infty)}$ \ considering different centerings as first \ $n$ \
 and then \ $N$ \ converge to infinity, where \ $(X^{(j)}_k)_{k\in\{0,1,\ldots\}}$, \ $j=1,2,\ldots$, \ are independent copies of \ $(X_k)_{k\in\{0,1,\ldots\}}$,
 \ see Theorem \ref{iterated_aggr_2}.

The present paper is organized as follows.
Sections \ref{results} and \ref{application} contain our above detailed results.
Section \ref{Proofs} is devoted to the proofs.
We close the paper with six appendices.
Here we mention that in Appendix \ref{App_reg_var} we recall the definitions of regularly varying functions, regularly varying random vectors and some
 related results that we use in the proofs, among others  Karamata's theorem for truncated moments.
Further, in Appendix \ref{mixing} we present an auxiliary lemma stating that a strongly stationary, subcritical Galton--Watson
 branching process \ $(X_k)_{k\in\{0,1,\ldots\}}$ \ with regularly varying immigration (for the precise assumptions, see the beginning of Section \ref{application})
 \ is strongly mixing with geometric rate and is \ $V$-uniformly ergodic with \ $V(x):=1+x^p$, \ $x\in\{0,1,\ldots,\}$, \
 where \ $p\in(0,\min(\alpha,1))$, \ see also Basrak et al.\ \cite[Remark 3.1]{BasKulPal}.

Finally, we summarize the novelties of the paper.
We give a generalization of Theorem 1 in Bartkiewicz et al.\ \cite{BarJakMikWin} considering weak convergence of finite dimensional distributions instead
 of the one dimensional ones detailed above.
As an application, we successfully utilize this generalization for describing the asymptotic behaviour of finite dimensional distributions
 of the above mentioned aggregation scheme for a strongly stationary, subcritical Galton--Watson
 branching process \ $(X_k)_{k\in\{0,1,\ldots\}}$ \ with regularly varying immigration (for the precise assumptions, see the beginning of Section \ref{application}).
\ Aggregation of branching processes with low moment conditions is a relatively new topic in the field of aggregation of (randomized) stochastic processes,
 so our previously mentioned application can be considered as one of the first steps.
In case of \ $\alpha\in(1,\frac{4}{3})$ \ we managed to get rid of the vanishing small value condition assumed in Basrak et al. \cite[Section 3.2]{BasKulPal}
 for the description of the asymptotic behaviour of \ $\sum_{k=1}^n X_k$ \ as \ $n\to\infty$.

In a companion paper Barczy et al. \cite{BarNedPap3} we studied the other iterated, idiosyncratic aggregation scheme, namely, when first taking the limit
 as the number of copies and then as the time scale tend to infinity.
It turned out that the limit distributions are the same for the appropriately scaled and centered
 aggregated processes as in the present paper.

\section{Convergence of partial sum processes to stable processes}
\label{results}

Let \ $\ZZ_+$, \ $\NN$, \ $\QQ$, \ $\RR$, \ $\RR_+$, \ $\RR_{++}$ \ and \ $\CC$ \ denote the set of non-negative integers, positive
 integers, rational numbers, real numbers, non-negative real numbers, positive real numbers and complex numbers, respectively.
For \ $x,y\in\RR$, \ let \ $x\wedge y:=\min(x,y)$ \ and \ $x\vee y := \max(x,y)$.
The imaginary unit is denoted by \ $\ii$.
\ Convergence in distribution and equality in distribution of random variables will be denoted by \ $\distr$ \ and by \ $\distre$, \ respectively.
We will use \ $\distrf$ \ or \ $\cD_\ff\text{-}\hspace*{-1mm}\lim$ \ for weak
 convergence of finite dimensional distributions.
The definition and some related results for regularly varying random vectors are recalled in Appendix \ref{App_reg_var}.

If \ $Y$ \ is a regularly varying random variable with index \ $\alpha \in \RR_{++}$, \ then there exists a sequence \ $(a_n)_{n\in\NN}$ \ in \ $\RR_{++}$ \ with \ $n \PP(|Y| > a_n) \to 1$ \ as \ $n \to \infty$, \ see, e.g., Lemma \ref{a_n}.
In fact, \ $a_n = n^{1/\alpha} L(n)$, \ $n \in \NN$, \ for some slowly varying continuous function
 \ $L : \RR_{++} \to \RR_{++}$, \ see, e.g., Araujo and Gin\'e \cite[Exercise 6 on page 90]{AraGin}.

\begin{Def}
A stochastic process \ $(Y_k)_{k\in\NN}$ \ is called jointly regularly varying with index \ $\alpha \in \RR_{++}$ \ if its finite dimensional distributions are regularly varying with index \ $\alpha$.
\end{Def}

Given a stochastic process \ $(Y_k)_{k\in\NN}$, \ let \ $S_n:=Y_1+\cdots+Y_n$, $n\in\NN$, \ and
 for any \ $d \in \NN$ \ and \ $j_1, \ldots, j_d \in \ZZ_+$ \ with \ $j_1 \leq \ldots \leq j_d$, \ let us define
 \[
      \varphi_{n,j_1,\ldots,j_d}(\vartheta_1,\ldots,\vartheta_d) := \EE\biggl(\exp\biggl\{\frac{\ii}{a_n} \sum_{\ell=1}^d \vartheta_\ell (S_{j_\ell} - S_{j_{\ell-1}})\biggr\}\biggr) , \qquad \vartheta_1, \ldots, \vartheta_d \in \RR ,
 \]
 where \ $j_0 := 0$ \ and \ $S_0 := 0$.

We derive a generalization of Theorem 1 in Bartkiewicz et al. \cite{BarJakMikWin} on weak convergence of partial sums of strongly stationary jointly regularly varying sequences with index in \ $(0, 2)$.

\begin{Thm}\label{BJMW1}
Assume that \ $(Y_k)_{k\in\NN}$ \ is a strongly stationary sequence satisfying the following conditions:
\renewcommand{\labelenumi}{{\rm(\roman{enumi})}}
 \begin{enumerate}
  \item
   $(Y_k)_{k\in\NN}$ \ is jointly regularly varying with some index \ $\alpha \in (0, 2)$, \ and \ $(a_n)_{n\in\NN}$ \ is a sequence in \ $\RR_{++}$ \ with \ $n \PP(|Y_1| > a_n) \to 1$ \ as \ $n \to \infty$,
  \item
   there exists a sequence \ $(m_n)_{n\in\NN}$ \ in \ $\NN$ \ with \ $m_n \to \infty$ \ and \ $\lfloor n/m_n\rfloor \to \infty$ \ as \ $n \to \infty$ \ such that for each \ $d \in \NN$, \ $\vartheta_1, \ldots, \vartheta_d \in \RR$ \ and \ $t_1, \ldots, t_d \in \RR_+$ \ with \ $0 =: t_0 < t_1 < \ldots < t_d$,
    \[
      \biggl|\varphi_{n,\lfloor nt_1\rfloor,\ldots,\lfloor nt_d\rfloor}(\vartheta_1,\ldots,\vartheta_d)
             - \prod_{\ell=1}^d \bigl(\varphi_{n,m_n}(\vartheta_\ell)\bigr)^{\lfloor(\lfloor nt_\ell\rfloor-\lfloor nt_{\ell-1}\rfloor)/m_n\rfloor}
      \biggr|
      \to 0 \qquad \text{as \ $n \to \infty$,}
    \]
  \item
   for all \ $\vartheta \in \RR$, \ we have
    \[
      \lim_{d\to\infty} \limsup_{n\to\infty} \frac{n}{m_n} \sum_{j=d+1}^{m_n} \EE\bigl(\bigl|\overline{\vartheta a_n^{-1} (S_j - S_d)} \; \cdot\;
                                                                                        \overline{\vartheta a_n^{-1} Y_1}\bigr|\bigr) = 0 ,
    \]
    where, for any random variable \ $Z$, \ we use the notation \ $\overline{Z} := (Z \land 2) \lor (-2)$,
  \item
   there exist the limits
    \[
      \lim_{d\to\infty} (b_+(d) - b_+(d - 1)) =: c_+ , \qquad
      \lim_{d\to\infty} (b_-(d) - b_-(d - 1)) =: c_- ,
    \]
    where
    \[
      b_+(d) := \lim_{n\to\infty} n \PP(S_d > a_n) , \qquad
      b_-(d) := \lim_{n\to\infty} n \PP(S_d \leq - a_n) ,\qquad d\in\NN,
    \]
  \item
   for \ $\alpha \in (1, 2)$ \ assume \ $\EE(Y_1) = 0$, \ and for \ $\alpha = 1$,
    \[
      \lim_{d\to\infty} \limsup_{n\to\infty} n |\EE(\sin(a_n^{-1} S_d))| = 0 .
    \]
 \end{enumerate}
Then \ $c_+, c_- \in \RR_+$ \ and
 \[
   \bigl(a_n^{-1} S_\nt\bigr)_{t\in\RR_+} \distrf (\cS_t)_{t\in\RR_+} \qquad \text{as \ $n \to \infty$,}
 \]
 where \ $(\cS_t)_{t\in\RR_+}$ \ is an \ $\alpha$-stable process such that the characteristic function of \ $\cS_1$ \ has the form
 \[
   \EE(\ee^{\ii\vartheta\cS_1})
   = \begin{cases}
      \exp\bigl\{- C_\alpha |\vartheta|^\alpha \bigl((c_+ + c_-) - \ii(c_+ - c_-) \tan\bigl(\frac{\pi\alpha}{2}\bigr) \sign(\vartheta)\bigr)\bigr\}
       & \text{if \ $\alpha \ne 1$,} \\
      \exp\bigl\{- C_1 |\vartheta| \bigl((c_+ + c_-) + \ii(c_+ - c_-) \frac{2}{\pi} \sign(\vartheta) \log(|\vartheta|)\bigr)\bigr\}
       & \text{if \ $\alpha = 1$,}
     \end{cases}
 \]
 for \ $\vartheta \in \RR$, \ where
 \begin{equation}\label{c_alpha}
   C_\alpha
   := \begin{cases}
       \frac{\Gamma(2 - \alpha)}{1-\alpha}
       \cos\left(\frac{\pi\alpha}{2}\right) ,
        & \text{if \ $\alpha \ne 1$,} \\
       \frac{\pi}{2} , & \text{if \ $\alpha = 1$.}
      \end{cases}
  \end{equation}
\end{Thm}

Next, we formulate a result similar to Lemma 3 in Bartkiewicz et al.\ \cite{BarJakMikWin}, which gives sufficient conditions under which condition
 \textup{(ii)} of Theorem \ref{BJMW1} holds.
First, recall that a strongly stationary sequence \ $(Y_k)_{k\in\NN}$ \ is called strongly mixing with a rate function \ $(\alpha_h)_{h\in\NN}$ \ if its strongly stationary extensions \ $(\ldots,Y_{-2},Y_{-1},Y_0,Y_1,Y_2,\ldots)$ \ admit this property, namely,
 \begin{align}\label{help_mixing_rate}
   \alpha_h := \sup_{A\in\cF^Y_{-\infty,0},\; B\in\cF^Y_{h,\infty}} |\PP(A \cap B) - \PP(A) \PP(B)|
   \to 0 \qquad \text{as \ $h \to \infty$,}
 \end{align}
 where \ $\cF^Y_{-\infty,0} := \sigma(\ldots, Y_{-1}, Y_0)$, \ $\cF^X_{h,\infty} := \sigma(Y_h, Y_{h+1}, \ldots)$, \ $h \in \NN$.
\ For another representation of the rate function \ $(\alpha_h)_{h\in\NN}$, \ see Remark \ref{Rem_strongly_mixing}.
Recall also that a strongly stationary sequence \ $(Y_k)_{k\in\NN}$ \ is strongly mixing with geometric rate if there exists a constant \ $q \in (0, 1)$ \ such that \ $\alpha_h = \OO(q^h)$ \ as \ $h\to\infty$, \ i.e., \ $\sup_{h\in\NN} q^{-h} \alpha_h < \infty$.
\ A strongly mixing, strongly stationary sequence \ $(Y_k)_{k\in\NN}$ \ with a rate function \ $(\alpha_h)_{h\in\NN}$ \ is called
 \ $M_0$-dependent with some \ $M_0\in\ZZ_+$, \ if \ $\alpha_h=0$  \ for all \ $h>M_0$.

\begin{Lem}\label{Lemma3}
Assume that \ $(Y_k)_{k\in\NN}$ \ is a strongly mixing, strongly stationary sequence with a rate function \ $(\alpha_h)_{h\in\NN}$ \ such that
 \ $Y_1$ \ is regularly varying with index \ $\alpha \in (0, 2)$, \ and that \ $(a_n)_{n\in\NN}$ \ is a sequence in \ $\RR_{++}$ \ with \ $n \PP(|Y_1| > a_n) \to 1$ \ as \ $n \to \infty$.
\renewcommand{\labelenumi}{{\rm(\roman{enumi})}}
 \begin{enumerate}
  \item
   If \ $(Y_k)_{k\in\NN}$ \ is \ $M_0$-dependent with some \ $M_0 \in \ZZ_+$, \ or there exists a sequence \ $(\vare_n)_{n\in\NN}$ \ in \ $\RR_{++}$
   \ such that \ $\vare_n \to 0$ \ as \ $n \to \infty$ \ and
    \begin{equation}\label{vare}
     n \alpha_{\lfloor\vare_n a_n^{2\beta}/n\rfloor} \to 0 \qquad \text{as \ $n \to \infty$}
    \end{equation}
    with some \ $\beta \in \bigl(\frac{\alpha}{2}, \alpha \land 1\bigr)$, \ then condition \textup{(ii)} of Theorem \ref{BJMW1} holds with some sequence \ $(m_n)_{n\in\NN}$ \ in \ $\NN$, \ and there exists a sequence \ $(r_n)_{n\in\NN}$ \ in \ $\NN$ \ with \ $r_n \to \infty$ \ and \ $r_n/m_n \to 0$ \ as \ $n \to \infty$, \ such that
    \begin{equation}\label{r}
     n \alpha_{r_n} \to 0 \qquad \text{as \ $n \to \infty$.}
    \end{equation}
  \item
   If \ $(Y_k)_{k\in\NN}$ \ is strongly mixing with geometric rate, then condition \textup{(ii)} of Theorem \ref{BJMW1}, \ $r_n \to \infty$ \ and \ $r_n/m_n \to 0$ \ as \ $n \to \infty$, \ and convergence \eqref{r} hold for \ $m_n = \lfloor n^{\gamma_2}\rfloor$ \ and \ $r_n = \lfloor n^{\gamma_1}\rfloor$, \ $n \in \NN$, \ with arbitrary \ $\gamma_1, \gamma_2 \in (0, 1)$ \ satisfying \ $\gamma_1 < \gamma_2$ \ and \ $\gamma_2 \in \bigl(\frac{1}{2}, \frac{1}{\alpha} \land 1\bigr)$.
 \end{enumerate}
\end{Lem}

\begin{Rem}\label{Rem_Lemma3}
(i). \ Note that condition \eqref{vare} is equivalent to \ $n \alpha_{\lfloor\vare_n ( \frac{a_n^{2\beta}}{n} \wedge n) \rfloor} \to 0$ \ as \ $n \to \infty$,
 \ since \ $\beta \in \bigl(\frac{\alpha}{2}, \alpha \land 1\bigr)$ \ yields that
 \begin{equation}\label{a_n_beta_n}
   \frac{a_n^{2\beta}}{n^2} = \frac{ (n^{\frac{1}{\alpha}}L(n))^{2\beta} }{n^2} = n^{\frac{2\beta}{\alpha} - 2} (L(n))^{2\beta} \to 0
    \qquad \text{as \ $n\to\infty$,}
 \end{equation}
 where \ $L:\RR_{++}\to\RR_{++}$ \ is some slowly varying continuous function,
 and hence \ $\frac{a_n^{2\beta}}{n} \wedge n = \frac{a_n^{2\beta}}{n}$ \ for large enough \ $n\in\NN$.

(ii). \ If \ $(Y_k)_{k\in\NN}$ \ is strongly mixing with geometric rate, then \eqref{vare} holds for \ $\vare_n = n^{-c}$, \ $n \in \NN$, \ with arbitrary \ $c \in \bigl(0, \bigl(\frac{2}{\alpha} - 1\bigr) \land 1\bigr)$.
\ For more details, see Section \ref{Proofs}.
\proofend
\end{Rem}

Next, we present a sufficient condition under which condition (iii) of Theorem \ref{BJMW1} holds.
This sufficient condition is nothing else but the anti-clustering type condition (8.3.13) in the forthcoming book of Mikosch and Wintenberger \cite{MikWin2}.

\begin{Lem}\label{Lemma2}
Assume that \ $(Y_k)_{k\in\NN}$ \ is a strongly stationary sequence such that conditions (i) and (ii) of Theorem \ref{BJMW1} hold and
 \begin{equation}\label{3.9}
  \lim_{d\to\infty} \limsup_{n\to\infty} n \sum_{j=d+1}^{m_n} \EE\big( (\vert a_n^{-1} Y_j\vert \wedge x) (\vert a_n^{-1} Y_0\vert \wedge x) \big)= 0
        \qquad \text{for all \ $x\in\RR_{++}$,}
 \end{equation}
 where \ $(a_n)_{n\in\NN}$ \ and \ $(m_n)_{n\in\NN}$ \ are the same as in conditions \textup{(i)} and \textup{(ii)} of Theorem \ref{BJMW1}, respectively.
Then condition \textup{(iii)} of Theorem \ref{BJMW1} holds.
\end{Lem}

\section{An application on aggregation of branching processes}
\label{application}

Let \ $(X_k)_{k\in\ZZ_+}$ \ be a Galton--Watson branching process with immigration.
For each \ $k, j \in \ZZ_+$, \ the number of individuals in the \ $k^\mathrm{th}$
 \ generation will be denoted by \ $X_k$, \ the number of offsprings produced by
 the \ $j^\mathrm{th}$ \ individual belonging to the \ $(k-1)^\mathrm{th}$
 \ generation will be denoted by \ $\xi_{k,j}$, \ and the number of immigrants in the
 \ $k^\mathrm{th}$ \ generation will be denoted by \ $\vare_k$.
\ Then we have
 \[
  X_k = \sum_{j=1}^{X_{k-1}} \xi_{k,j} + \vare_k , \qquad k \in \NN ,
 \]
 where we define \ $\sum_{j=1}^0 := 0$.
\ Here \ $\bigl\{X_0, \, \xi_{k,j}, \, \vare_k : k, j \in \NN\bigr\}$ \ are supposed
 to be independent non-negative integer-valued random variables.
Moreover, \ $\{\xi_{k,j} : k, j \in \NN\}$ \ and
 \ $\{\vare_k : k \in \NN\}$ \ are supposed to consist of identically
 distributed random variables, respectively.
For notational convenience, let \ $\xi$ \ and \ $\vare$ \ be random variables such that \ $\xi \distre \xi_{1,1}$ \ and \ $\vare \distre \vare_1$.

If \ $m_\xi := \EE(\xi) \in [0, 1)$ \ and \ $\sum_{\ell=1}^\infty \log(\ell) \PP(\vare = \ell) < \infty$, \ then the Markov chain \ $(X_k)_{k\in\ZZ_+}$ \ admits a unique stationary distribution \ $\pi$, \ see, e.g., Quine \cite{Qui} (for more details, see the proof of Lemma \ref{lemma:strong_mixing}).
Note that if \ $m_\xi \in [0, 1)$ \ and \ $\PP(\vare = 0) = 1$, \ then \ $\sum_{\ell=1}^\infty \log(\ell) \PP(\vare = \ell) = 0$ \ and \ $\pi$ \ is the Dirac measure \ $\delta_0$ \ concentrated at the point \ $0$.
\ In fact, \ $\pi = \delta_0$ \ if and only if \ $\PP(\vare = 0) = 1$.
\ Moreover, if \ $m_\xi = 0$ \ (which is equivalent to \ $\PP(\xi = 0) = 1$), \ then \ $\pi$ \ is the distribution of \ $\vare$.

In what follows, we formulate our assumptions valid for the whole section.
We assume that
 \begin{itemize}
  \item $m_\xi \in [0, 1)$ \ (so-called subcritical case),
  \item $\vare$ \ is regularly varying with index \ $\alpha \in (0, 1)\cup (1,\frac{4}{3})$, \ i.e.,
        \[
         \lim_{x\to\infty} \frac{\PP(\vare > qx)}{\PP(\vare > x)} = q^{-\alpha} \qquad
         \text{for all \ $q \in \RR_{++}$,}
        \]
        which yields that \ $\PP(\vare = 0) < 1$ \ and \ $\sum_{\ell=1}^\infty \log(\ell) \PP(\vare = \ell) < \infty$ \ (see, e.g.,
        Barczy et al.\ \cite[Lemma A.3]{BarBosPap}), and hence the Markov process \ $(X_k)_{k\in\ZZ_+}$ \ admits a unique stationary distribution \ $\pi$,
  \item $X_0 \distre \pi$, \ yielding that the Markov chain \ $(X_k)_{k\in\ZZ_+}$ \ is strongly stationary,
  \item $\EE(\xi^2) < \infty$ \ in case of \ $\alpha \in (1, \frac{4}{3})$.
 \end{itemize}
By Basrak et al. \cite[Theorem 2.1.1]{BasKulPal} (see also Theorem \ref{Xtail}), \ $X_0$ \ is regularly varying with index
 \ $\alpha$, \ yielding the existence of a sequence \ $(a_n)_{n\in\NN}$ \ in \ $\RR_{++}$ \ with
 \ $n \PP(X_0 > a_n) \to 1$ \ as \ $n \to \infty$, \ see, e.g., Lemma \ref{a_n}.
Let us fix an arbitrary sequence \ $(a_n)_{n\in\NN}$ \ in \ $\RR_{++}$ \ with this property.
Let \ $X^{(j)} = (X^{(j)}_k)_{k\in\ZZ_+}$, \ $j \in \NN$, \ be a sequence of
 independent copies of \ $(X_k)_{k\in\ZZ_+}$.

For a temporal aggregation, we have the following theorem considering different centerings.

\begin{Thm}\label{aggr_time}
\renewcommand{\labelenumi}{{\rm(\roman{enumi})}}
 \begin{enumerate}
  \item
   For each \ $\alpha \in (0, 1) \cup (1, \frac{4}{3})$,
 \begin{align*}
  \biggl(\frac{1}{a_n}
          \sum_{k=1}^\nt
           \Bigl(X_k
                 - \EE\bigl(X_k
                            \bbone_{\{X_k\leq a_n\}}\bigr)\Bigr)
   \biggr)_{t\in\RR_+}
  &= \biggl(\frac{1}{a_n} \sum_{k=1}^\nt X_k
            - \frac{\nt}{a_n} \EE\bigl(X_0
                            \bbone_{\{X_0\leq a_n\}}\bigr)
   \biggr)_{t\in\RR_+} \\
  &\distrf \bigl(\cZ_t^{(\alpha)}\bigr)_{t\in\RR_+} \qquad \text{as \ $n \to \infty$,}
 \end{align*}
  \item
   in case of \ $\alpha \in (0, 1)$,
 \[
   \biggl(\frac{1}{a_n} \sum_{k=1}^\nt X_k\biggr)_{t\in\RR_+}
   \distrf \Bigl(\cZ_t^{(\alpha)} + \frac{\alpha}{1-\alpha} t\Bigr)_{t\in\RR_+} \qquad \text{as \ $n \to \infty$,}
 \]
  \item
   in case of \ $\alpha \in (1, \frac{4}{3})$,
 \begin{align*}
  \biggl(\frac{1}{a_n} \sum_{k=1}^\nt (X_k - \EE(X_k))\biggr)_{t\in\RR_+}
  &= \biggl(\frac{1}{a_n} \sum_{k=1}^\nt X_k - \frac{\nt}{a_n} \EE(X_0) \biggr)_{t\in\RR_+} \\
  &\distrf \Bigl(\cZ_t^{(\alpha)} + \frac{\alpha}{1-\alpha} t\Bigr)_{t\in\RR_+} \qquad \text{as \ $n \to \infty$,}
 \end{align*}
 \end{enumerate}
\noindent
 where \ $\bigl(\cZ_t^{(\alpha)}\bigr)_{t\in\RR_+}$ \ is an \ $\alpha$-stable process such that the characteristic function of the distribution of \ $\cZ_1^{(\alpha)}$ \ has the form
 \[
   \EE\bigl(\ee^{\ii\vartheta\cZ_1^{(\alpha)}}\bigr)
   = \exp\biggl\{\ii b_\alpha \vartheta + \frac{1-m_\xi^\alpha}{(1-m_\xi)^\alpha}
                 \int_0^\infty
                  (\ee^{\ii\vartheta u} - 1
                   - \ii \vartheta u \bbone_{(0,1]}(u)) \alpha
                  u^{-1-\alpha} \, \dd u\biggr\} , \qquad \vartheta \in \RR ,
 \]
 where
 \[
   b_\alpha := \biggl(\frac{1-m_\xi^\alpha}{(1-m_\xi)^\alpha} - 1\biggr) \frac{\alpha}{1-\alpha} , \qquad \alpha \in (0, 1) \cup \left(1, \frac{4}{3}\right) ,
 \]
 and \ $\bigl(\cZ_t^{(\alpha)} + \frac{\alpha}{1-\alpha} t\bigr)_{t\in\RR_+}$ \ is an \ $\alpha$-stable process such that the characteristic function of the distribution of \ $\cZ_1^{(\alpha)} + \frac{\alpha}{1-\alpha}$ \ has the form
 \begin{align*}
  &\EE\Bigl(\exp\Bigl\{\ii\vartheta\Bigl(\cZ_1^{(\alpha)}+\frac{\alpha}{1-\alpha}\Bigr)\Bigr\}\Bigr)
   = \begin{cases}
      \exp\Bigl\{\frac{1-m_\xi^\alpha}{(1-m_\xi)^\alpha}
                 \int_0^\infty
                  (\ee^{\ii\vartheta u} - 1) \alpha
                  u^{-1-\alpha} \, \dd u \Bigr\} ,
       & \text{if \ $\alpha \in (0, 1)$,} \\[3mm]
      \exp\Bigl\{\frac{1-m_\xi^\alpha}{(1-m_\xi)^\alpha}
                 \int_0^\infty
                  (\ee^{\ii\vartheta u} - 1
                   - \ii \vartheta u) \alpha
                  u^{-1-\alpha} \, \dd u \Bigr\} ,
       & \text{if \ $\alpha \in (1, \frac{4}{3})$,}
     \end{cases} \\
  &= \exp\biggl\{- C_\alpha \frac{1-m_\xi^\alpha}{(1-m_\xi)^\alpha}
                   |\vartheta|^\alpha
                    \left(1 - \ii \tan\left(\frac{\pi\alpha}{2}\right)
                                  \sign(\vartheta)\right)\biggr\} \quad \text{if \ $\alpha \in (0, 1) \cup \left(1, \frac{4}{3}\right)$,}
 \end{align*}
 for \ $\vartheta \in \RR$.
\end{Thm}

In the next remark we explain why we do not have any result in case of \ $\alpha=1$ \ and \ $\alpha\in[\frac{4}{3},2)$.

\begin{Rem}\label{Rem_support}
(i). The proof of Theorem \ref{aggr_time} is based on Theorem \ref{BJMW1}, and Theorem \ref{BJMW1} can not be applied in case of \ $\alpha = 1$ \ for \ $(X_n)_{n\in\ZZ_+}$.
\ This is the reason for not having any result in case of \ $\alpha=1$.
\ Indeed, if all the conditions of Theorem \ref{BJMW1} were satisfied, then the limit 1-stable process \ $(\cS_t)_{t\in\RR_+}$ \ would be nonnegative.
Clearly, \ $c_- = 0$, \ thus, by (14.20) in Sato \cite{Sato}, the L\'evy measure of \ $\cS_1$ \ has the form \ $c_+ u^{-2} \bbone_{\RR_{++}}(u) \, \dd u$, \ hence, by Sato \cite[Definition 11.9 and part (i) of Theorem 24.10]{Sato}, the support of \ $\cS_1$ \ is \ $\RR$, \ leading to a contradiction.

\noindent (ii).  Our present technique is not suitable for handling the case \ $\alpha\in[\frac{4}{3},2)$.
\ The reason for it is that we were not able to check the corresponding anti-clustering type condition (iii) of Theorem \ref{BJMW1}
 nor the condition \eqref{3.9} which would imply condition (iii) of Theorem \ref{BJMW1}.
To give more technical details, in Step 2/(b) of the proof of Theorem \ref{aggr_time}, we need \ $\frac{nm_n}{a_n^2}\to 0$ \ as \ $n\to\infty$,
 \ where, by part (ii) of Lemma \ref{Lemma3}, \ $m_n=\lfloor n^{\gamma_2}\rfloor$ \ with some
 \ $\gamma_2 \in \bigl(\frac{1}{2}, \frac{1}{\alpha} \land 1\bigr)$.
\ So \ $\frac{nm_n}{a_n^2}$ \ is asymptotically equivalent to \ $n^{1+\gamma_2 - \frac{2}{\alpha}} (L(n))^{-2}$ \ as \ $n\to\infty$ \
 with some slowly varying function \ $L$, \ and hence we need \ $\gamma_2<\frac{2}{\alpha}-1$. \
In case of \ $\alpha\in(1,2)$, \ it can be fulfilled provided that
 \ $\frac{1}{2} < \frac{2}{\alpha}-1 < \frac{1}{\alpha} \Longleftrightarrow \alpha\in(1,\frac{4}{3})$.
\ We note that the condition  \ $\frac{nm_n}{a_n^2}\to 0$ \ as \ $n\to\infty$ \ also appears in the forthcoming book
 of Mikosch and Wintenberger \cite[page 324]{MikWin2}.
\proofend
\end{Rem}

\begin{Rem}
Note that, in accordance with Basrak and Segers \cite[Remark 4.8]{BasSeg} and Mikosch and Wintenberger \cite[page 171]{MikWin},
 in case of \ $\alpha \in (0, 1)$, \ we have
 \begin{align}\label{help_kar_theta}
 \begin{split}
 &\EE\Bigl(\exp\Bigl\{\ii \vartheta \Bigl(\cZ_1^{(\alpha)} + \frac{\alpha}{1-\alpha}\Bigr)\Bigr\}\Bigr) \\
 &= \exp\left\{ - \int_0^\infty \EE\left[ \exp\Big( \ii u \vartheta \sum_{\ell=1}^\infty \Theta_\ell \Big)
                               - \exp\Big( \ii u \vartheta \sum_{\ell=0}^\infty \Theta_\ell \Big) \right] \alpha u^{-\alpha-1}\,\dd u \right\}
 \end{split}
 \end{align}
 for \ $\vartheta \in \RR$, \ where \ $(\Theta_\ell)_{\ell\in\ZZ_+}$ \ is the (forward) spectral tail process of \ $(X_\ell)_{\ell\in\ZZ_+}$, \
  see Barczy et al.\ \cite[Remark 2]{BarNedPap3}.
We also remark that, using (14.19) in Sato \cite{Sato}, one can check that \eqref{help_kar_theta} does not hold in case of \ $\alpha \in (1, \frac{4}{3})$, \ which is somewhat unexpected in view of page 171 in Mikosch and Wintenberger \cite{MikWin}.
In fact, the above mentioned mistake on page 171 in Mikosch and Wintenberger \cite{MikWin} was already corrected
 in Mikosch and Wintenberger \cite[Theorem 4.3]{MikWin3}.
Namely, in case of \ $\alpha \in (1, \frac{4}{3})$, \ we have
 \begin{align*}
 &\EE\Bigl(\exp\Bigl\{\ii \vartheta \Bigl(\cZ_1^{(\alpha)} + \frac{\alpha}{1-\alpha}\Bigr)\Bigr\}\Bigr) \\
 &= \exp\left\{ - \int_0^\infty \EE\left[ \exp\Big( \ii u \vartheta \sum_{\ell=1}^\infty \Theta_\ell \Big)
                               - \exp\Big( \ii u \vartheta \sum_{\ell=0}^\infty \Theta_\ell \Big) + \ii u \vartheta \Theta_0\right] \alpha u^{-\alpha-1}\,\dd u \right\}
 \end{align*}
 for \ $\vartheta \in \RR$.
\ Indeed, by using (14.19) in Sato \cite{Sato}, we have
 \begin{align*}
   \int_0^\infty
    (\ee^{\ii\vartheta u} - 1 - \ii \vartheta u) u^{-1-\alpha}
    \, \dd u
   = - \frac{C_\alpha}{\alpha} |\vartheta|^\alpha
     \left(1 - \ii \tan\left(\frac{\pi\alpha}{2}\right) \sign(\vartheta)\right) , \qquad \vartheta\in\RR,
 \end{align*}
 and hence using that \ $\Theta_\ell= m_\xi^\ell$, \ $\ell\in\ZZ_+$ \ (see, e.g., (3.7) and (3.8) in Barczy et al.\ \cite{BarNedPap3}),
 we get
 \begin{align*}
  &\exp\left\{ - \int_0^\infty  \EE\left[ \exp\Big( \ii u \vartheta \sum_{\ell=1}^\infty \Theta_\ell \Big)
                               - \exp\Big( \ii u \vartheta \sum_{\ell=0}^\infty \Theta_\ell \Big) + \ii u \vartheta \Theta_0 \right] \alpha u^{-\alpha-1}\,\dd u \right\}  \\
  & = \exp\left\{ - \int_0^\infty  \EE\left[ \exp\Big( \ii u \vartheta \sum_{\ell=1}^\infty m_\xi^\ell \Big)
                               - \exp\Big( \ii u \vartheta \sum_{\ell=0}^\infty m_\xi^\ell \Big) + \ii u \vartheta \right] \alpha u^{-\alpha-1}\,\dd u \right\} \\
  & = \exp\left\{ - \int_0^\infty  \left( \exp\Big( \ii u \vartheta \frac{m_\xi}{1-m_\xi} \Big)
                               - \exp\Big( \ii u \vartheta \frac{1}{1-m_\xi} \Big) + \ii u \vartheta  \right) \alpha u^{-\alpha-1}\,\dd u \right\} \\
  &= \exp\left\{ - \int_0^\infty  \left( \exp\Big( \ii u \frac {\vartheta m_\xi}{1-m_\xi} \Big) - 1
                 - \ii u \frac{\vartheta m_\xi}{1-m_\xi} \right) \alpha u^{-\alpha-1}\,\dd u \right.\\
  &\phantom{= \exp\Big\{\,}
         \left. +  \int_0^\infty  \left( \exp\Big( \ii u \frac{\vartheta}{1-m_\xi} \Big) - 1 -  \ii u \frac{\vartheta}{1-m_\xi} \right) \alpha u^{-\alpha-1}\,\dd u \right\} \\
  & = \exp\Bigg\{ C_\alpha \left\vert \frac{\vartheta m_\xi}{1-m_\xi} \right\vert^\alpha
                   \left(1 - \ii \tan\left(\frac{\pi\alpha}{2}\right)
                                  \sign\left(\frac{\vartheta m_\xi}{1-m_\xi}\right)\right)  \\
  &\phantom{= \exp\Bigg\{ \;}
                   - C_\alpha \left\vert \frac{\vartheta }{1-m_\xi} \right\vert^\alpha
                   \left(1 - \ii \tan\left(\frac{\pi\alpha}{2}\right)
                                  \sign\left(\frac{\vartheta}{1-m_\xi}\right)\right) \Bigg\} \\
  & = \exp\left\{ - C_\alpha \frac{1-m_\xi^\alpha}{(1-m_\xi)^\alpha} \vert\vartheta\vert^\alpha
                        \left(1 - \ii \tan\left(\frac{\pi\alpha}{2}\right)
                                  \sign(\vartheta)\right) \right\} ,
 \end{align*}
 as desired.
\proofend
\end{Rem}

\begin{Rem}\label{cZ_properties}
If \ $\alpha \in (0, 1)$, \ then the drift of the distribution of \ $\cZ_1^{(\alpha)} + \frac{\alpha}{1-\alpha}$ \ is \ $0$, \ hence the process \ $\bigl(\cZ_t^{(\alpha)} + \frac{\alpha}{1-\alpha} t\bigr)_{t\in\RR_+}$ \ is strictly \ $\alpha$-stable, see, e.g., Sato \cite[Theorem 14.7 (iv) and Definition 13.2]{Sato}.

If \ $\alpha \in (1, \frac{4}{3})$, \ then the center, i.e., the expectation of
 \ $\cZ_1^{(\alpha)} + \frac{\alpha}{1-\alpha}$ \ is $0$, hence the process
 $\bigl(\cZ_t^{(\alpha)} + \frac{\alpha}{1-\alpha} t\bigr)_{t\in\RR_+}$ is strictly $\alpha$-stable see, e.g., Sato \cite[Theorem 14.7 (vi) and Definition 13.2]{Sato}.

All in all, the process \ $\bigl(\cZ_t^{(\alpha)} + \frac{\alpha}{1-\alpha}t\bigr)_{t\in\RR_+}$ \ is strictly \ $\alpha$-stable for any
 \ $\alpha \in (0,1)\cup(1,\frac{4}{3})$.
\proofend
\end{Rem}

In the end, limit theorems will be presented for finite dimensional distributions of the aggregated stochastic process
 \ $\bigl(\sum_{j=1}^N \sum_{k=1}^{\lfloor nt \rfloor} X^{(j)}_k\bigr)_{t\in\RR_+}$ \
 considering different centerings as first \ $n$ \ and then \ $N$ \ converges to infinity.
The other iterated case of taking the limit first as \ $N$ \ and then as \ $n$ \ converges to infinity has been studied in Barczy et al. \cite[Theorem 2]{BarNedPap3}.

\begin{Thm}\label{iterated_aggr_2}
\renewcommand{\labelenumi}{{\rm(\roman{enumi})}}
 \begin{enumerate}
  \item For each \ $\alpha \in (0, 1)$,
 \begin{equation}\label{iterated_aggr_2_1}
 \begin{aligned}
  &\cD_\ff\text{-}\hspace*{-1mm}\lim_{N\to\infty} \,
   \cD_\ff\text{-}\hspace*{-1mm}\lim_{n\to\infty} \,
    \biggl(\frac{1}{a_nN^{\frac{1}{\alpha}}}
           \sum_{k=1}^\nt \sum_{j=1}^N
            \Bigl(X^{(j)}_k
                  - \EE\bigl(X^{(j)}_k
                             \bbone_{\{X^{(j)}_k\leq a_n\}}\bigr)
            \Bigr)
   \biggr)_{t\in\RR_+} \\
  &=\cD_\ff\text{-}\hspace*{-1mm}\lim_{N\to\infty} \,
    \cD_\ff\text{-}\hspace*{-1mm}\lim_{n\to\infty} \,
    \biggl(\frac{1}{a_nN^{\frac{1}{\alpha}}}
           \sum_{k=1}^\nt \sum_{j=1}^N X^{(j)}_k
           - \frac{\nt N}{a_nN^{\frac{1}{\alpha}}}
             \EE\bigl(X_0 \bbone_{\{X_0\leq a_n\}}\bigr)\biggr)_{t\in\RR_+} \\
  &= \Bigl(\cZ_t^{(\alpha)} + \frac{\alpha}{1-\alpha} t\Bigr)_{t\in\RR_+} ,
 \end{aligned}
 \end{equation}
 and
  \begin{equation}\label{iterated_aggr_2_2}
   \cD_\ff\text{-}\hspace*{-1mm}\lim_{N\to\infty} \,
   \cD_\ff\text{-}\hspace*{-1mm}\lim_{n\to\infty} \,
    \biggl(\frac{1}{a_nN^{\frac{1}{\alpha}}}
           \sum_{k=1}^\nt \sum_{j=1}^N
            X^{(j)}_k
            \Bigr)
   \biggr)_{t\in\RR_+}
   = \Bigl(\cZ_t^{(\alpha)} + \frac{\alpha}{1-\alpha} t\Bigr)_{t\in\RR_+} .
 \end{equation}
\item For each \ $\alpha\in(1,\frac{4}{3})$,
 \begin{equation}\label{iterated_aggr_2_4}
 \begin{aligned}
  &\cD_\ff\text{-}\hspace*{-1mm}\lim_{N\to\infty} \,
   \cD_\ff\text{-}\hspace*{-1mm}\lim_{n\to\infty} \,
    \biggl(\frac{1}{a_nN^{\frac{1}{\alpha}}}
           \sum_{k=1}^\nt \sum_{j=1}^N
            (X^{(j)}_k - \EE(X^{(j)}_k))
            \Bigr)
   \biggr)_{t\in\RR_+} \\
  &=\cD_\ff\text{-}\hspace*{-1mm}\lim_{N\to\infty} \,
    \cD_\ff\text{-}\hspace*{-1mm}\lim_{n\to\infty} \,
    \biggl(\frac{1}{a_nN^{\frac{1}{\alpha}}}
           \sum_{k=1}^\nt \sum_{j=1}^N X^{(j)}_k
           - \frac{\nt N}{a_nN^{\frac{1}{\alpha}}} \EE(X_0)
   \biggr)_{t\in\RR_+} \\
  &= \Bigl(\cZ_t^{(\alpha)} + \frac{\alpha}{1-\alpha} t\Bigr)_{t\in\RR_+} .
 \end{aligned}
 \end{equation}
 \end{enumerate}
\end{Thm}

\section{Proofs}
\label{Proofs}

\noindent{\bf Proof of Theorem \ref{BJMW1}.}
By the continuous mapping theorem, it is enough to check that for each \ $d \in \NN$ \ and \ $t_1, \ldots, t_d \in \RR_{++}$ \ with \ $t_1 < \ldots < t_d$, \ we have
 \[
   a_n^{-1} \bigl(S_{\lfloor nt_1\rfloor}, S_{\lfloor nt_2\rfloor} - S_{\lfloor nt_1\rfloor}, \ldots, S_{\lfloor nt_d\rfloor} - S_{\lfloor nt_{d-1}\rfloor}\bigr)
   \distr (\cS_{t_1}, \cS_{t_2} - \cS_{t_1}, \ldots, \cS_{t_d} - \cS_{t_{d-1}})
 \]
 as \ $n \to \infty$.
\ In view of condition (ii), taking also into account that the process \ $(\cS_t)_{t\in\RR_+}$ \ has independent and stationary increments
 (being a L\'evy process), this is proved if we can show that for each \ $\ell \in \{1, \ldots, d\}$, \ we have
 \begin{equation}\label{BJMW2}
  a_n^{-1} \sum_{j=1}^{\lfloor(\lfloor nt_\ell\rfloor-\lfloor nt_{\ell-1}\rfloor)/m_n\rfloor} S_{m_n}^{(j)} \distr \cS_{t_\ell} - \cS_{t_{\ell-1}} \distre \cS_{t_\ell-t_{\ell-1}} \qquad \text{as \ $n \to \infty$,}
 \end{equation}
 where \ $(S_n^{(j)})_{n\in\NN}$, \ $j \in \NN$, \ are independent copies of \ $(S_n)_{n\in\NN}$.
\ Since \ $(Y_k)_{k\in\NN}$ \ is strongly stationary, by the condition (ii) with \ $d=1$ \ and \ $t_1=\frac{\lfloor nt_\ell\rfloor - \lfloor nt_{\ell-1}\rfloor}{n}$,
 \  \eqref{BJMW2} is proved if we can show that for each \ $\ell \in \{1, \ldots, d\}$,
 \[
  a_n^{-1} \bigl(S_{\lfloor nt_\ell\rfloor} - S_{\lfloor nt_{\ell-1}\rfloor}\bigr)
  \distre a_n^{-1} S_{\lfloor nt_\ell\rfloor-\lfloor nt_{\ell-1}\rfloor}
  \distr \cS_{t_\ell-t_{\ell-1}} \qquad \text{as \ $n \to \infty$.}
 \]
Condition (ii) with \ $d = 1$ \ and \ $t_1 = 1$ \ implies the mixing condition (2.8) of Theorem 1 in Bartkiewicz et al.\ \cite{BarJakMikWin}, hence we may apply Theorem 1 in Bartkiewicz et al.\ \cite{BarJakMikWin}, thus we obtain \ $a_n^{-1} S_n \distr \cS_1$ \ as \ $n \to \infty$, \ and we get
 \begin{equation}\label{BJMW4}
  a_{\lfloor nt_\ell\rfloor-\lfloor nt_{\ell-1}\rfloor}^{-1} S_{\lfloor nt_\ell\rfloor-\lfloor nt_{\ell-1}\rfloor} \distr \cS_1 \qquad \text{as \ $n \to \infty$,}
 \end{equation}
 since \ $(\lfloor nt_\ell\rfloor - \lfloor nt_{\ell-1}\rfloor)/n \to t_\ell - t_{\ell-1} \in \RR_{++}$ \ as \ $n \to \infty$ \ yields \ $\lfloor nt_\ell\rfloor-\lfloor nt_{\ell-1}\rfloor \to \infty$ \ as \ $n \to \infty$.
\ There exists a slowly varying continuous function \ $L : \RR_{++} \to \RR_{++}$ \ such that \ $a_n = n^{1/\alpha} L(n)$, \ $n \in \NN$,
\ thus, by the uniform convergence theorem for slowly varying functions (see, e.g., Bingham et al.\ \cite[Theorem 1.5.2]{BinGolTeu}),
 \begin{equation}\label{BJMW5}
   \frac{a_{\lfloor nt_\ell\rfloor-\lfloor nt_{\ell-1}\rfloor}}{a_n}
   = \frac{(\lfloor nt_\ell\rfloor-\lfloor nt_{\ell-1}\rfloor)^{1/\alpha}}{n^{1/\alpha}}
     \frac{L(\lfloor nt_\ell\rfloor-\lfloor nt_{\ell-1}\rfloor)}{L(n)}
   \to (t_\ell - t_{\ell-1})^{1/\alpha}
 \end{equation}
 as \ $n \to \infty$ \ for every \ $\ell \in \{1, \ldots d\}$, \ since \ $(\lfloor nt_\ell\rfloor - \lfloor nt_{\ell-1}\rfloor)/n \in [(t_\ell -t_{\ell-1})/2, t_\ell]$ \ for sufficiently large \ $n \in \NN$.
\ By \eqref{BJMW4}, \eqref{BJMW5} and Slutsky's lemma, we obtain
 \[
   a_n^{-1} S_{\lfloor nt_\ell\rfloor-\lfloor nt_{\ell-1}\rfloor}
   \distr (t_\ell - t_{\ell-1})^{1/\alpha} \cS_1
   \distre \cS_{t_\ell-t_{\ell-1}} \qquad \text{as \ $n \to \infty$,}
 \]
 hence we showed \eqref{BJMW2}, as desired.
\proofend

\noindent{\bf Proof of Lemma \ref{Lemma3}.}
(i). \ We start by showing that condition (ii) of Theorem \ref{BJMW1} holds for a suitable sequence \ $(m_n)_{n\in\NN}$.
\ First, we derive sufficient conditions on a sequence \ $(m_n)_{n\in\NN}$ \ under which condition (ii) of Theorem \ref{BJMW1} holds, and then we show that
 these sufficient conditions can be fulfilled.
Let us fix \ $d \in \NN$, \ $\vartheta_1, \ldots, \vartheta_d \in \RR$ \ and \ $t_1, \ldots, t_d \in \RR_+$ \ with \ $0 =: t_0 < t_1 < \ldots < t_d$.
\ Let us introduce the notations \ $k_{n,\ell} := \lfloor(\lfloor n t_\ell\rfloor - \lfloor n t_{\ell-1}\rfloor)/m_n\rfloor$ \ for \ $n \in \NN$, \ $\ell \in \{1, \ldots, d\}$.
\ By the inequalities
 \[
  \biggl|\prod_{\ell=1}^d z_\ell - \prod_{\ell=1}^d z_\ell'\biggr|
  \leq \sum_{\ell=1}^d |z_\ell - z_\ell'|
 \]
 for \ $z_1, z_1', \ldots, z_d, z_d' \in D_1 := \{z \in \CC : |z| \leq 1\}$, \ $|\ee^{\ii u} - 1| \leq |u| \land 2$ \ for \ $u \in \RR$, \
 $|u + v|^\gamma \leq |u|^\gamma + |v|^\gamma$ \ for \ $u, v \in \RR$ \ and \ $\gamma \in (0, 1]$, \ and \ $y\leq y^\gamma$ \ for \ $y\in[0,1]$ \ and \ $\gamma \in (0, 1]$, \ and using the strong stationarity of \ $(Y_k)_{k\in\NN}$, \ we have
 \begin{align*}
  &\biggl|\varphi_{n,\lfloor nt_1\rfloor,\ldots,\lfloor nt_d\rfloor}(\vartheta_1,\ldots,\vartheta_d)
         - \EE\biggl(\exp\biggl\{\frac{\ii}{a_n} \sum_{\ell=1}^d \vartheta_\ell \bigl(S_{\lfloor n t_{\ell-1}\rfloor+k_{n,\ell}m_n}
         - S_{\lfloor n t_{\ell-1}\rfloor}\bigr)\biggr\}\biggr)\biggr| \\
  &= \biggl|\EE\biggl(\exp\biggl\{\frac{\ii}{a_n} \sum_{\ell=1}^d \vartheta_\ell \bigl(S_{\lfloor n t_{\ell-1}\rfloor+k_{n,\ell}m_n}
               - S_{\lfloor n t_{\ell-1}\rfloor}\bigr) \biggr\} \\
  &\phantom{= \biggl|\EE\biggl(}
                      \times
                      \biggl(\exp\biggl\{\frac{\ii}{a_n} \sum_{\ell=1}^d \vartheta_\ell \bigl(S_{\lfloor n t_\ell\rfloor} - S_{\lfloor n t_{\ell-1}\rfloor+k_{n,\ell}m_n}\bigr)\biggr\} - 1\biggr)\biggr)\biggr|\leq
 \end{align*}
 \begin{align*}
  &\leq \EE\biggl(\biggl|\exp\biggl\{\frac{\ii}{a_n} \sum_{\ell=1}^d \vartheta_\ell \bigl(S_{\lfloor n t_\ell\rfloor} - S_{\lfloor n t_{\ell-1}\rfloor+k_{n,\ell}m_n}\bigr)\biggr\} - 1\biggr|\biggr)\\
  &= \EE\biggl(\biggl|\prod_{\ell=1}^d\exp\biggl\{\frac{\ii\vartheta_\ell}{a_n} \bigl(S_{\lfloor n t_\ell\rfloor} - S_{\lfloor n t_{\ell-1}\rfloor+k_{n,\ell}m_n}\bigr)\biggr\} - 1\biggr|\biggr) \\
  &\leq \EE\biggl(\sum_{\ell=1}^d \biggl|\exp\biggl\{\frac{\ii\vartheta_\ell}{a_n} \bigl(S_{\lfloor n t_\ell\rfloor} - S_{\lfloor n t_{\ell-1}\rfloor+k_{n,\ell}m_n}\bigr)\biggr\} - 1\biggr|\biggr) \\
  &= \sum_{\ell=1}^d \EE\biggl(\biggl|\exp\biggl\{\frac{\ii\vartheta_\ell}{a_n} S_{\lfloor n t_\ell\rfloor-\lfloor n t_{\ell-1}\rfloor-k_{n,\ell}m_n}\biggr\} - 1\biggr|\biggr)\\
  &\leq \sum_{\ell=1}^d \EE\biggl(\biggl|\frac{\vartheta_\ell}{a_n} S_{\lfloor n t_\ell\rfloor-\lfloor n t_{\ell-1}\rfloor-k_{n,\ell}m_n}\biggr| \land 2\biggr)
   \leq 2 \sum_{\ell=1}^d \EE\Biggl(\Biggl(\frac{|\vartheta_\ell|}{2a_n} \sum_{i=1}^{\lfloor n t_\ell\rfloor-\lfloor n t_{\ell-1}\rfloor-k_{n,\ell}m_n} |Y_i|\Biggr) \land 1\Biggr) \\
  &\leq 2 \sum_{\ell=1}^d \EE\Biggl(\Biggl(\frac{|\vartheta_\ell|}{2a_n} \sum_{i=1}^{\lfloor n t_\ell\rfloor-\lfloor n t_{\ell-1}\rfloor-k_{n,\ell}m_n} |Y_i|\Biggr)^\beta \land 1\Biggr)
   \leq 2 \sum_{\ell=1}^d \EE\Biggl(\frac{|\vartheta_\ell|^\beta}{(2a_n)^\beta} \sum_{i=1}^{\lfloor n t_\ell\rfloor-\lfloor n t_{\ell-1}\rfloor-k_{n,\ell}m_n} |Y_i|^\beta\Biggr) \\
  &= \frac{2}{(2a_n)^\beta} \sum_{\ell=1}^d |\vartheta_\ell|^\beta \bigl(\lfloor n t_\ell\rfloor-\lfloor n t_{\ell-1}\rfloor-k_{n,\ell}m_n\bigr)
  \EE\bigl(|Y_1|^\beta\bigr)
   \leq \frac{2^{1-\beta}m_n}{a_n^\beta} \EE\bigl(|Y_1|^\beta\bigr) \sum_{\ell=1}^d |\vartheta_\ell|^\beta ,
 \end{align*}
 since \ $\lfloor n t_\ell\rfloor-\lfloor n t_{\ell-1}\rfloor-k_{n,\ell}m_n \leq  \lfloor n t_\ell\rfloor - \lfloor n t_{\ell-1}\rfloor
   - \Big(\frac{\lfloor n t_\ell\rfloor  - \lfloor n t_{\ell-1}\rfloor }{m_n} - 1\Big)m_n = m_n$.
To be precise, if \ $(Y_k)_{k\in\NN}$ \ is \ $M_0$-dependent, then choose a \ $\beta\in(\frac{\alpha}{2}, \alpha\wedge 1)$ \ in the above
 calculations, otherwise such a \ $\beta$ \ is given via \eqref{vare}, and especially in both cases \ $\beta\in(0,1)$.
\ The right-hand side of the above inequality approaches zero if \ $m_n/a_n^\beta \to 0$ \ as \ $n \to \infty$, \
 since \ $\EE(|Y_1|^\beta)<\infty$ \ due to the facts that \ $Y_1$ \ is regularly varying with index \ $\alpha$ \ and \ $\beta < \alpha$.

In a similar way, for an arbitrary sequence \ $(\delta_n)_{n\in\NN}$ \ in \ $\NN$ \ satisfying \ $\delta_n < m_n$ \ for sufficently large \ $n \in \NN$, \ we obtain
 \begin{align*}
  &\biggl|\EE\biggl(\exp\biggl\{\frac{\ii}{a_n} \sum_{\ell=1}^d \vartheta_\ell \bigl(S_{\lfloor n t_{\ell-1}\rfloor+k_{n,\ell}m_n} - S_{\lfloor n t_{\ell-1}\rfloor}\bigr)\biggr\} \biggl) \\
  &\,  - \EE\biggl(\exp\biggl\{\frac{\ii}{a_n} \sum_{\ell=1}^d \vartheta_\ell \sum_{j=1}^{k_{n,\ell}} \bigl(S_{\lfloor n t_{\ell-1}\rfloor+jm_n-\delta_n} - S_{\lfloor n t_{\ell-1}\rfloor+(j-1)m_n}\bigr)\biggr\}\biggl)\biggr| \\
  &= \biggl|\EE\biggl(\exp\biggl\{\frac{\ii}{a_n} \sum_{\ell=1}^d \vartheta_\ell \sum_{j=1}^{k_{n,\ell}} \bigl(S_{\lfloor n t_{\ell-1}\rfloor+jm_n-\delta_n} - S_{\lfloor n t_{\ell-1}\rfloor+(j-1)m_n}\bigr)\biggr\} \\
  &\phantom{= \biggl|\EE\biggl(}
                      \times
                      \biggl(\exp\biggl\{\frac{\ii}{a_n} \sum_{\ell=1}^d \vartheta_\ell \sum_{j=1}^{k_{n,\ell}} \bigl(S_{\lfloor n t_{\ell-1}\rfloor+jm_n} - S_{\lfloor n t_{\ell-1}\rfloor+jm_n-\delta_n}\bigr)\biggr\} - 1\biggr)\biggr)\biggr| \leq
 \end{align*}
 \begin{align*}
  &\leq \EE\biggl(\biggl|\exp\biggl\{\frac{\ii}{a_n} \sum_{\ell=1}^d \vartheta_\ell \sum_{j=1}^{k_{n,\ell}} \bigl(S_{\lfloor n t_{\ell-1}\rfloor+jm_n} - S_{\lfloor n t_{\ell-1}\rfloor+jm_n-\delta_n}\bigr)\biggr\} - 1\biggr|\biggr)\\
  &= \EE\biggl(\biggl|\prod_{\ell=1}^d \prod_{j=1}^{k_{n,\ell}} \exp\biggl\{\frac{\ii\vartheta_\ell}{a_n} \bigl(S_{\lfloor n t_{\ell-1}\rfloor+jm_n} - S_{\lfloor n t_{\ell-1}\rfloor+jm_n-\delta_n}\bigr)\biggr\} - 1\biggr|\biggr) \\
  &\leq \EE\biggl(\sum_{\ell=1}^d \sum_{j=1}^{k_{n,\ell}} \biggl|\exp\biggl\{\frac{\ii\vartheta_\ell}{a_n} \bigl(S_{\lfloor n t_{\ell-1}\rfloor+jm_n} - S_{\lfloor n t_{\ell-1}\rfloor+jm_n-\delta_n}\bigr)\biggr\} - 1\biggr|\biggr) \\
  &= \sum_{\ell=1}^d k_{n,\ell} \EE\biggl(\biggl|\exp\biggl\{\frac{\ii\vartheta_\ell}{a_n} S_{\delta_n}\biggr\} - 1\biggr|\biggr)
   \leq \sum_{\ell=1}^d k_{n,\ell} \EE\biggl(\biggl|\frac{\vartheta_\ell}{a_n} S_{\delta_n}\biggr| \land 2\biggr)\\
  &\leq 2 \sum_{\ell=1}^d k_{n,\ell} \EE\biggl(\biggl(\frac{|\vartheta_\ell|}{2a_n} \sum_{i=1}^{\delta_n} |Y_i|\biggr) \land 1\biggr)
   \leq 2 \sum_{\ell=1}^d k_{n,\ell} \EE\biggl(\biggl(\frac{|\vartheta_\ell|}{2a_n} \sum_{i=1}^{\delta_n} |Y_i|\biggr)^\beta \land 1\biggr)\\
  &\leq 2 \sum_{\ell=1}^d k_{n,\ell} \EE\biggl(\frac{|\vartheta_\ell|^\beta}{(2a_n)^\beta} \sum_{i=1}^{\delta_n} |Y_i|^\beta\biggr)
   = \frac{2}{(2a_n)^\beta} \sum_{\ell=1}^d k_{n,\ell} |\vartheta_\ell|^\beta \delta_n \EE\bigl(|Y_1|^\beta\bigr)\\
  &\leq \frac{2^{1-\beta}\delta_n}{a_n^\beta m_n} \EE\bigl(|Y_1|^\beta\bigr) \sum_{\ell=1}^d \bigl((t_\ell - t_{\ell-1}) n + 1\bigr) |\vartheta_\ell|^\beta .
 \end{align*}
The right-hand side of the above inequality approaches zero if \ $\delta_n n/(a_n^\beta m_n) \to 0$ \ as \ $n \to \infty$.

In a similar way, we get
 \begin{align*}
  &\biggl|\prod_{\ell=1}^d \bigl(\varphi_{n,m_n}(\vartheta_\ell)\bigr)^{k_{n,\ell}}
          - \prod_{\ell=1}^d \bigl(\varphi_{n,m_n-\delta_n}(\vartheta_\ell)\bigr)^{k_{n,\ell}}\biggr|
   \leq \sum_{\ell=1}^d k_{n,\ell} |\varphi_{n,m_n}(\vartheta_\ell) - \varphi_{n,m_n-\delta_n}(\vartheta_\ell)| \\
  &= \sum_{\ell=1}^d k_{n,\ell} \biggl|\EE\biggl(\exp\biggl\{\frac{\ii\vartheta_\ell}{a_n} S_{m_n}\biggr\}\biggr) - \EE\biggl(\exp\biggl\{\frac{\ii\vartheta_\ell}{a_n} S_{m_n-\delta_n}\biggr\}\biggr)\biggr| \\
  &= \sum_{\ell=1}^d k_{n,\ell} \biggl|\EE\biggl(\exp\biggl\{\frac{\ii\vartheta_\ell}{a_n} S_{m_n-\delta_n}\biggr\} \biggl(\exp\biggl\{\frac{\ii\vartheta_\ell}{a_n} \bigl(S_{m_n} - S_{m_n-\delta_n}\bigr)\biggr\} - 1\biggr)\biggr)\biggr| \\
  &\leq \sum_{\ell=1}^d k_{n,\ell} \EE\biggl(\biggl|\exp\biggl\{\frac{\ii\vartheta_\ell}{a_n} \bigl(S_{m_n} - S_{m_n-\delta_n}\bigr)\biggr\} - 1\biggr|\biggr) \\
  &= \sum_{\ell=1}^d k_{n,\ell} \EE\biggl(\biggl|\exp\biggl\{\frac{\ii\vartheta_\ell}{a_n} S_{\delta_n}\biggr\} - 1\biggr|\biggr)
   \leq \frac{2^{1-\beta}\delta_n}{a_n^\beta m_n} \EE\bigl(|Y_1|^\beta\bigr) \sum_{\ell=1}^d \bigl((t_\ell - t_{\ell-1}) n + 1\bigr) |\vartheta_\ell|^\beta .
 \end{align*}
The right-hand side of the above inequality approaches zero if \ $\delta_n n / (a_n^\beta m_n) \to 0$ \ as \ $n \to \infty$.

Further, using the strong stationarity of \ $(Y_k)_{k\in\NN}$, \ we can write
 \[
  \EE\biggl(\exp\biggl\{\frac{\ii}{a_n} \sum_{\ell=1}^d \vartheta_\ell \sum_{j=1}^{k_{n,\ell}} \bigl(S_{\lfloor n t_{\ell-1}\rfloor+jm_n-\delta_n} - S_{\lfloor n t_{\ell-1}\rfloor+(j-1)m_n}\bigr)\biggr\}\biggr)
   = \EE\Biggl(\prod_{\ell=1}^d \prod_{j=1}^{k_{n,\ell}} U_{\ell,j}\Biggr) ,
  \]
 and
 \[
  \prod_{\ell=1}^d \bigl(\varphi_{n,m_n-\delta_n}(\vartheta_\ell)\bigr)^{k_{n,\ell}}
   = \prod_{\ell=1}^d \prod_{j=1}^{k_{n,\ell}} \EE(U_{\ell,j})
 \]
 with
 \[
   U_{\ell,j} := \exp\biggl\{\frac{\ii\vartheta_\ell}{a_n} \bigl(S_{\lfloor n t_{\ell-1}\rfloor+jm_n-\delta_n} - S_{\lfloor n t_{\ell-1}\rfloor+(j-1)m_n}\bigr)\biggr\} .
 \]
We have
 \begin{align*}
  &\Biggl|\EE\Biggl(\prod_{\ell=1}^d \prod_{j=1}^{k_{n,\ell}} U_{\ell,j}\Biggr)
          - \prod_{\ell=1}^d \prod_{j=1}^{k_{n,\ell}} \EE(U_{\ell,j})\Biggr|\\
  &\leq \Biggl|\EE\Biggl(\prod_{\ell=1}^d \prod_{j=1}^{k_{n,\ell}} U_{\ell,j}\Biggr)
               - \EE(U_{1,1})
                 \EE\Biggl(\prod_{j=2}^{k_{n,1}} U_{1,j}
                           \prod_{\ell=2}^d \prod_{j=1}^{k_{n,\ell}} U_{\ell,j}\Biggr)\Biggr| \\
  &\phantom{\leq}
         + |\EE(U_{1,1})|
           \Biggl|\EE\Biggl(\prod_{j=2}^{k_{n,1}} U_{1,j}
                            \prod_{\ell=2}^d \prod_{j=1}^{k_{n,\ell}} U_{\ell,j}\Biggr)
                  - \Biggl(\prod_{j=2}^{k_{n,1}} \EE(U_{1,j})\Biggr)
                    \Biggl(\prod_{\ell=2}^d \prod_{j=1}^{k_{n,\ell}} \EE(U_{\ell,j})\Biggr)\Biggr| ,
 \end{align*}
 where we check that the first term on the right-hand side is bounded by \ $16 \alpha_{\delta_n +1}$.
\ Namely, the strong stationarity of \ $(Y_k)_{k\in\NN}$ \ implies that
 \[
   |\EE(V_1 V_2) - \EE(V_1) \EE(V_2)| \leq 4 C_1 C_2 \alpha_h
 \]
 for any \ $\cF^Y_{1,j}$-measurable (real-valued) random variable \ $V_1$ \ and any \ $\cF^Y_{j+h,\infty}$-measurable (real-valued)
 random variable \ $V_2$ \ with \ $j, h \in \NN$, \ $|V_1| \leq C_1$ \ and \ $|V_2| \leq C_2$ \ (see, e.g., Lemma 1.2.1 in Lin and Lu \cite{LinLu}).
Hence, for any \ $\cF^Y_{1,j}$-measurable complex-valued random variable \ $V_1$ \ and any \ $\cF^Y_{j+h,\infty}$-measurable complex-valued random variable \ $V_2$ \ with \ $j, h \in \NN$, \ $|V_1| \leq C_1$ \ and \ $|V_2| \leq C_2$, \ we get
 \begin{align*}
  &|\EE(V_1 V_2) - \EE(V_1) \EE(V_2)| \\
  &\leq |\EE(\Re(V_1) \Re(V_2)) - \EE(\Re(V_1)) \EE(\Re(V_2))|
       + |\EE(\Im(V_1) \Im(V_2)) - \EE(\Im(V_1)) \EE(\Im(V_2))| \\
  &\phantom{\leq}
       + |\EE(\Re(V_1) \Im(V_2)) - \EE(\Re(V_1)) \EE(\Im(V_2))|
       + |\EE(\Im(V_1) \Re(V_2)) - \EE(\Im(V_1)) \EE(\Re(V_2))| \\
  &\leq 16 C_1 C_2 \alpha_h .
 \end{align*}
Then one can apply this inequality with the choices \ $V_1 = U_{1,1}$, \ $V_2 = \prod_{j=2}^{k_{n,1}} U_{1,j} \prod_{\ell=2}^d \prod_{j=1}^{k_{n,\ell}} U_{\ell,j}$ \ and \ $C_1 = C_2 = 1$.

Iterative use of the previous argument, recursively on distinct blocks, shows that
  \begin{align*}
  &\biggl|\EE\biggl(\exp\biggl\{\frac{\ii}{a_n} \sum_{\ell=1}^d \vartheta_\ell \sum_{j=1}^{k_{n,\ell}} \bigl(S_{\lfloor n t_{\ell-1}\rfloor+jm_n-\delta_n}
        - S_{\lfloor n t_{\ell-1}\rfloor+(j-1)m_n}\bigr)\biggr\}\biggr)
           - \prod_{\ell=1}^d \bigl(\varphi_{n,m_n-\delta_n}(\vartheta_\ell)\bigr)^{k_{n,\ell}}\biggr|  \\
    & \leq 16 \alpha_{\delta_n + 1} \sum_{\ell=1}^d k_{n,\ell}
      \leq 16  \alpha_{\delta_n +1} \sum_{\ell=1}^d \frac{\lfloor n t_\ell\rfloor - \lfloor n t_{\ell-1}\rfloor}{m_n}
      \leq  \frac{16\alpha_{\delta_n + 1}}{m_n} \sum_{\ell=1}^d (n t_\ell - n t_{\ell-1} +1) ,
  \end{align*}
 hence, taking into account the fact that the function \ $(\alpha_h)_{h\in\NN}$ \ is non-increasing,
 the right-hand side of the above inequality approaches zero if \ $\alpha_{\delta_n} n / m_n \to 0$ \ as \ $n \to \infty$.

Thus we proved that condition (ii) of Theorem \ref{BJMW1} is satisfied if
 \begin{equation}\label{ii}
  m_n\to\infty, \qquad \frac{n}{m_n} \to\infty, \qquad
   \frac{n\alpha_{\delta_n}}{m_n} \to 0 , \qquad \frac{n\delta_n}{a_n^\beta m_n} \to 0 \qquad \text{and} \qquad \frac{m_n}{a_n^\beta} \to 0
 \end{equation}
 as \ $n \to \infty$.

If \ $(Y_k)_{k\in\NN}$ \ is \ $M_0$-dependent with some \ $M_0 \in \ZZ_+$, \ then \ $\alpha_h = 0$ \ for all \ $h > M_0$.
\ Choose \ $\beta \in \bigl(\frac{\alpha}{2}, \alpha \land 1\bigr)$ \ and \ $\gamma \in \bigl(\frac{1}{2}, \frac{\beta}{\alpha}\bigr)$ \ and put \ $m_n := \lfloor n^\gamma\rfloor$, \ $\delta_n := \lfloor m_n^2/n\rfloor$ \ and \ $r_n := \lfloor m_n^2/a_n^\beta\rfloor$ \ for all \ $n \in \NN$.
\ Then \eqref{ii}, \ $\delta_n < m_n$ \ for sufficently large \ $n \in \NN$, \ $r_n\to\infty$ \ and \ $r_n/m_n\to 0$ \ as \ $n\to\infty$, \ and \eqref{r} hold.
Indeed, using the inequality \ $x - 1 < \lfloor x\rfloor \leq x$ \ for \ $x \in \RR$, \ we have
 \[
   \frac{n}{m_n} \geq \frac{n}{n^\gamma} = n^{1-\gamma} \to \infty \qquad \text{as \ $n \to \infty$,}
 \]
 since \ $\gamma < \frac{\beta}{\alpha} < 1$, \ and
 \begin{equation}\label{HELP}
  \begin{aligned}
   \frac{n\delta_n}{a_n^\beta m_n}
   &\leq \frac{n(m_n^2/n)}{a_n^\beta(n^\gamma-1)}
    \leq \frac{n^{2\gamma}}{a_n^\beta(n^\gamma-1)}
    = \frac{n^\gamma}{a_n^\beta(1-n^{-\gamma})} \\
   &= \frac{n^\gamma}{(n^{1/\alpha}L(n))^\beta(1-n^{-\gamma})}
    = \frac{n^{\gamma-\frac{\beta}{\alpha}}L(n)^{-\beta}}{1-n^{-\gamma}}
    \to 0 \qquad \text{as \ $n \to \infty$,}
  \end{aligned}
 \end{equation}
 since \ $0 < \gamma < \frac{\beta}{\alpha}$ \ and \ $L$ \ (and hence \ $L^{-\beta}$) \ is slowly varying, and
 \[
   \delta_n
   \geq \frac{m_n^2}{n} - 1
   \geq \frac{(n^\gamma-1)^2}{n} - 1
   = \bigl(n^{\gamma-\frac{1}{2}} - n^{-\frac{1}{2}}\bigr)^2 - 1 \to \infty \qquad \text{as \ $n \to \infty$,}
 \]
 since \ $\gamma > \frac{1}{2}$, \ yielding
 \[
   \frac{n\alpha_{\delta_n}}{m_n}
   = \frac{n\cdot0}{m_n} = 0 \qquad \text{if \ $\delta_n > M_0$ \ (which holds for large enough \ $n\in\NN$),}
 \]
 and, by \eqref{HELP},
 \[
   \frac{m_n}{a_n^\beta}
   \leq \frac{n^\gamma}{a_n^\beta}
   \to 0 \qquad \text{as \ $n \to \infty$.}
 \]
Further,
 \[
   \frac{\delta_n}{m_n} \leq \frac{m_n^2/n}{m_n} = \frac{m_n}{n} \to 0 \qquad \text{as \ $n \to \infty$,}
 \]
 hence \ $\delta_n < m_n$ \ for sufficently large \ $n \in \NN$, \ and
 \begin{align*}
  r_n
  &\geq \frac{m_n^2}{a_n^\beta} - 1
   \geq \frac{(n^\gamma-1)^2}{a_n^\beta} - 1
   = \frac{(n^\gamma-1)^2}{(n^{1/\alpha}L(n))^\beta} - 1 \\
  &= (1-n^{-\gamma})^2 n^{2(\gamma-\frac{\beta}{2\alpha})} L(n)^{-\beta} - 1
   \to \infty \qquad \text{as \ $n \to \infty$,}
 \end{align*}
 since \ $\gamma - \frac{\beta}{2\alpha} > \frac{1}{2} - \frac{\beta}{2\alpha} > 0$ \ and \ $L$ \ is slowly varying, and
 \[
   \frac{r_n}{m_n}
   \leq \frac{m_n^2/a_n^\beta}{m_n}
   = \frac{m_n}{a_n^\beta} \to 0 \qquad \text{as \ $n \to \infty$,}
 \]
 and
 \[
   n \alpha_{r_n} = n \cdot 0 = 0 \qquad \text{if \ $r_n > M_0$ \ (which holds for large enough \ $n \in \NN$).}
 \]

If \ $(Y_k)_{k\in\NN}$ \ is not \ $M_0$-dependent for any \ $M_0 \in \ZZ_+$, \ then choose \ $m_n := \lfloor\vare_n^{1/2} a_n^\beta\rfloor$ \ and put \ $\delta_n := \lfloor\vare_n a_n^{2\beta}/n\rfloor$ \ and \ $r_n := \lfloor\vare_n a_n^\beta\rfloor$ \ for all \ $n \in \NN$, \ where \ $\beta\in(\frac{\alpha}{2},\alpha\wedge 1)$ \
 is given via  \eqref{vare}.
Then we show that condition \eqref{vare} implies that \eqref{ii}, \ $\delta_n < m_n$ \ for sufficiently
 large \ $n \in \NN$, \ $r_n\to\infty$ \ and \ $r_n/m_n\to 0$ \ as \ $n\to\infty$, \ and \eqref{r} hold.
Condition \eqref{vare} implies \ $\alpha_{\lfloor\vare_n a_n^{2\beta}/n\rfloor} = n^{-1} (n \alpha_{\lfloor\vare_n a_n^{2\beta}/n\rfloor}) \to 0$ \ as \ $n\to\infty$.
\ We check that \ $\delta_n = \lfloor\vare_n a_n^{2\beta}/n\rfloor \to \infty$ \ as \ $n \to \infty$.
\ If, on the contrary, we suppose that \ $\lfloor\vare_n a_n^{2\beta}/n\rfloor$ \ does not converge to \ $\infty$ \ as \ $n \to \infty$, \ then there exists a \ $K\in\RR_{++}$ \ such that for all \ $N\in\NN$ \ we have \ $\lfloor\vare_n a_n^{2\beta}/n\rfloor<K$ \ with some \ $n\geq N$, \
 i.e., there exists a sequence \ $(n_k)_{k\in\NN}$ \ such that \ $n_k\to\infty$ \ as \ $k\to\infty$ \ and \ $\lfloor \vare_{n_k} a_{n_k}^{2\beta}/n_k\rfloor<K$, \ $k\in\NN$. \
Since the function \ $\alpha$ \ is non-increasing and strictly positive (due to the fact that \ $(Y_k)_{k\in\NN}$ \ is not \ $M_0$-dependent for any \ $M_0\in\ZZ_+$),
 \ we have \ $\alpha_{\lfloor \vare_{n_k} a_{n_k}^{2\beta}/n_k\rfloor}\geq \alpha_K>0$, \ $k\in\NN$, \
 which leads us to a contradiction, since \ $\alpha_{\lfloor\vare_{n_k} a_{n_k}^{2\beta}/n_k\rfloor}\to 0$ \ as \ $k\to\infty$.
\ Clearly, we also obtain \ $\vare_n a_n^{2\beta}/n \to \infty$ \ as \ $n \to \infty$.
\ Consequently,
 \[
   m_n \geq \vare_n^{1/2} a_n^\beta - 1
   = \left(\frac{\vare_n a_n^{2\beta}}{n}\right)^{1/2} n^{1/2} - 1 \to\infty
    \qquad \text{as \ $n \to \infty$,}
 \]
 and, by \eqref{a_n_beta_n} (being equivalent to \eqref{vare}),
 \[
   \frac{m_n}{n} \leq \frac{\vare_n^{1/2}a_n^\beta}{n}
   \to 0 \qquad \text{as \ $n \to \infty$,}
 \]
 and, by \eqref{vare} and \ $m_n \to \infty$ \ as \ $n \to \infty$,
 \[
   \frac{n\alpha_{\delta_n}}{m_n}
   = \frac{n\alpha_{\lfloor\vare_n a_n^{2\beta}/n\rfloor}}{m_n}
   \to 0 \qquad \text{as \ $n \to \infty$,}
 \]
 and
 \[
   \frac{n\delta_n}{a_n^\beta m_n}
   \leq \frac{n(\vare_n a_n^{2\beta}/n)}{a_n^\beta(\vare_n^{1/2} a_n^\beta-1)}
   = \frac{\vare_n^{1/2}}{1-\vare_n^{-1/2}a_n^{-\beta}} \to 0  \qquad \text{as \ $n \to \infty$}
 \]
 since \ $\vare_n^{-1/2}a_n^{-\beta} \leq m_n^{-1} \to 0$ \ as \ $n \to \infty$, \ and
 \[
   \frac{m_n}{a_n^\beta}
   \leq \frac{\vare_n^{1/2}a_n^\beta}{a_n^\beta}
   = \vare_n^{1/2}
   \to 0 \qquad \text{as \ $n\to\infty$.}
 \]
Further, by \eqref{a_n_beta_n},
 \[
   \frac{\delta_n}{m_n} = \frac{n\delta_n}{a_n^\beta m_n} \frac{a_n^\beta}{n} \to 0  \qquad \text{as \ $n \to \infty$,}
 \]
 hence \ $\delta_n < m_n$ \ for sufficently large \ $n \in \NN$, \ and, by \ $\vare_n a_n^{2\beta}/n \to \infty$ \ as \ $n \to \infty$ \ and \eqref{a_n_beta_n},
 \[
   r_n \geq \vare_n a_n^\beta - 1
   = \frac{\vare_n a_n^{2\beta}}{n} \frac{n}{a_n^\beta} - 1
   \to \infty \qquad \text{as \ $n \to \infty$,}
 \]
 and
 \[
   \frac{r_n}{m_n}
   \leq \frac{\vare_n a_n^\beta}{\vare^{1/2} a_n^\beta-1}
   = \frac{\vare_n^{1/2}}{1-\vare_n^{-1/2}a_n^{-\beta}} \to 0  \qquad \text{as \ $n\to\infty$,}
 \]
 and
 \[
   n\alpha_{r_n}
   = n \alpha_{\lfloor \vare_n a_n^\beta\rfloor}
   \leq n \alpha_{\lfloor \vare_n a_n^{2\beta}/n\rfloor} \to 0 \qquad \text{as \ $n\to\infty$,}
 \]
 since the function \ $\alpha$ \ is non-increasing, \ $a_n^\beta\geq a_n^{2\beta}/n$ \ for sufficiently large \ $n\in\NN$ \ (see \eqref{a_n_beta_n}),
 and one can use \eqref{vare}.

(ii). \ Let us suppose that \ $(Y_k)_{k\in\NN}$ \ is strongly mixing with geometric rate.
Recall that \ $\alpha_h \leq C_{\mathrm{mix}} q^h$ \ for all \ $h \in \NN$ \ with some \ $q\in(0,1)$ \ and \ $C_{\mathrm{mix}}:=\sup_{h\in\NN} q^{-h} \alpha_h<\infty$.
\ Next, let us choose \ $\gamma_1, \gamma_2 \in (0, 1)$ \ satisfying \ $\gamma_1 < \gamma_2$ \ and \ $\gamma_2 \in \bigl(\frac{1}{2}, \frac{1}{\alpha} \land 1\bigr)$, \ choose \ $\beta \in \bigl(\gamma_2 \alpha, \alpha \land 1\bigr)$, \ and put \ $\delta_n := \lfloor n^{2\gamma_2-1} \rfloor$, \ $n \in \NN$.
\ Then \eqref{ii}, \ $\delta_n < m_n$ \ for sufficiently large \ $n \in \NN$, \ $r_n \to\infty$ \ and \ $r_n/m_n\to 0$ \ as \ $n\to\infty$, \ and \eqref{r} hold.
Indeed,
 \[
   \frac{n}{m_n} \leq \frac{n}{n^{\gamma_2}-1} = \frac{n^{1-\gamma_2}}{1-n^{-\gamma_2}} \to \infty \qquad \text{as \ $n \to \infty$,}
 \]
 since \ $0 < \gamma_2 < 1$, \ and
 \[
   \frac{n\alpha_{\delta_n}}{m_n}
   \leq \frac{C_{\mathrm{mix}}nq^{\delta_n}}{n^{\gamma_2}-1}
   \leq \frac{C_{\mathrm{mix}}n^{1-\gamma_2} q^{n^{2\gamma_2-1}-1}}{1-n^{-\gamma_2}} \to 0  \qquad \text{as \ $n \to \infty$,}
 \]
 since \ $\frac{1}{2} < \gamma_2 < 1$, \ and \ $\lim_{n\to\infty} n^{\kappa_1} q^{n^{\kappa_2}} = 0$ \ for any \ $\kappa_1, \kappa_2 \in \RR_{++}$, \ and
 \[
   \frac{n\delta_n}{a_n^\beta m_n}
   \leq \frac{n\cdot n^{2\gamma_2-1}}{(n^{1/\alpha}L(n))^\beta(n^{\gamma_2}-1)}
   = \frac{n^{\gamma_2-\frac{\beta}{\alpha}}}{L(n)^\beta(1-n^{-\gamma_2})}
   \to 0 \qquad \text{as \ $n \to \infty$,}
 \]
 since \ $0 < \gamma_2 < \frac{\beta}{\alpha}$ \ and \ $L$ \ is slowly varying, and
 \[
   \frac{m_n}{a_n^\beta}
   \leq \frac{n^{\gamma_2}}{(n^{1/\alpha}L(n))^\beta}
   = \frac{n^{\gamma_2-\frac{\beta}{\alpha}}}{L(n)^\beta}
   \to 0 \qquad \text{as \ $n \to \infty$,}
 \]
 and
 \[
   \frac{\delta_n}{m_n} \leq \frac{n^{2\gamma_2-1}}{n^{\gamma_2}-1} = \frac{n^{\gamma_2-1}}{1-n^{-\gamma_2}} \to 0 \qquad \text{as \ $n \to \infty$,}
 \]
 since \ $0 < \gamma_2 < 1$, \ hence \ $\delta_n < m_n$ \ for sufficently large \ $n \in \NN$, \ and
 \[
   \frac{r_n}{m_n}
   \leq \frac{n^{\gamma_1}}{n^{\gamma_2}-1}
   = \frac{n^{\gamma_1-\gamma_2}}{1-n^{-\gamma_2}} \to 0 \qquad \text{as \ $n \to \infty$,}
 \]
 since \ $0 < \gamma_1 < \gamma_2$, \ and
 \[
   n \alpha_{r_n}
   \leq C_{\mathrm{mix}} n q^{r_n}
   \leq C_{\mathrm{mix}} n q^{n^{\gamma_1}-1}
   \to 0 \qquad \text{as \ $n \to \infty$,}
 \]
 since \ $\gamma_1 > 0$.
\proofend

\noindent{\bf Proof of part (ii) of Remark \ref{Rem_Lemma3}.}
Indeed, there exists a constant \ $q \in (0, 1)$ \ such that \ $\alpha_h = \OO(q^h)$ \ as \ $h \to \infty$, \ i.e., \ $C_{\mathrm{mix}} := \sup_{h\in\NN} q^{-h} \alpha_h < \infty$, \ yielding \ $\alpha_h \leq C_{\mathrm{mix}} q^h$ \ for all \ $h \in \NN$.
\ Since \ $c\in(0,\left(\frac{2}{\alpha} -1\right)\wedge 1)$, \ we can choose \ $\beta \in \bigl(\frac{(c+1)\alpha}{2}, \alpha \land 1\bigr)$,
 \ hence \ $c \in \big(0, \frac{2\beta}{\alpha} -1\big)$.
\ Then, for each \ $n \in \NN$, \ we have
 \[
   n \alpha_{\lfloor\vare_n a_n^{2\beta}/n\rfloor}
   = n \alpha_{\lfloor n^{-1-c} a_n^{2\beta}\rfloor}
   \leq C_{\mathrm{mix}} n q^{\lfloor n^{-1-c} a_n^{2\beta}\rfloor} .
 \]
For sufficiently large \ $n \in \NN$, \ we have
 \[
   \lfloor n^{-1-c} a_n^{2\beta}\rfloor
   \geq n^{-1-c} a_n^{2\beta} - 1
   = n^{-1-c} (n^{1/\alpha} L(n))^{2\beta} - 1
   = n^{\frac{2\beta}{\alpha}-1-c} L(n)^{2\beta} - 1
   \geq n^{(\frac{2\beta}{\alpha}-1-c)/2} - 1 ,
 \]
 since \ $\frac{2\beta}{\alpha} - 1 - c > 0$ \ and \ $L$ \ is slowly varying, hence
 \[
   n \alpha_{\lfloor\vare_n a_n^{2\beta}/n\rfloor}
   \leq C_{\mathrm{mix}} n q^{n^{(\frac{2\beta}{\alpha}-1-c)/2}-1}
   \to 0 \qquad \text{as \ $n\to\infty$,}
 \]
 since \ $\lim_{n\to\infty} nq^{n^\kappa} = 0$ \ for any \ $\kappa \in \RR_{++}$, \ yielding \eqref{vare}.
\proofend

\noindent{\bf Proof of Lemma \ref{Lemma2}.}
For all \ $d\in\NN$, \ $\vartheta\in\RR$, \ $\vartheta\ne 0$, \ and large enough \ $n\in\NN$ \ satisfying \ $m_n\geq d+1$ \
  (which can be assumed since \ $m_n\to\infty$ \ as \ $n\to\infty$), using that \ $\vert \overline{Z}\vert = \vert Z\vert\wedge 2$ \
  for any random variable \ $Z$, \ we have
 \begin{align*}
   &\frac{n}{m_n} \sum_{j=d+1}^{m_n} \EE\bigl(\bigl|\overline{\vartheta a_n^{-1} (S_j - S_d)} \;\cdot\; \overline{\vartheta a_n^{-1} Y_1}\bigr|\bigr)\\
   &\qquad = \frac{n}{m_n} \sum_{j=d+1}^{m_n} \EE\bigl( (\vert \vartheta a_n^{-1} (S_j - S_d) \vert\wedge 2 )  ( \vert \vartheta a_n^{-1} Y_1 \vert \wedge 2) \bigr)\\
   &\qquad \leq \frac{n}{m_n} \sum_{j=d+1}^{m_n} \sum_{\ell=d+1}^j
                \EE\bigl( (\vert \vartheta a_n^{-1} Y_\ell \vert\wedge 2 )  ( \vert \vartheta a_n^{-1} Y_1 \vert \wedge 2) \bigr)\\
   &\qquad = \frac{n}{m_n} \sum_{\ell=d+1}^{m_n} \sum_{j=\ell}^{m_n}
                 \EE\bigl( (\vert \vartheta a_n^{-1} Y_\ell \vert\wedge 2 )  ( \vert \vartheta a_n^{-1} Y_1 \vert \wedge 2) \bigr)\\
   & \qquad = \frac{n}{m_n} \sum_{\ell=d+1}^{m_n} (m_n-\ell+1)
                   \EE\bigl( (\vert \vartheta a_n^{-1} Y_\ell \vert\wedge 2 )  ( \vert \vartheta a_n^{-1} Y_1 \vert \wedge 2) \bigr)\\
   & \qquad \leq n \frac{m_n+1}{m_n}\sum_{\ell=d+1}^{m_n}
                   \EE\bigl( (\vert \vartheta a_n^{-1} Y_\ell \vert\wedge 2 )  ( \vert \vartheta a_n^{-1} Y_1 \vert \wedge 2) \bigr)\\
   & \qquad = n \frac{m_n+1}{m_n} \vartheta^2 \sum_{\ell=d+1}^{m_n}
                   \EE\bigl( (\vert a_n^{-1} Y_\ell \vert\wedge (2 \vert \vartheta \vert^{-1}) )
                              ( \vert a_n^{-1} Y_1 \vert \wedge (2\vert \vartheta \vert^{-1}) ) \bigr),
 \end{align*}
 where \ $\frac{m_n+1}{m_n}\to 1$ \ as \ $n\to\infty$.
\ So first taking \ $\limsup_{n\to\infty}$ \ and then \ $\lim_{d\to\infty}$, \  condition \eqref{3.9} yields condition (iii) of
  Theorem \ref{BJMW1} in case of \ $\vartheta\in\RR$, \ $\vartheta\ne 0$.
\ For \ $\vartheta=0$, \ condition (iii) of Theorem \ref{BJMW1} trivially holds.
\proofend

\noindent{\bf Proof of Theorem \ref{aggr_time}.}
We are going to apply Theorem \ref{BJMW1} for \ $(X_n)_{n\in\NN}$ \ if \ $\alpha \in (0, 1)$, \
 and for \ $(X_n - \EE(X_n))_{n\in\NN}$ \ if \ $\alpha \in (1, \frac{4}{3})$.

{\sl Step 1 (checking conditions (i) and (ii) of Theorem \ref{BJMW1}).}
Condition (i) of Theorem \ref{BJMW1} holds, since \ $(X_n)_{n\in\NN}$ \ is strongly stationary and jointly regularly varying with index \ $\alpha \in (0, 2)$, \ see Basrak et al.\ \cite{BasKulPal} (or Theorem \ref{Xtailprocess}), and, by parts (i) and (ii) of Lemma \ref{Lem_process_shift}, \ $(X_n - \EE(X_n) )_{n\in\NN}$ \ is also strongly stationary and jointly regularly varying with index \ $\alpha \in (0, 2)$.
\ Lemma \ref{Lem_an_tan} shows that the sequence \ $(a_n)_{n\in\NN}$ \ appearing in condition (i) of Theorem \ref{BJMW1} satisfies both \ $n\PP(X_1 > a_n) \to 1$ \ as \ $n \to \infty$ \ and \ $n\PP(|X_1 - \EE(X_1)| > a_n) \to 1$ \ as \ $n \to \infty$.

By Lemma \ref{lemma:strong_mixing} and part (iii) of Lemma  \ref{Lem_process_shift}, \ $(X_n)_{n\in\NN}$ \ and \ $(X_n-\EE(X_n))_{n\in\NN}$ \ are strongly stationary and strongly mixing processes with the same geometric rate function.
Hence, by Lemma \ref{Lemma3}, condition (ii) of Theorem \ref{BJMW1} holds for \ $(X_n)_{n\in\NN}$ \ if \ $\alpha\in(0,1)$, \ and for \ $(X_n - \EE(X_n))_{n\in\NN}$ \
 if \ $\alpha\in(1,\frac{4}{3})$ \ with \ $m_n = \lfloor n^{\gamma_2}\rfloor$,
 \ where \ $\gamma_2 \in \bigl(\frac{1}{2}, \frac{1}{\alpha} \land 1\bigr)$ \ can be chosen arbitrarily.

{\sl Step 2 (checking condition (iii) of Theorem \ref{BJMW1}).}
We check that the anti-cluster type condition (iii) of Theorem \ref{BJMW1} holds for \ $(X_n)_{n\in\NN}$ \
 if \ $\alpha \in (0, 1)$ \ and for  \ $(X_n - \EE(X_n))_{n\in\NN}$ \  if \ $\alpha \in (1, \frac{4}{3})$.
By Step 1 and Lemma \ref{Lemma2}, to check condition (iii) of Theorem \ref{BJMW1}, it is enough to verify \eqref{3.9}.
Namely, it is enough to verify that if \ $\alpha \in (0, 1)$, \ then for all \ $x\in\RR_{++}$,
 \begin{align}\label{BJMM1}
  \lim_{d\to\infty} \limsup_{n\to\infty} n \sum_{j=d+1}^{m_n} \EE\big( (\vert a_n^{-1} X_j\vert \wedge x) (\vert a_n^{-1} X_0\vert \wedge x) \big)= 0
 \end{align}
 and if \ $\alpha \in (1, \frac{4}{3})$, \ then for all \ $x\in\RR_{++}$,
 \begin{align}\label{BJMM2}
  \lim_{d\to\infty} \limsup_{n\to\infty} n \sum_{j=d+1}^{m_n} \EE\big( (\vert a_n^{-1}(X_j - \EE(X_j))\vert \wedge x) (\vert a_n^{-1} (X_0 - \EE(X_0)) \vert \wedge x) \big)= 0 .
 \end{align}

First, we recall a representation of \ $(X_0,X_1,\ldots,X_n)$, \ which will be useful for checking \eqref{BJMM1} and \eqref{BJMM2}.
Introducing the notation \ $\Pi^{(0)}_n := \theta^{(0)}_n \circ \cdots \circ \theta^{(0)}_1$, \ $n \in \NN$, \ by Lemma \ref{lem:reprXn}, we have
 \begin{align}\label{help_GWI_rep}
   (X_0,X_1,\ldots,X_n)
    \distre
    (X_0,\kappa_1+ \Pi^{(0)}_1\circ X_0,\ldots, \kappa_n + \Pi^{(0)}_n\circ X_0),
    \qquad n\in\NN,
 \end{align}
 where $\kappa_n$, $\Pi^{(0)}_n\circ j$ and $X_0$ are independent for each $n,j\in\NN$,
 and $\kappa_n$ is given in Lemma \ref{lem:reprXn}.

{\sl Step 2/(a).}  We check \eqref{BJMM1}.
Let \ $x\in\RR_{++}$ \ be fixed.
For each \ $d\in\NN$ \ and large enough \ $n\in\NN$ \ satisfying \ $m_n\geq d+1$ \ (which can be assumed since
 \ $m_n\to\infty$ \ as \ $n\to\infty$), \ we have
 \begin{align*}
  & n \sum_{j=d+1}^{m_n} \EE\big( (\vert a_n^{-1} X_j\vert \wedge x) (\vert a_n^{-1} X_0\vert \wedge x) \big)  \\
  & \qquad = x^2 n \sum_{j=d+1}^{m_n} \PP( \vert X_j\vert > a_n x, \vert X_0\vert > a_n x )
             + n \sum_{j=d+1}^{m_n} \EE\big( \vert a_n^{-1} X_j\vert \vert a_n^{-1} X_0\vert \bone_{\{ \vert X_j\vert \leq a_n x,  \vert X_0\vert \leq a_n x  \}}\big) \\
  & \phantom{\qquad =\,} + xn \sum_{j=d+1}^{m_n} \EE\big( \vert a_n^{-1} X_j\vert \bone_{\{ \vert X_j\vert \leq a_n x,  \vert X_0\vert > a_n x  \}}\big)
                         + xn \sum_{j=d+1}^{m_n} \EE\big( \vert a_n^{-1} X_0\vert \bone_{\{ \vert X_j\vert > a_n x,  \vert X_0\vert \leq a_n x  \}}\big) \\
  & \qquad =: \mathrm{I}_{d,n,1}+ \mathrm{I}_{d,n,2} + \mathrm{I}_{d,n,3} + \mathrm{I}_{d,n,4}.
 \end{align*}
We check that \ $\lim_{d\to\infty} \limsup_{n\to\infty} \mathrm{I}_{d,n,i} = 0$, \ $i=1,2,3,4$, \ yielding \eqref{BJMM1}.
Here, by \eqref{help_GWI_rep}, we have
 \begin{align*}
   \mathrm{I}_{d,n,1} &= x^2 n \PP(X_0 > a_n x) \sum_{j=d+1}^{m_n} \PP( X_j > a_n x \mid X_0 > a_n x ) \\
             &\leq x^2 n \PP( X_0 > a_n x) \sum_{j=d+1}^{m_n}
                      \PP( \Pi^{(0)}_j \circ X_0 > a_n x/2 \mid X_0 > a_n x )\\
             &\phantom{\leq\,}
                  + x^2 n \PP( X_0 > a_n x) \sum_{j=d+1}^{m_n}
                              \PP(\kappa_j > a_n x/2 \mid X_0 > a_nx)\\
             &=: x^2 n \PP( X_0 > a_n x) \mathrm{I}_{d,n,1}^* + x^2 n \PP( X_0 > a_n x) \mathrm{I}_{d,n,1}^{**}.
 \end{align*}
Since \ $X_0$ \ is regularly varying with index \ $\alpha\in(0,1)$, \ we have
 \begin{align}\label{help_a_n_x}
    \lim_{n\to\infty} n \PP( X_0 > a_n x) = \lim_{n\to\infty} n \PP( X_0 > a_n)  \lim_{n\to\infty} \frac{\PP( X_0 > a_n x)}{\PP( X_0 > a_n)}
                                          = 1\cdot x^{-\alpha} = x^{-\alpha}.
 \end{align}
Further, for any \ $\beta\in(0,\alpha\wedge 1) = (0,\alpha)$, \ using the independence of \ $\Pi^{(0)}_j\circ i$ \ and \ $X_0$ \ for each \ $i\in\ZZ_+$ \ and \ $j\in\NN$,
 \ by Markov's and conditional Jensen's inequalities, and the equality \ $\EE(\Pi^{(0)}_j \circ i) = m_\xi^j i$, \ $i\in\ZZ_+$, $j\in\NN$, \  we have
 \begin{align*}
   \mathrm{I}_{d,n,1}^*
   & = \sum_{j=d+1}^{m_n}  \frac{\PP\left(\Pi^{(0)}_j\circ X_0 > a_nx/2, X_0> a_nx \right)}{\PP( X_0> a_n x)}
     = \sum_{j=d+1}^{m_n} \frac{\PP((\Pi^{(0)}_j\circ X_0)\bone_{\{X_0>a_nx\}} > a_nx/2)}{\PP( X_0> a_n x)}\\
   & \leq \sum_{j=d+1}^{m_n} \frac{\PP\left( (\Pi^{(0)}_j\circ X_0)^\beta \bbone_{\{X_0> a_nx\}}> (a_nx/2)^\beta \right)}{\PP( X_0> a_nx)}
    \leq \sum_{j=d+1}^{m_n} \frac{\EE\left((\Pi^{(0)}_j\circ X_0)^\beta \bbone_{\{X_0> a_nx\}}\right)}{(a_nx/2)^\beta \PP( X_0> a_nx)} \\
   &= \sum_{j=d+1}^{m_n} \frac{\EE\left(\EE\left((\Pi^{(0)}_j\circ X_0)^\beta \bbone_{\{X_0> a_nx\}} \,\big|\, X_0\right)\right)}{(a_nx/2)^\beta \PP( X_0> a_nx)}
    \leq \sum_{j=d+1}^{m_n}  \frac{\EE\left(\left(\EE\left((\Pi^{(0)}_j\circ X_0) \bbone_{\{X_0> a_nx\}} \,\big|\, X_0\right)\right)^\beta\right)}{(a_nx/2)^\beta \PP( X_0> a_nx)} =
 \end{align*}
 \begin{align*}
   &= \sum_{j=d+1}^{m_n} \frac{\EE\left((m_\xi^j X_0)^\beta \bbone_{\{X_0> a_nx\}}\right)}{(a_nx/2)^\beta \PP( X_0> a_nx)}
    \leq \sum_{j=d+1}^\infty \frac{\EE\left((m_\xi^j X_0)^\beta \bbone_{\{X_0> a_nx\}}\right)}{(a_nx/2)^\beta \PP( X_0> a_nx)} \\
   & = 2^\beta \frac{m_\xi^{(d+1)\beta}}{1-m_\xi^\beta} \frac{\EE(X_0^\beta \bbone_{\{X_0 > a_nx\}}) }{(a_nx)^\beta \PP( X_0> a_nx)}.
 \end{align*}
By Karamata's theorem (see Theorem \ref{truncated_moments}), we have
 \[
   0\leq \limsup_{n\to\infty} \mathrm{I}_{d,n,1}^*  \leq   2^\beta \frac{m_\xi^{(d+1)\beta}}{1-m_\xi^\beta} \frac{\alpha}{\alpha-\beta},
   \qquad d\in\NN,
 \]
 yielding that \ $\lim_{d\to\infty}\limsup_{n\to\infty} \mathrm{I}_{d,n,1}^*=0$ \ by the squeeze theorem.

\noindent Further, for any \ $\beta\in(0,\alpha\wedge 1)=(0,\alpha)$, \ using the independence of \ $\kappa_j$, \ $j\in\NN$, \ and \ $X_0$, \
 and Markov's inequality, we have
 \begin{align*}
   \mathrm{I}_{d,n,1}^{**}
       = \sum_{j=d+1}^{m_n} \PP(\kappa_j > a_nx/2)
       \leq \sum_{j=d+1}^{m_n} \frac{\EE(\kappa_j^\beta)}{(a_nx/2)^\beta}
       \leq \sum_{j=d+1}^{m_n} \frac{ \EE(\kappa_{j,\infty}^\beta)}{(a_nx/2)^\beta}
 \end{align*}
 with \ $\kappa_{j,\infty} := \vare_j + \sum_{i=1}^\infty \theta_j^{(j-i)}\circ \cdots\circ \theta_{j-i+1}^{(j-i)} \circ \vare_{j-i}$, \ $j \in \NN$, \ where, by Lemma \ref{lem:reprX}, the series is convergent almost surely, and we have \ $\PP(\kappa_j\leq \kappa_{j,\infty})=1$, \ $j\in\NN$.
\ Using that \ $(\kappa_{j,\infty})_{j\in\NN}$ \ is strongly stationary and \ $\kappa_{j,\infty}\distre X_j$, \ $j\in\NN$ \ (see Lemma \ref{lem:reprX}), we have
 \begin{align*}
  \mathrm{I}_{d,n,1}^{**}
    \leq \frac{(m_n-d) \EE(X_0^\beta)}{(a_nx/2)^\beta}
    &\leq (2/x)^\beta \EE(X_0^\beta) \frac{m_n}{a_n^\beta}
     = (2/x)^\beta \EE(X_0^\beta) \frac{\lfloor n^{\gamma_2} \rfloor }{(n^{1/\alpha}L(n))^\beta} \\
    & \leq (2/x)^\beta \EE(X_0^\beta)\frac{n^{\gamma_2-\frac{\beta}{\alpha}}}{(L(n))^\beta}
     \to 0\qquad \text{as \ $n\to\infty$,}
 \end{align*}
 provided that \ $\gamma_2$ \ is chosen such that \ $\gamma_2 \in (0, \frac{\beta}{\alpha})$ \
  (fulfilling also \ $\gamma_2\in (\frac{1}{2}, \frac{1}{\alpha}\wedge 1) = (\frac{1}{2}, 1) $),
 \ since \ $\EE(X_0^ \beta) < \infty$ \ and \ $L$ \ is slowly varying.
One can choose such a \ $\gamma_2$ \ in the following way: first let us choose \ $\beta$ \
 such that \ $\frac{\alpha}{2} <\beta < \alpha\wedge 1 = \alpha$, \ i.e.\ $\frac{1}{2} < \frac{\beta}{\alpha} < \frac{1}{\alpha}\wedge 1 = 1$, \
 and then let us choose \ $\gamma_2$ \ in a way that \ $\gamma_2\in (\frac{1}{2}, \frac{\beta}{\alpha})$.
\ Hence for each \ $d\in\NN$, \ we have \ $\lim_{n\to\infty}\mathrm{I}_{d,n,1}^{**} =0$,
 \ yielding that \ $\lim_{d\to\infty}\limsup_{n\to\infty} \mathrm{I}_{d,n,1}^{**}=0$.

\noindent Putting parts together, we have \ $\lim_{d\to\infty}\limsup_{n\to\infty} \mathrm{I}_{d,n,1}=0$.

We turn to check \ $\lim_{d\to\infty}\limsup_{n\to\infty} \mathrm{I}_{d,n,2}=0$.
\ For each \ $d\in\NN$ \ and large enough \ $n\in\NN$ \ satisfying \ $m_n\geq d+1$, \ by the law of total expectation, we have
 \begin{align*}
  \mathrm{I}_{d,n,2}
   = n \sum_{j=d+1}^{m_n} \sum_{\ell=0}^{\lfloor a_n x\rfloor} a_n^{-1}\ell \PP(X_0=\ell)
                     \EE(a_n^{-1} X_j \bone_{\{X_j\leq a_nx\}} \mid   X_0=\ell).
 \end{align*}
Since the function \ $g:\ZZ_+\to\RR$, \ $g(y):=\frac{y}{(a_nx)^{1-p}} \bone_{\{y\leq a_nx\}}$, \ $y\in\ZZ_+$, \
 satisfies \ $g(y)\leq V(y) := 1 + y^p$, \ $y\in\ZZ_+$, \ for any \ $p\in(0,\alpha\wedge 1) = (0,\alpha)$, \
 by the \ $V$-uniformly ergodicity of \ $(X_k)_{k\in\ZZ_+}$ \ (see Corollary \ref{Cor_ergod}), we have
 \[
    \left\vert  \EE\left( \frac{X_j}{(a_nx)^{1-p}} \bone_{\{X_j\leq a_nx\}} \,\Big \vert\,   X_0=\ell\right)
                - \EE\left( \frac{X_0}{(a_nx)^{1-p}} \bone_{\{X_0\leq a_nx\}} \right)  \right\vert
    \leq C (1+\ell^p)\varrho^j,
    \qquad j\in\ZZ_+, \ell\in\ZZ_+,
 \]
 with some \ $C>0$ \ and \ $\varrho\in(0,1)$.
\ So
 \begin{align*}
  \mathrm{I}_{d,n,2}
   &= \frac{n}{a_n^2} \sum_{j=d+1}^{m_n} \sum_{\ell=0}^{\lfloor a_n x\rfloor} \ell \PP(X_0=\ell)
                        (a_nx)^{1-p}\EE\left(\frac{X_j}{(a_nx)^{1-p}} \bone_{\{X_j\leq a_nx\}} \,\Big \vert\,  X_0=\ell\right)\\
   &\leq \frac{n}{a_n^2} \sum_{j=d+1}^{m_n} \sum_{\ell=0}^{\lfloor a_n x\rfloor} \ell \PP(X_0=\ell)
                         (a_nx)^{1-p} \left(  \EE\left( \frac{X_0}{(a_nx)^{1-p}} \bone_{\{X_0\leq a_nx\}} \right) +  C (1+\ell^p)\varrho^j \right)\\
   &\leq \frac{n}{a_n^2} \sum_{j=d+1}^{m_n}
                       \Bigg[(\EE(X_0 \bone_{\{ X_0\leq a_nx\}}))^2  + 2C\varrho^j (a_nx)^{1-p} \EE(X_0^{p+1} \bone_{\{X_0\leq a_nx\}})\Bigg]\\
   &\leq \frac{nm_n}{a_n^2}(\EE(X_0 \bone_{\{ X_0\leq a_nx\}}))^2 + \frac{n}{a_n^2} 2C (a_nx)^{1-p}\EE(X_0^{p+1} \bone_{\{X_0\leq a_nx\}}) \sum_{j=d+1}^\infty \varrho^j\\
   &= \frac{nm_n}{a_n^2}(\EE(X_0 \bone_{\{ X_0\leq a_nx\}}))^2  + 2Cx^{1-p} \frac{\varrho^{d+1}}{1-\varrho} n a_n^{-1-p}\EE(X_0^{p+1} \bone_{\{X_0\leq a_nx\}})\\
   &=: \mathrm{I}_{d,n,2}^* +  \mathrm{I}_{d,n,2}^{**}.
 \end{align*}
Here
 \begin{align*}
    \mathrm{I}_{d,n,2}^{*}
      = x^2(n\PP(X_0 > a_nx))^2 \frac{m_n}{n} \left( \frac{\EE(X_0 \bone_{\{ X_0\leq a_nx\}})}{a_nx \PP(X_0> a_nx)} \right)^2
      \to 0 \qquad \text{as \ $n\to\infty$,}
 \end{align*}
 since \ $\lim_{n\to\infty} n \PP( X_0 > a_n x)=x^{-\alpha}$ \ (see \eqref{help_a_n_x}), \ $\lim_{n\to\infty}\frac{m_n}{n}=0$, \ and, by the assumption
 \ $\alpha\in(0,1)$ \ and Karamata's theorem,
 \[
   \lim_{n\to\infty} \frac{\EE(X_0 \bone_{\{ X_0\leq a_nx\}})}{a_nx \PP(X_0> a_nx)} = \frac{\alpha}{1-\alpha}.
 \]
This yields \ $\lim_{d\to\infty}\limsup_{n\to\infty}\mathrm{I}_{d,n,2}^{*}=0$.
\ Further, by \eqref{help_a_n_x} and Karamata's theorem (which can be applied since \ $\alpha<1<p+1$),
 \begin{align*}
     \mathrm{I}_{d,n,2}^{**}
       & =2Cx^{1-p} \frac{\varrho^{d+1}}{1-\varrho} n\PP(X_0>a_nx) \frac{\EE(X_0^{p+1} \bone_{\{X_0\leq a_nx\}})}{a_n^{p+1} \PP(X_0 > a_nx)} \\
       & \to 2Cx^{1-p} \frac{\varrho^{d+1}}{1-\varrho} x^{-\alpha} \frac{\alpha}{p+1-\alpha}
         \qquad \text{as \ $n\to\infty$.}
 \end{align*}
This yields \ $\lim_{d\to\infty}\limsup_{n\to\infty}\mathrm{I}_{d,n,2}^{**}=0$.

\noindent Putting parts together, we have \ $\lim_{d\to\infty}\limsup_{n\to\infty} \mathrm{I}_{d,n,2}=0$.
\ For historical fidelity, we note that for handling the term \ $\mathrm{I}_{d,n,2}$ \ we used the same idea as in the proof
 of Lemma 6 in Basrak and Kevei \cite{BasKev}.

We turn to check \ $\lim_{d\to\infty}\limsup_{n\to\infty} \mathrm{I}_{d,n,3}=0$, \ which can be done similarly as in case of \ $\mathrm{I}_{d,n,1}$.
\ Namely, for each \ $d\in\NN$, \ large enough \ $n\in\NN$ \ satisfying \ $m_n\geq d+1$, \ and for any \ $\beta\in(0,\alpha\wedge 1) = (0,\alpha)$, \
 using also the inequality \ $\vert x+y\vert^\beta \leq \vert x\vert^\beta + \vert y\vert^\beta$, \ $x,y\in\RR$, \ we have
 \begin{align*}
   \mathrm{I}_{d,n,3}
    &  = x^2 n \PP(X_0>a_nx) \sum_{j=d+1}^{m_n} \frac{\EE\left( \frac{X_j}{a_nx}\bone_{\{ X_j \leq a_nx\}} \bone_{\{ X_0>a_nx\}} \right)}{\PP(X_0>a_nx)}\\
    & \leq x^2 n \PP(X_0>a_nx) \sum_{j=d+1}^{m_n} \frac{\EE\left( X_j^\beta \bone_{\{ X_0>a_nx\}} \right)}{(a_nx)^\beta\PP(X_0>a_nx)}\\
    & = x^2 n \PP(X_0>a_nx)\sum_{j=d+1}^{m_n} \frac{\EE\left( (\Pi_j^{(0)}\circ X_0 + \kappa_j)^\beta \bone_{\{ X_0>a_nx\}} \right)}{(a_nx)^\beta\PP(X_0>a_nx)}\\
    & \leq x^2 n \PP(X_0>a_nx)\sum_{j=d+1}^{m_n} \frac{\EE\left( (\Pi_j^{(0)}\circ X_0)^\beta \bone_{\{ X_0>a_nx\}} + \kappa_j^\beta \bone_{\{ X_0>a_nx\}}  \right)}
                                                    {(a_nx)^\beta\PP(X_0>a_nx)}\\
    &= x^2 n \PP(X_0>a_nx)\sum_{j=d+1}^{m_n}\frac{\EE\left( (\Pi_j^{(0)}\circ X_0)^\beta \bone_{\{ X_0>a_nx\}}  \right)}{(a_nx)^\beta\PP(X_0>a_nx)}
       + x^2 n \PP(X_0>a_nx)\sum_{j=d+1}^{m_n} \frac{\EE(\kappa_j^\beta)}{(a_nx)^\beta}\\
    &= x^2 n \PP(X_0>a_nx) \mathrm{I}_{d,n,3}^* + x^2 n \PP(X_0>a_nx) \mathrm{I}_{d,n,3}^{**}.
 \end{align*}
By the calculations for \ $\mathrm{I}_{d,n,1}^*$, \ we have
 \[
    \mathrm{I}_{d,n,3}^* \leq \frac{m_\xi^{(d+1)\beta}}{1-m_\xi^\beta} \frac{\EE(X_0^\beta \bbone_{\{X_0 > a_nx\}}) }{(a_nx)^\beta \PP( X_0> a_nx)},
 \]
 and, as we have seen, an application of Karamata's theorem yields \ $\lim_{d\to\infty}\limsup_{n\to\infty} \mathrm{I}_{d,n,3}^*=0$.
\ Further, using now the calculations for \ $\mathrm{I}_{d,n,1}^{**}$, \ we have
 \[
    \mathrm{I}_{d,n,3}^{**} \leq x^{-\beta} \EE(X_0^\beta)\frac{n^{\gamma_2-\frac{\beta}{\alpha}}}{(L(n))^\beta}  \to 0
     \qquad \text{as \ $n\to\infty$,}
 \]
 with the same choice of \ $\gamma_2$ \ and \ $\beta$ \ that were chosen in the calculations for \ $\mathrm{I}_{d,n,1}^{**}$,
 \ i.e. with \ $\gamma_2\in(\frac{1}{2},\frac{\beta}{\alpha})$ \ and \ $\beta\in(\frac{\alpha}{2},\alpha)$.
\ This yields \ $\lim_{d\to\infty}\limsup_{n\to\infty} \mathrm{I}_{d,n,3}^{**}=0$, and
 putting parts together, we have $\lim_{d\to\infty}\limsup_{n\to\infty} \mathrm{I}_{d,n,3}=0$.

We turn to check \ $\lim_{d\to\infty}\limsup_{n\to\infty} \mathrm{I}_{d,n,4}=0$.
\ For each \ $d\in\NN$ \ and large enough \ $n\in\NN$ \ satisfying \ $m_n\geq d+1$, \ by \eqref{help_GWI_rep}, we have
 \begin{align*}
  \mathrm{I}_{d,n,4}
    &\leq x\frac{n}{a_n} \sum_{j=d+1}^{m_n} \EE\Big( X_0 \bone_{\{X_0\leq a_nx\}} \bone_{\{ \Pi_j^{(0)}\circ X_0 > a_nx/2\}}\Big)
         + x\frac{n}{a_n} \sum_{j=d+1}^{m_n} \EE\Big( X_0 \bone_{\{X_0\leq a_nx\}} \bone_{\{ \kappa_j > a_nx/2\}}\Big)\\
    &=:\mathrm{I}_{d,n,4}^* + \mathrm{I}_{d,n,4}^{**}.
 \end{align*}
Here
 \begin{align*}
   \mathrm{I}_{d,n,4}^*
     = x \frac{n}{a_n} \sum_{j=d+1}^{m_n} \EE\Big( X_0 \bone_{\{X_0\leq a_nx\}} \EE(\bone_{\{ \Pi_j^{(0)}\circ X_0 > a_nx/2\}} \mid X_0 )  \Big),
 \end{align*}
 where, by the conditional Markov's inequality (see, e.g., Klenke \cite[Exercise 8.2.5]{Kle})
  and the independence of \ $\Pi_j^{(0)}\circ i$ \ and \ $X_0$ \ for each \ $i\in\ZZ_+$ \ and \ $j\in\NN$, \ we have
 \begin{align*}
     \EE(\bone_{\{ \Pi_j^{(0)}\circ X_0 > a_nx/2\}} \mid X_0 )
      & = \PP(\Pi_j^{(0)}\circ X_0 > a_nx/2 \mid X_0)
         \leq \frac{\EE(\Pi_j^{(0)}\circ X_0 \mid X_0)}{a_nx/2} \\
     &= \frac{  (\EE(\Pi_j^{(0)}\circ k))\vert_{k=X_0}}{a_nx/2}
       =\frac{2}{a_nx}m_\xi^j X_0.
 \end{align*}
So, by \eqref{help_a_n_x} and Karamata's theorem,
 \begin{align*}
 \mathrm{I}_{d,n,4}^*
   &\leq x \frac{n}{a_n} \sum_{j=d+1}^{m_n} \EE\Big( X_0 \bone_{\{X_0\leq a_nx\}} \frac{2}{a_nx}m_\xi^j X_0  \Big)
   = \frac{2n}{a_n^2} \EE( X_0^2 \bone_{\{X_0\leq a_nx\}} ) \sum_{j=d+1}^{m_n} m_\xi^j\\
   &\leq 2x^2\cdot n\PP(X_0>a_nx) \cdot \frac{\EE( X_0^2 \bone_{\{X_0\leq a_nx\}} )}{(a_nx)^2 \PP(X_0>a_nx)}\cdot \frac{m_\xi^{d+1}}{1-m_\xi}
   \to 2x^2 \cdot x^{-\alpha} \cdot \frac{\alpha}{2-\alpha}  \cdot \frac{m_\xi^{d+1}}{1-m_\xi}
 \end{align*}
 as \ $n\to\infty$, \ yielding \ $\lim_{d\to\infty}\limsup_{n\to\infty} \mathrm{I}_{d,n,4}^{*}=0$.
\ Further, by \eqref{help_a_n_x} and the independence of \ $X_0$ \ and \ $\kappa_j$, \ $j\in\NN$, \ we have
 \begin{align*}
   \mathrm{I}_{d,n,4}^{**} &= x\frac{n}{a_n} \sum_{j=d+1}^{m_n} \EE( X_0 \bone_{\{X_0\leq a_nx\}}) \PP(\kappa_j > a_nx/2)
                            \leq x\frac{n}{a_n} \sum_{j=d+1}^{m_n} \EE( X_0 \bone_{\{X_0\leq a_nx\}}) \PP(\kappa_{j,\infty} > a_nx/2)\\
                           & = x \frac{n}{a_n} \sum_{j=d+1}^{m_n} \EE( X_0 \bone_{\{X_0\leq a_nx\}}) \PP( X_0 > a_nx/2)
                             \leq x \frac{nm_n}{a_n} \EE( X_0 \bone_{\{X_0\leq a_nx\}}) \PP( X_0 > a_nx/2)\\
                           &=x^2\cdot \frac{m_n}{n} \cdot \frac{\EE( X_0 \bone_{\{X_0\leq a_nx\}})}{a_nx \PP(X_0>a_nx)}\cdot n\PP(X_0>a_nx) \cdot n\PP(X_0>a_nx/2)\\
                           &\to x^2\cdot 0\cdot \frac{\alpha}{1-\alpha} \cdot x^{-\alpha} \cdot (x/2)^{-\alpha} =0
                           \qquad \text{as \ $n\to\infty$,}
 \end{align*}
 yielding $\lim_{d\to\infty}\limsup_{n\to\infty} \mathrm{I}_{d,n,4}^{**}=0$,
 and putting parts together, $\lim_{d\to\infty}\limsup_{n\to\infty} \mathrm{I}_{d,n,4}=0$.

One can realize that in Step 2/(a) we used the assumption \ $\alpha\in(0,1)$ \ only in handling the term \ $\mathrm{I}_{d,n,2}$,
 \ the presented technique for \ $\mathrm{I}_{d,n,1}$, \ $\mathrm{I}_{d,n,3}$ \ and \ $\mathrm{I}_{d,n,4}$ \ would work in case of
 \ $\alpha\in(1,\frac{4}{3})$ \ as well.

{\sl Step 2/(b).} We check \eqref{BJMM2}.
We can proceed similarly as in case of \eqref{BJMM1}.
For each \ $j\in\ZZ_+$, \ let us introduce the notation \ $\tX_j := X_j - \EE(X_j) =  X_j - \EE(X_0)$.
\ Let \ $x\in\RR_{++}$ \ be fixed.
For each \ $d\in\NN$ \ and large enough \ $n\in\NN$ \ satisfying \ $m_n\geq d+1$, \ we have
 \begin{align*}
  & n \sum_{j=d+1}^{m_n} \EE\big( (\vert a_n^{-1} \tX_j\vert \wedge x) (\vert a_n^{-1} \tX_0\vert \wedge x) \big)  \\
  & \qquad = x^2 n \sum_{j=d+1}^{m_n} \PP( \vert \tX_j\vert > a_n x, \vert \tX_0\vert > a_n x )
             + n \sum_{j=d+1}^{m_n} \EE\big( \vert a_n^{-1} \tX_j\vert \vert a_n^{-1} \tX_0\vert \bone_{\{ \vert \tX_j\vert \leq a_n x,  \vert \tX_0\vert \leq a_n x  \}}\big) \\
  & \phantom{\qquad =\,} + x n \sum_{j=d+1}^{m_n} \EE\big( \vert a_n^{-1} \tX_j\vert \bone_{\{ \vert \tX_j\vert \leq a_n x,  \vert \tX_0\vert > a_n x  \}}\big)
                         + x n \sum_{j=d+1}^{m_n} \EE\big( \vert a_n^{-1} \tX_0\vert \bone_{\{ \vert \tX_j\vert > a_n x,  \vert \tX_0\vert \leq a_n x  \}}\big) =:
 \end{align*}
 \begin{align*}
  & \qquad =: \mathrm{\tI}_{d,n,1} + \mathrm{\tI}_{d,n,2} + \mathrm{\tI}_{d,n,3} + \mathrm{\tI}_{d,n,4}.
 \end{align*}
We check that \ $\lim_{d\to\infty} \limsup_{n\to\infty} \mathrm{\tI}_{d,n,i} = 0$, \ $i=1,2,3,4$, \ yielding \eqref{BJMM2}.
Since \ $a_n\to\infty$ \ as \ $n\to\infty$, \ for large enough \ $n\in\NN$, \ we have \ $\frac{a_n x}{2}>\EE(X_0)$, \
 and hence for large enough \ $n\in\NN$, \
 \begin{align*}
 \mathrm{\tI}_{d,n,1}
   &\leq x^2 n \sum_{j=d+1}^{m_n} \PP\Big( (\{X_j > a_nx/2\}\cup\{\EE(X_0) > a_nx/2\}) \cap (\{X_0 > a_nx/2\}\cup\{\EE(X_0) > a_nx/2\}) \Big)\\
   &= x^2 n \sum_{j=d+1}^{m_n} \PP(X_j>a_nx/2, X_0>a_nx/2)\\
   &= x^2 n \PP(X_0>a_nx/2) \sum_{j=d+1}^{m_n} \PP(X_j>a_nx/2 \mid X_0>a_nx/2),
 \end{align*}
 which coincides with \ $\mathrm{I}_{d,n,1}$ \ in Step 2/(a) replacing \ $a_n$ \ by \ $a_n/2$.
\ Since, as we noted at the end of Step 2/(a), our technique for \ $\mathrm{I}_{d,n,1}$ \ works in case of \ $\alpha \in (1,\frac{4}{3})$ \ as well,
 we have \ $\lim_{d\to\infty}\limsup_{n\to\infty} \mathrm{\tI}_{d,n,1}=0$.

We turn to check \ $\lim_{d\to\infty}\limsup_{n\to\infty} \mathrm{\tI}_{d,n,2}=0$.
\ As we will see, here we will effectively use the assumption \ $\alpha<\frac{4}{3}$, \ meaning that our technique does not work for \ $\alpha\geq \frac{4}{3}$.
\ For each \ $d\in\NN$ \ and large enough \ $n\in\NN$ \ satisfying \ $m_n\geq d+1$, \ by the law of total expectation, we have
 \begin{align*}
  \mathrm{\tI}_{d,n,2}
   = n \sum_{j=d+1}^{m_n} \sum_{\{\ell\in\ZZ_+ : \vert \ell - \EE(X_0)\vert  \leq a_nx\}} \vert a_n^{-1}(\ell - \EE(X_0))\vert \PP(X_0=\ell)
                     \EE( \vert a_n^{-1} \tX_j\vert \bone_{\{ \vert \tX_j\vert   \leq a_nx\}} \mid   X_0=\ell).
 \end{align*}
The function \ $\tg:\ZZ_+\to\RR$, \ $\tg(y):=\frac{\vert y - \EE(X_0)\vert}{(a_nx)^{1-p}} \bone_{\{\vert y-\EE(X_0)\vert  \leq a_nx\}}$, \ $y\in\ZZ_+$, \
 satisfies \ $\tg(y)\leq V(y) = 1 + y^p$, \ $y\in\ZZ_+$, \ for any \ $p\in(0,\alpha\wedge 1) = (0,1)$ \ for large enough \ $n\in\NN$.
\ Indeed, given \ $p\in(0,1)$, \ we have \ $\tg(0) = \frac{\EE(X_0)}{(a_nx)^{1-p}}<1=V(0)$ \ for large enough \ $n\in\NN$,
 \ $\tg$ \ is monotone decreasing on \ $[0,\EE(X_0)]$, \ $\tg(\EE(X_0))=0$, \ $\tg$ \ is monotone increasing on \ $[\EE(X_0), a_nx+\EE(X_0)]$,
 \ and \ $\tg(a_nx+\EE(X_0))=(a_nx)^p<V(a_nx+\EE(X_0))$.
\ So, by the \ $V$-uniformly ergodicity of \ $(X_k)_{k\in\ZZ_+}$ \ (see Corollary \ref{Cor_ergod}), we have
 \[
    \left\vert  \EE\left( \frac{\vert \tX_j\vert}{(a_nx)^{1-p}} \bone_{\{\vert \tX_j\vert  \leq a_nx\}} \,\Big \vert\,   X_0=\ell\right)
                - \EE\left( \frac{\vert \tX_0\vert}{(a_nx)^{1-p}} \bone_{\{\vert \tX_0\vert  \leq a_nx\}} \right)  \right\vert
    \leq C (1+\ell^p)\varrho^j,
    \quad j\in\ZZ_+, \ell\in\ZZ_+,
 \]
 with some \ $C>0$ \ and \ $\varrho\in(0,1)$.
\ So
 \begin{align*}
  \mathrm{\tI}_{d,n,2}
   &= \frac{n}{a_n^2} \sum_{j=d+1}^{m_n} \sum_{\{\ell\in\ZZ_+ : \vert \ell - \EE(X_0)\vert \leq a_nx\}} \vert \ell - \EE(X_0)\vert \PP(X_0=\ell)
                        (a_nx)^{1-p}\EE\left(\frac{\vert \tX_j\vert}{(a_nx)^{1-p}} \bone_{\{\vert \tX_j\vert  \leq a_nx\}} \,\Big \vert\,  X_0=\ell\right)\\
   &\leq \frac{n}{a_n^2} \sum_{j=d+1}^{m_n} \sum_{\{\ell\in\ZZ_+ : \vert \ell - \EE(X_0)\vert  \leq a_nx\}} \vert \ell - \EE(X_0)\vert \PP(X_0=\ell)
                         (a_nx)^{1-p} \\
  &\phantom{\leq \frac{n}{a_n^2} \sum_{j=d+1}^{m_n} \sum_{\{\ell\in\ZZ_+ : \vert \ell - \EE(X_0)\vert <a_nx\}}}
                 \times \left(  \EE\left( \frac{\vert \tX_0\vert}{(a_nx)^{1-p}} \bone_{\{\vert \tX_0\vert  \leq a_nx\}} \right)  + C(1+\ell^p)\varrho^j \right).
 \end{align*}
Here the term corresponding to \ $\ell=0$ \ in \ $C(1+\ell^p)\varrho^j$ \ tends to 0 as \ $n\to\infty$, \ namely,
 \begin{align*}
   &\frac{n}{a_n^2} \sum_{j=d+1}^{m_n} \EE(X_0) \PP(X_0=0) a_n^{1-p}C\varrho^j
     \leq C\EE(X_0) \PP(X_0=0) \cdot n a_n^{-1-p} \cdot \frac{\varrho^{d+1}}{1-\varrho}\\
   &\qquad =C\EE(X_0) \PP(X_0=0) \cdot n^{1-\frac{1+p}{\alpha}}(L(n))^{-1-p} \cdot \frac{\varrho^{d+1}}{1-\varrho}
     \to 0 \qquad \text{as \ $n\to\infty$,}
 \end{align*}
 since one can choose \ $p\in (0,\alpha\wedge 1)=(0,1)$ \ in a way that \ $\alpha<p+1$, \ i.e., \ $\alpha-1<p<1$ \ should be satisfied.
Consequently,
 \begin{align*}
  \mathrm{\tI}_{d,n,2}
   &\leq \frac{n}{a_n^2} \sum_{j=d+1}^{m_n}
                       \Bigg[(\EE(\vert \tX_0 \vert \bone_{\{ \vert \tX_0\vert \leq a_nx\}}))^2  + 2C\varrho^j (a_nx)^{1-p} \EE( \vert \tX_0\vert X_0^p
                       \bone_{\{ \vert X_0\vert  \leq a_nx\}})\Bigg]\\
   &\leq \frac{nm_n}{a_n^2}(\EE( \vert \tX_0\vert \bone_{\{ \vert \tX_0\vert \leq a_nx\}}))^2
                + \frac{n}{a_n^2} 2C (a_nx)^{1-p}\EE( \vert \tX_0\vert X_0^p \bone_{\{ \vert \tX_0\vert  \leq a_nx\}}) \sum_{j=d+1}^\infty \varrho^j\\
   &= \frac{nm_n}{a_n^2}(\EE( \vert \tX_0\vert \bone_{\{ \vert \tX_0\vert  \leq a_nx\}}))^2
                + \frac{n}{a_n^2} 2C (a_nx)^{1-p}\EE( \vert \tX_0\vert X_0^p \bone_{\{ \vert \tX_0\vert  \leq a_nx\}}) \frac{\varrho^{d+1}}{1-\varrho} \\
   &=: \mathrm{\tI}_{d,n,2}^* +  \mathrm{\tI}_{d,n,2}^{**}.
 \end{align*}
Here \ $\lim_{d\to\infty}\limsup_{n\to\infty} \mathrm{\tI}_{d,n,2}^*=0$, \ since, due to the assumption \ $\alpha\in(1,\frac{4}{3})$, \ we have
 \begin{align*}
   \EE( \vert \tX_0\vert \bone_{\{ \vert \tX_0\vert <a_nx\}})
      \to  \EE( \vert \tX_0\vert )\in\RR_+ \qquad \text{as \ $n\to\infty$,}
 \end{align*}
 and, since \ $m_n = \lfloor n^{\gamma_2}\rfloor$ \ with \ $\gamma_2\in(\frac{1}{2}, \frac{1}{\alpha}\wedge 1) = (\frac{1}{2}, \frac{1}{\alpha})$,
 \ $\frac{nm_n}{a_n^2}\leq n^{1+\gamma_2-\frac{2}{\alpha}}(L(n))^{-2}\to 0 \qquad \text{as \ $n\to\infty$,}$
 \ provided that \ $\gamma_2$ \ is chosen such that \ $\gamma_2<\frac{2}{\alpha}-1$.
\ One can choose such a \ $\gamma_2$, \ since, under the assumption \ $\alpha\in(1,\frac{4}{3})$ \ both conditions
 \ $\gamma_2\in(\frac{1}{2},\frac{1}{\alpha}\wedge 1) = (\frac{1}{2},\frac{1}{\alpha})$ \ and \ $\gamma_2<\frac{2}{\alpha}-1$ \ can be fulfilled
 following from \ $\frac{1}{2} < \frac{2}{\alpha}-1 < \frac{1}{\alpha} \Longleftrightarrow \alpha\in(1,\frac{4}{3})$.
\ Briefly, we need to suppose \ $\alpha\in(1,\frac{4}{3})$ \ to be able to derive \ $\frac{nm_n}{a_n^2}\to 0$ \ as \ $n\to\infty$.
\ We note that such a condition also appears in the forthcoming book of Mikosch and Wintenberger \cite[page 324]{MikWin2}.
\ Further, using that
 \[
   X_0^p \leq \vert X_0 - \EE(X_0)\vert^p + (\EE(X_0))^p =  \vert \tX_0\vert^p + (\EE(X_0))^p \qquad \text{for \ $p\in(\alpha-1,1)$,}
 \]
 by Karamata's theorem, we have
 \begin{align*}
  \mathrm{\tI}_{d,n,2}^{**}
     &\leq 2C x^{1-p} \frac{\varrho^{d+1}}{1-\varrho} na_n^{-1-p} \EE(\vert \tX_0\vert^{p+1}\bone_{\{\vert \tX_0\vert \leq  a_nx\}} )\\
     &\phantom{=\,} + 2C x^{1-p} \frac{\varrho^{d+1}}{1-\varrho} na_n^{-1-p} \EE(\vert \tX_0\vert (\EE(X_0))^p \bone_{\{\vert \tX_0\vert \leq a_nx\}} )\\
     &= 2C x^2 \frac{\varrho^{d+1}}{1-\varrho} \cdot n \PP(\vert \tX_0\vert > a_nx)
             \cdot \frac{\EE(\vert \tX_0\vert^{p+1}\bone_{\{\vert \tX_0\vert \leq a_nx\}} )}{(a_nx)^{p+1} \PP(\vert \tX_0\vert > a_nx)}\\
     &\phantom{=\,} + 2C x^{1-p} \frac{\varrho^{d+1}}{1-\varrho} \cdot na_n^{-1-p} \cdot (\EE(X_0))^p \cdot \EE(\vert \tX_0\vert \bone_{\{\vert \tX_0\vert \leq a_nx\}} )
  \end{align*}
  \begin{align*}
     &\to 2C x^2 \frac{\varrho^{d+1}}{1-\varrho} \cdot x^{-\alpha} \cdot \frac{\alpha}{p+1-\alpha}
          + 2C x^{1-p} \frac{\varrho^{d+1}}{1-\varrho}  \cdot 0 \cdot (\EE(X_0))^p \cdot \EE(\vert \tX_0\vert)
        =0 \qquad \text{as \ $n\to\infty$,}
 \end{align*}
 where we used \ $a_n^{-1-p} = n^{1-\frac{1+p}{\alpha}}(L(n))^{-1-p}\to 0$ \ as \ $n\to\infty$.
\ This yields \ $\lim_{d\to\infty}\limsup_{n\to\infty}\mathrm{\tI}_{d,n,2}^{**}=0$.

\noindent Putting parts together, we have \ $\lim_{d\to\infty}\limsup_{n\to\infty} \mathrm{\tI}_{d,n,2}=0$.
\ For historical fidelity, we note that for handling the term \ $\mathrm{\tI}_{d,n,2}$ \ we used the same idea as in the proof
 of Lemma 6 in Basrak and Kevei \cite{BasKev}.

We turn to check \ $\lim_{d\to\infty}\limsup_{n\to\infty} \mathrm{\tI}_{d,n,3}=0$.
\ Since \ $a_n\to\infty$ \ as \ $n\to\infty$, \ for large enough \ $n\in\NN$, \ we have \ $\frac{a_nx}{2}>\EE(X_0)$,
 \ and hence \ $\{ \vert \tX_0\vert > a_nx\} \leq \{ X_0> a_nx/2\} \cup\{ \EE(X_0) > a_nx/2\} = \{ X_0> a_nx/2\}$.
\ Further, for all \ $\beta\in(0,\alpha\wedge 1) =(0,1)$, \ we have
 \[
    \vert a_n^{-1} \tX_j\vert \bone_{\{ \vert \tX_j\vert \leq a_n x \}}
    \leq x\left(\frac{\vert \tX_j\vert}{a_nx}\right)^\beta
    \leq x^{1-\beta}\frac{X_j^\beta}{a_n^\beta} + x^{1-\beta}\frac{(\EE(X_0))^\beta}{a_n^\beta}.
 \]
Hence for large enough \ $n\in\NN$ \ and any \ $\beta\in(0,\alpha\wedge 1)=(0,1)$, \ we have
 \begin{align*}
  \mathrm{\tI}_{d,n,3}
   &\leq  x^{2-\beta} n \sum_{j=d+1}^{m_n} \EE\left( \left( \frac{X_j^\beta}{a_n^\beta} + \frac{(\EE(X_0))^\beta}{a_n^\beta} \right)\bone_{\{ X_0> a_nx/2\}} \right)\\
   &= x^2 2^{-\beta} \cdot n\PP(X_0 > a_nx/2) \sum_{j=d+1}^{m_n} \frac{\EE(X_j^\beta \bone_{\{X_0 > a_nx/2\}})}{(a_nx/2)^\beta \PP(X_0>a_nx/2)}\\
   &\phantom{=\;} + x^{2-\beta} (\EE(X_0))^\beta \frac{m_n-d}{a_n^\beta}  n\PP(X_0 > a_nx/2)\\
   &=: \mathrm{\tI}_{d,n,3}^* + \mathrm{\tI}_{d,n,3}^{**}.
 \end{align*}
Since, as we noted at the end of Step 2/(a), our technique for \ $\mathrm{I}_{d,n,3}$ \ works in case of \ $\alpha \in (1,\frac{4}{3})$ \ as well,
 we have $\lim_{d\to\infty}\limsup_{n\to\infty} \mathrm{\tI}_{d,n,3}^*=0$.
\ Further,
 \begin{align*}
    \mathrm{\tI}_{d,n,3}^{**}
      & = x^{2-\beta} (\EE(X_0))^\beta \frac{m_n-d}{a_n^\beta}  n \PP(X_0 > a_n)  \frac{\PP(X_0 > a_nx/2)}{\PP(X_0 > a_n)}\\
      & \leq  x^{2-\beta} (\EE(X_0))^\beta \cdot \frac{n^{\gamma_2}}{n^{\frac{\beta}{\alpha}} (L(n))^\beta} \cdot n \PP(X_0 > a_n) \cdot \frac{\PP(X_0 > a_nx/2)}{\PP(X_0 > a_n)}\\
      & \to x^{2-\beta} (\EE(X_0))^\beta \cdot 0 \cdot 1 \cdot (x/2)^{-\alpha} =0
      \qquad \text{as \ $n\to\infty$}
 \end{align*}
 with the same choice of \ $\gamma_2$ \ and \ $\beta$ \ that were chosen in the calculations for \ $\mathrm{I}_{d,n,1}^{**}$,
 \ i.e. with \ $\gamma_2\in(\frac{1}{2},\frac{\beta}{\alpha})$ \ and \ $\beta\in(\frac{\alpha}{2},1)$.
\ This yields $\lim_{d\to\infty}\limsup_{n\to\infty} \mathrm{\tI}_{d,n,3}^{**}=0$, and putting parts together
 we have $\lim_{d\to\infty}\limsup_{n\to\infty} \mathrm{\tI}_{d,n,3}=0$.

We turn to check \ $\lim_{d\to\infty}\limsup_{n\to\infty} \mathrm{\tI}_{d,n,4}=0$.
 \ For each \ $d\in\NN$ \ and large enough \ $n\in\NN$ \ satisfying \ $m_n\geq d+1$, \ by \eqref{help_GWI_rep}, we have
 \begin{align*}
  \mathrm{\tI}_{d,n,4}
    &=x\frac{n}{a_n} \sum_{j=d+1}^{m_n} \EE\Big( \vert \tX_0\vert \bone_{\{ \vert \tX_0\vert \leq a_nx\}} \bone_{\{ X_j > a_nx +\EE(X_0)\}} \Big)\\
    &\leq x\frac{n}{a_n} \sum_{j=d+1}^{m_n} \EE\Big( \vert \tX_0\vert \bone_{\{ \vert \tX_0\vert \leq a_nx\}} \bone_{\{ \Pi_j^{(0)}\circ X_0 > (a_nx +\EE(X_0))/2\}}\Big)\\
    &\phantom{\leq} + x\frac{n}{a_n} \sum_{j=d+1}^{m_n} \EE\Big( \vert \tX_0\vert \bone_{\{\vert \tX_0\vert \leq a_nx\}} \bone_{\{ \kappa_j > (a_nx +\EE(X_0))/2\}}\Big)\\
    &=:\mathrm{\tI}_{d,n,4}^* + \mathrm{\tI}_{d,n,4}^{**}.
 \end{align*}
Here
 \begin{align*}
   \mathrm{\tI}_{d,n,4}^*
     = x\frac{n}{a_n} \sum_{j=d+1}^{m_n} \EE\Big( \vert \tX_0 \vert \bone_{\{\vert \tX_0\vert\leq a_nx\}}
                                       \EE(\bone_{\{ \Pi_j^{(0)}\circ X_0 > (a_nx+\EE(X_0))/2\}} \mid X_0 )  \Big),
 \end{align*}
 where
 \begin{align*}
     &\EE(\bone_{\{ \Pi_j^{(0)}\circ X_0 > (a_nx+\EE(X_0))/2\}} \mid X_0 )
       = \PP(\Pi_j^{(0)}\circ X_0 > (a_nx+\EE(X_0))/2 \mid X_0)\\
     &\qquad  \leq \frac{  (\EE(\Pi_j^{(0)}\circ k))\vert_{k=X_0}}{(a_nx+\EE(X_0))/2}
       =\frac{2}{a_nx+\EE(X_0)}m_\xi^j X_0
       \leq \frac{2}{a_nx}m_\xi^j X_0.
 \end{align*}
So, by Karamata's theorem,
 \begin{align*}
 \mathrm{\tI}_{d,n,4}^*
   &\leq x\frac{n}{a_n} \sum_{j=d+1}^{m_n} \EE\Big( \vert \tX_0\vert \bone_{\{\vert \tX_0\vert \leq a_nx\}} \frac{2}{a_nx}m_\xi^j X_0  \Big)
   = \frac{2n}{a_n^2} \EE( \vert \tX_0\vert X_0 \bone_{\{\vert \tX_0\vert \leq a_nx\}} ) \sum_{j=d+1}^{m_n} m_\xi^j\\
   &\leq \frac{m_\xi^{d+1}}{1-m_\xi} \frac{2n}{a_n^2} \EE( \vert \tX_0\vert^2 \bone_{\{\vert \tX_0\vert \leq a_nx\}} )\\
   &\phantom{\leq\,} + \frac{m_\xi^{d+1}}{1-m_\xi} \EE(X_0)\frac{2n}{a_n^2} \EE( \vert \tX_0\vert \bone_{\{\vert \tX_0\vert \leq a_nx\}} ) \\
   &=2x^2\frac{m_\xi^{d+1}}{1-m_\xi} \cdot n\PP(\vert \tX_0\vert > a_nx)
        \cdot \frac{\EE( \vert \tX_0\vert^2 \bone_{\{\vert \tX_0\vert \leq a_nx\}})}{(a_nx)^2 \PP(\vert \tX_0\vert > a_nx)}\\
   &\phantom{=\,}
     + 2\frac{m_\xi^{d+1}}{1-m_\xi}\EE(X_0) \cdot n^{1-\frac{2}{\alpha}} (L(n))^{-2} \EE(\vert \tX_0\vert \bone_{\{\vert \tX_0\vert \leq a_nx\}})\\
   &\to 2x^2\frac{m_\xi^{d+1}}{1-m_\xi} \cdot1\cdot \frac{\alpha}{2-\alpha}
       + 2\frac{m_\xi^{d+1}}{1-m_\xi}\EE(X_0) \cdot 0\cdot \EE(\vert \tX_0\vert)
    = 2x^2\frac{m_\xi^{d+1}}{1-m_\xi} \frac{\alpha}{2-\alpha}
 \end{align*}
 as \ $n\to\infty$, \ yielding \ $\lim_{d\to\infty}\limsup_{n\to\infty} \mathrm{\tI}_{d,n,4}^{*}=0$.
\ Further, by the independence of \ $X_0$ \ and \ $\kappa_j$, \ $j\in\NN$, \ we have
 \begin{align*}
   \mathrm{\tI}_{d,n,4}^{**} &= x\frac{n}{a_n} \sum_{j=d+1}^{m_n} \EE( \vert \tX_0\vert \bone_{\{\vert \tX_0\vert \leq a_nx\}}) \PP(\kappa_j > (a_nx +\EE(X_0))/2)\\
                            &\leq x\frac{n}{a_n} \sum_{j=d+1}^{m_n} \EE( \vert \tX_0\vert \bone_{\{ \vert \tX_0\vert \leq a_nx\}}) \PP(\kappa_{j,\infty} > (a_nx +\EE(X_0))/2)\\
                           & = x\frac{n}{a_n} \sum_{j=d+1}^{m_n} \EE( \vert \tX_0\vert \bone_{\{\vert \tX_0\vert \leq a_nx\}}) \PP( X_0 > (a_nx +\EE(X_0))/2)\\
                           & \leq x\frac{nm_n}{a_n} \EE( \vert \tX_0\vert  \bone_{\{\vert \tX_0\vert \leq a_nx\}}) \PP( X_0 > (a_nx +\EE(X_0))/2)\\
                           &= x \cdot \frac{m_n}{a_n} \cdot \EE( \vert \tX_0\vert  \bone_{\{\vert \tX_0\vert \leq a_nx\}}) \cdot n \PP(X_0> a_n)
                               \cdot \frac{\PP(X_0> (a_nx+\EE(X_0))/2)}{\PP(X_0> a_n)} \\
                           &\to x\cdot 0\cdot \EE(\vert \tX_0\vert)\cdot 1\cdot 1 =0
                           \qquad \text{as \ $n\to\infty$,}
 \end{align*}
 since \ $\frac{m_n}{a_n}\leq n^{\gamma_2 -\frac{1}{\alpha}}(L(n))^{-1}\to 0$ \ as \ $n\to\infty$ \ due to \ $\gamma_2<\frac{1}{\alpha}\wedge 1 = \frac{1}{\alpha}$
 \ and \ $x/2 +\EE(X_0)/(2a_n)\in[x/2, x/2+1]$ \ for large enough \ $n\in\NN$ \ and one can use uniform convergence theorem for
 regularly varying functions ( see, e.g., Theorem 1.5.2 in Bingham et al.\ \cite{BinGolTeu}).
This yields \ $\lim_{d\to\infty}\limsup_{n\to\infty} \mathrm{\tI}_{d,n,4}^{**}=0$,
 \ and putting parts together, \ $\lim_{d\to\infty}\limsup_{n\to\infty} \mathrm{\tI}_{d,n,4}=0$.

{\sl Step 3 (checking conditions (iv) and (v) of Theorem \ref{BJMW1}).}
Next, we check that condition (iv) of Theorem \ref{BJMW1} holds and we determine the constants \ $c_+$ \ and \ $c_{-}$ \ explicitly as well.
First, note that for each \ $d\in\NN$ \ we have \ $b_{-}(d) = \lim_{n\to\infty} n \PP(S_d \leq -a_n) = 0$, \ since \ $S_d$ \ is non-negative
 and \ $a_n \in \RR_{++}$, \ $n \in \NN$, \ yielding that \ $c_{-} = 0$.
\ We show that
 \[
   c_{+} =  \frac{1-m_\xi^\alpha}{(1-m_\xi)^\alpha}.
 \]
For each \ $d \in \ZZ_+$, \ by part (ii) of Proposition D1 in Barczy et al.\ \cite{BarNedPap3}, we obtain
 \begin{align*}
  b_{+}(d+1)
  &= \lim_{n\to\infty} n \PP(X_0+\cdots+X_d>a_n)
   = \lim_{n\to\infty}
      n \PP(X_0 > a_n)
      \frac{\PP(X_0+\cdots+X_d>a_n)}{\PP(X_0>a_n)} \\
  &= \frac{1-m_\xi^\alpha}{(1-m_\xi)^\alpha} \biggl(\frac{ (1-m_\xi^{d+1})^\alpha }{1-m_\xi^\alpha}  + \sum_{j=1}^d (1-m_\xi^{d-j+1})^\alpha\biggr) .
 \end{align*}
Hence for each \ $d \in \NN$,
 \begin{align}\label{help_bd}
  b_+(d+1) - b_+(d)
    = \frac{1-m_\xi^\alpha}{(1-m_\xi)^\alpha}
      \biggl(\frac{ (1-m_\xi^{d+1})^\alpha  - (1-m_\xi^d)^\alpha}{1-m_\xi^\alpha}  + (1-m_\xi^d)^\alpha\biggr) ,
 \end{align}
 which tends to \ $(1-m_\xi^\alpha)/(1-m_\xi)^\alpha$ \ as  \ $d\to\infty$, \ as desired.
Note also that, in accordance with Lemma 3.1 in Mikosch and Wintenberger \cite{MikWin}, we have
 \[
   b_+(d+1) - b_+(d)
      = \EE \left( \left( \sum_{\ell=0}^d \Theta_\ell \right)^\alpha_+ -   \left( \sum_{\ell=1}^d \Theta_\ell \right)^\alpha_+ \right),
     \qquad d\in\ZZ_+,
 \]
 where \ $(\Theta_\ell)_{\ell\in\ZZ_+}$ \ is the (forward) spectral tail process of \ $(X_\ell)_{\ell\in\ZZ_+}$ \
 given by \ $\Theta_\ell = m_\xi^\ell$, \ $\ell \in \ZZ_+$ \ (see, e.g., Theorem \ref{Xtailprocess}), and \ $x_+ := x \lor 0$, \ $x \in \RR$.
\ Indeed, for each \ $d\in\ZZ_+$,
 \begin{align*}
   \EE \left( \left( \sum_{\ell=0}^d \Theta_\ell \right)^\alpha_+ -   \left( \sum_{\ell=1}^d \Theta_\ell \right)^\alpha_+ \right)
    & = \EE \left( \left( \sum_{\ell=0}^d m_\xi^\ell \right)^\alpha_+ -   \left( \sum_{\ell=1}^d m_\xi^\ell \right)^\alpha_+ \right) \\
    & = \left( \frac{1-m_\xi^{d+1}}{1-m_\xi} \right)^\alpha  - \left(m_\xi \frac{1-m_\xi^d}{1-m_\xi} \right)^\alpha ,
 \end{align*}
 yielding \eqref{help_bd}.

Condition (v) of Theorem \ref{BJMW1} holds trivially, since \ $\EE(X_1 - \EE(X_1)) = 0$ \
 in case of \ $\alpha\in(1, \frac{4}{3})$.

{\sl Step 4 (application of Theorem \ref{BJMW1}).}
All in all, by Theorem \ref{BJMW1}, in case of \ $\alpha\in(0,1)$, \  we have
 \begin{align}\label{help_Bur_1}
  \biggl(\frac{1}{a_n} \sum_{k=1}^\nt X_k\biggr)_{t\in\RR_+}
    \distrf (S_t)_{t\in\RR_+}  \qquad \text{as \ $n\to\infty$,}
 \end{align}
 and,  in case of \ $\alpha\in(1, \frac{4}{3})$, \ we have
 \begin{align}\label{help_Bur_2}
  \biggl(\frac{1}{a_n} \sum_{k=1}^\nt (X_k - \EE(X_k))\biggr)_{t\in\RR_+}
   \distrf (S_t)_{t\in\RR_+}  \qquad \text{as \ $n\to\infty$,}
 \end{align}
 where \ $(S_t)_{t\in\RR_+}$ \ is an \ $\alpha$-stable process such that the characteristic function of \ $S_1$ \ has the form
 \begin{align}\label{help_char_func}
 \begin{split}
  \EE(\ee^{\ii\vartheta\cS_1})
  &= \exp\Bigl\{- C_\alpha |\vartheta|^\alpha \Bigl((c_+ + c_-) - \ii(c_+ - c_-) \tan\Bigl(\frac{\pi\alpha}{2}\Bigr) \sign(\vartheta)\Bigr)\Bigr\} \\
   &= \exp\Bigl\{- C_\alpha \frac{1-m_\xi^\alpha}{(1-m_\xi)^\alpha}
                  |\vartheta|^\alpha \bigl(1 - \ii \tan\Bigl(\frac{\pi\alpha}{2}\Bigr) \sign(\vartheta)\bigr)\Bigr\}
  \end{split}
 \end{align}
 for \ $\theta\in\RR$, \ where \ $C_\alpha$, \ $\alpha\in (0,1)\cup (1,\frac{4}{3})$, \ is given in \eqref{c_alpha}.

{\sl Step 5 (limit theorems for truncated moments).}
For each \ $\alpha \in (0, 1)$ \ and \ $t \in \RR_{++}$, \ by Karamata's theorem (see, e.g., Lemma \ref{truncated_moments}),
 \begin{equation}\label{centering01T}
   \frac{\lfloor nt\rfloor}{a_n} \EE(X_0 \bbone_{\{X_0\leq a_n\}})
  = \frac{\lfloor nt\rfloor}{n}
     \frac{\EE(X_0 \bbone_{\{X_0\leq a_n\}})}{a_n\PP(X_0 > a_n)}
     n \PP(X_0> a_n)
  \to  \frac{\alpha}{1-\alpha} t
   \qquad \text{as \ $n \to \infty$.}
 \end{equation}
Consequently, the decomposition
 \[
   \frac{1}{a_n} \sum_{k=1}^\nt (X_k - \EE(X_k \bbone_{\{ X_k \leq a_n\}}))
   = \frac{1}{a_n} \sum_{k=1}^\nt X_k - \frac{\lfloor nt\rfloor}{a_n}\EE(X_0 \bbone_{\{X_0\leq a_n\}}) ,
   \qquad n \in \NN , \quad t \in \RR_+ ,
 \]
 \eqref{help_Bur_1} and Slutsky's lemma yield that
 \[
   \Biggl(\frac{1}{a_n} \sum_{k=1}^\nt (X_k - \EE(X_k \bbone_{\{ X_k \leq a_n\}}))\Biggr)_{t\in\RR_+}
   \distrf \biggl(S_t - \frac{\alpha}{1-\alpha} t\biggr)_{t\in\RR_+}
   \qquad \text{as \ $n\to\infty$.}
 \]
Further, for each \ $\alpha \in (1, \frac{4}{3})$ \ and \ $t \in \RR_{++}$, \ by Karamata's theorem (see, e.g., Lemma \ref{truncated_moments}),
 \begin{align*}
   \frac{\lfloor nt\rfloor}{a_n} \EE(X_0 \bbone_{\{X_0>a_n\}})
   \to \frac{\alpha}{\alpha-1} t
      \qquad \text{as \ $n \to \infty$.}
 \end{align*}
Consequently, the decomposition
 \[
   \frac{1}{a_n} \sum_{k=1}^\nt (X_k - \EE(X_k \bbone_{\{ X_k \leq a_n\}}))
   = \frac{1}{a_n} \sum_{k=1}^\nt (X_k - \EE(X_k)) + \frac{\lfloor nt\rfloor}{a_n}\EE(X_0 \bbone_{\{X_0>a_n\}}) ,
   \qquad n \in \NN , \quad t \in \RR_+ ,
 \]
 \eqref{help_Bur_2}, and Slutsky's lemma yield that
 \[
   \Biggl(\frac{1}{a_n} \sum_{k=1}^\nt (X_k - \EE(X_k \bbone_{\{X_k\leq a_n\}}))\Biggr)_{t\in\RR_+}
   \distrf \biggl(S_t - \frac{\alpha}{1-\alpha}t\biggr)_{t\in\RR_+}
   \qquad \text{as \ $n\to\infty$.}
 \]

{\sl Step 6 (characteristic functions for the limit processes).}
Note that \ $(S_t)_{t\in\RR_+}$ \ and \ $\bigl(\cZ^{(\alpha)}_t + \frac{\alpha}{1-\alpha}t\bigr)_{t\in\RR_+}$ \ coincide in law for each
 \ $\alpha \in (0,1)\cup (1,\frac{4}{3})$, \ and, by \eqref{help_char_func}, the second expression for
 \ $\EE\Bigl(\exp\Bigl\{\ii\vartheta\Bigl(\cZ_1^{(\alpha)}+\frac{\alpha}{1-\alpha}\Bigr)\Bigr\}\Bigr)$, \
 $\vartheta \in \RR$,
 \ readily follows for each \ $\alpha \in (0, 1) \cup (1,\frac{4}{3})$.
\ For the first expression for \ $\EE\Bigl(\exp\Bigl\{\ii\vartheta\Bigl(\cZ_1^{(\alpha)}+\frac{\alpha}{1-\alpha}\Bigr)\Bigr\}\Bigr)$, \ $\vartheta \in \RR$, \
 see (14.18) and (14.19) in Sato \cite{Sato}.
Consequently, in case of \ $\alpha \in(0, 1)$,
 \begin{align*}
  \EE(\ee^{\ii \vartheta \cZ^{(\alpha)}_1})
     &= \exp\Bigl\{ -\ii \vartheta \frac{\alpha}{1-\alpha} + \frac{1-m_\xi^\alpha}{(1-m_\xi)^\alpha}
                                                            \int_0^\infty (\ee^{\ii\vartheta u} - 1)\alpha u^{-1-\alpha}\,\dd u\Bigr\}\\
     &= \exp\Biggl\{ \frac{1-m_\xi^\alpha}{(1-m_\xi)^\alpha}
                      \int_0^\infty (\ee^{\ii\vartheta u} - 1 - \ii\vartheta u \bbone_{(0,1]}(u)  )\alpha u^{-1-\alpha}\,\dd u \\
     &\phantom{= \exp\Biggl\{ \;}
                   + \ii\vartheta  \frac{1-m_\xi^\alpha}{(1-m_\xi)^\alpha} \int_0^1 \alpha u^{-\alpha} \,\dd u  - \ii\vartheta \frac{\alpha}{1-\alpha} \Biggr\},
     \qquad \vartheta\in\RR,
 \end{align*}
 yielding the expression for \ $\EE\bigl(\ee^{\ii \vartheta \cZ^{(\alpha)}_1}\bigr)$, \ $\vartheta \in \RR$, \ in case of \ $\alpha \in (0, 1)$.
\ Similarly, in case of \ $\alpha\in(1, \frac{4}{3})$, \ we have
\begin{align*}
  \EE(\ee^{\ii \vartheta \cZ^{(\alpha)}_1})
     &= \exp\Bigl\{ -\ii \vartheta \frac{\alpha}{1-\alpha} + \frac{1-m_\xi^\alpha}{(1-m_\xi)^\alpha}
                                                            \int_0^\infty (\ee^{\ii\vartheta u} - 1 - \ii\vartheta u)\alpha u^{-1-\alpha}\,\dd u\Bigr\}\\
     &= \exp\Biggl\{ \frac{1-m_\xi^\alpha}{(1-m_\xi)^\alpha}
                      \int_0^\infty (\ee^{\ii\vartheta u} - 1 - \ii\vartheta u \bbone_{(0,1]}(u)  )\alpha u^{-1-\alpha}\,\dd u \\
     &\phantom{= \exp\Biggl\{ \;}
                   - \ii\vartheta  \frac{1-m_\xi^\alpha}{(1-m_\xi)^\alpha} \int_1^\infty \alpha u^{-\alpha} \,\dd u  - \ii\vartheta \frac{\alpha}{1-\alpha} \Biggr\},
     \qquad \vartheta\in\RR,
 \end{align*}
 yielding the expression for \ $\EE\bigl(\ee^{\ii \vartheta \cZ^{(\alpha)}_1}\bigr)$, \ $\vartheta \in \RR$, \ in case of \ $\alpha \in (1,\frac{4}{3})$, \ as well.
\proofend

\noindent{\bf Proof of Theorem \ref{iterated_aggr_2}.}
In case of \ $\alpha \in (0, 1)$, \ by \eqref{centering01T}, for each \ $t\in\RR_+$  we have
 \[
   \lim_{n\to\infty} \frac{\nt N}{a_nN^{\frac{1}{\alpha}}} \EE(X_0 \bbone_{\{X_0\leq a_n\}})
   = t \frac{\alpha}{1-\alpha} N^{1-\frac{1}{\alpha}}
   \to 0 \qquad \text{as \ $N \to \infty$,}
 \]
 hence, by Slutsky's lemma, \eqref{iterated_aggr_2_1} will be a consequence of \eqref{iterated_aggr_2_2}.

For each \ $N \in \NN$, \ by Theorem \ref{aggr_time} and by the continuity theorem, we obtain
 \[
   \biggl(\frac{1}{a_nN^{\frac{1}{\alpha}}}
         \sum_{k=1}^\nt \sum_{j=1}^N X^{(j)}_k\biggr)_{t\in\RR_+}
   \distrf \biggl(\frac{1}{N^{\frac{1}{\alpha}}}
                  \sum_{j=1}^N
                   \Bigl(\cZ^{(j,\alpha)}_t + \frac{\alpha}{1-\alpha} t\Bigr)\biggr)_{t\in\RR_+}
                   \qquad \text{as \ $n\to\infty$}
 \]
 in case of \ $\alpha \in (0, 1)$, \ and
 \[
   \biggl(\frac{1}{a_nN^{\frac{1}{\alpha}}}
         \sum_{k=1}^\nt \sum_{j=1}^N \bigl(X^{(j)}_k - \EE(X^{(j)}_k) \bigr)\biggr)_{t\in\RR_+}
   \distrf \biggl(\frac{1}{N^{\frac{1}{\alpha}}}
                  \sum_{j=1}^N
                   \Bigl(\cZ^{(j,\alpha)}_t + \frac{\alpha}{1-\alpha} t\Bigr)\biggr)_{t\in\RR_+}
                   \qquad \text{as \ $n\to\infty$}
 \]
 in case of \ $\alpha \in (1, \frac{4}{3})$, \ where \ $(\cZ^{(j,\alpha)}_t)_{t\in\RR_+}$, \ $j \in \NN$, \ are independent copies of \ $(\cZ^{(\alpha)}_t)_{t\in\RR_+}$
 given in Theorem \ref{aggr_time}.
Consequently, in order to prove \eqref{iterated_aggr_2_2} and \eqref{iterated_aggr_2_4}, we need to show that for each
 \ $\alpha \in (0, 1) \cup (1, \frac{4}{3})$, \ we have
 \begin{equation}\label{conv_Z_2}
  \biggl(\frac{1}{N^{\frac{1}{\alpha}}}
                 \sum_{j=1}^N
                  \biggl(\cZ^{(j,\alpha)}_t + \frac{\alpha}{1-\alpha} t\biggr)\biggr)_{t\in\RR_+}
  \distrf \Bigl(\cZ^{(\alpha)}_t + \frac{\alpha}{1-\alpha} t\Bigr)_{t\in\RR_+} \qquad \text{as \ $N \to \infty$.}
 \end{equation}
Since the processes \ $(\cZ^{(\alpha)}_t + \frac{\alpha}{1-\alpha} t)_{t\in\RR_+}$ \ and
 \ $(\cZ^{(j,\alpha)}_t + \frac{\alpha}{1-\alpha} t)_{t\in\RR_+}$, \ $j\in\NN$, \ have independent and stationary increments (being L\'evy processes)
 and \ $(\cZ^{(j,\alpha)}_t)_{t\in\RR_+}$, \ $j\in\NN$, \ are independent,
 in order to check \eqref{conv_Z_2}, it is enough to show that for each \ $t \in \RR_+$, \ we have
 \begin{equation}\label{conv_Z_2t}
  \frac{1}{N^{\frac{1}{\alpha}}}
                 \sum_{j=1}^N
                  \Bigl(\cZ^{(j,\alpha)}_t + \frac{\alpha}{1-\alpha} t\Bigr)
  \distr \cZ^{(\alpha)}_t + \frac{\alpha}{1-\alpha} t \qquad \text{as \ $N \to \infty$.}
 \end{equation}
In fact, for each \ $N \in \NN$ \ and \ $t \in \RR_+$, \ we have
 \[
   \frac{1}{N^{\frac{1}{\alpha}}}
                 \sum_{j=1}^N
                  \Bigl(\cZ^{(j,\alpha)}_t + \frac{\alpha}{1-\alpha} t\Bigr)
  \distre \cZ^{(\alpha)}_t + \frac{\alpha}{1-\alpha} t ,
 \]
 since, by Remark \ref{cZ_properties}, in case of \ $\alpha\in(0,1)\cup (1, \frac{4}{3})$,
 \ the distribution of \ $\cZ^{(\alpha)}_t + \frac{\alpha}{1-\alpha} t$ \ is strictly \ $\alpha$-stable, hence we obtain \eqref{conv_Z_2t},
  and hence \eqref{conv_Z_2}, thus we completed the proofs of \eqref{iterated_aggr_2_1},
 \eqref{iterated_aggr_2_2} and \eqref{iterated_aggr_2_4}.
\proofend

\vspace*{5mm}

\appendix

\vspace*{5mm}

\noindent{\bf\Large Appendices}

\section{The underlying space and vague convergence}
\label{App_vague}

For each \ $d \in \NN$, \ put \ $\RR_0^d := \RR^d \setminus \{\bzero\}$, \ and denote by
 \ $\cB(\RR_0^d)$ \ the Borel \ $\sigma$-algebra of \ $\RR_0^d$ \ induced by the metric \ $\varrho : \RR^d_0 \times \RR^d_0 \to \RR_+$ \ given by
 \begin{align}\label{S_metric}
  \varrho(\bx, \by) := \min\{\|\bx- \by\|, 1\} + \biggl|\frac{1}{\|\bx\|} - \frac{1}{\|\by\|}\biggr| , \qquad \bx, \by \in \RR^d_0 .
 \end{align}

The proof of the following lemma can be found in Barczy et al. \cite[Lemma B.1]{BarNedPap3_a}.

\begin{Lem}\label{Lem_S_top}
The set \ $\RR^d_0$ \ furnished with the metric \ $\varrho$ \ given in \eqref{S_metric} is a complete separable metric space, and \ $B \subset \RR^d_0$ \ is bounded with respect to the metric \ $\varrho$ \ if and only if \ $B$ \ is separated from the origin \ $\bzero \in \RR^d$, \ i.e., there exists \ $\vare \in \RR_{++}$ \ such that \ $B \subset \{\bx \in \RR^d_0 : \|\bx\| > \vare\}$.
Moreover, the topology and the Borel \ $\sigma$-algebra \ $\cB(\RR_0^d )$ \ on \ $\RR_0^d$ \ induced by the metric \ $\varrho$ \ coincides with
 the topology and the Borel \ $\sigma$-algebra on \ $\RR_0^d$ \ induced by the usual metric \ $d(\bx, \by) := \|\bx-\by\|$, \ $\bx, \by \in \RR_0^d$,
 \ respectively.
\end{Lem}

Since \ $\RR^d_0$ \ is locally compact, second countable and Hausdorff, one could choose a metric such that the relatively compact sets
 are precisely the bounded ones, see Kallenberg \cite[page 18]{kallenberg:2017}.
The metric \ $\varrho$ \ does not have this property, but we do not need it.

A measure \ $\nu$ \ on \ $(\RR_0^d, \cB(\RR_0^d))$ \ is said to be locally finite if \ $\nu(B) < \infty$ \ for every bounded Borel set \ $B$ \ with respect to the metric
 \ $\varrho$ \ given in \eqref{S_metric}, and write \ $\cM(\RR^d_0)$ \ for the class of locally finite measures on \ $(\RR_0^d, \cB(\RR_0^d))$.

The vague topology on \ $\cM(\RR^d_0)$ \ is constructed as in Chapter 4 in Kallenberg \cite{kallenberg:2017}.
The associated notion of vague convergence of a sequence \ $(\nu_n)_{n\in\NN}$ \ in \ $\cM(\RR^d_0)$ \ towards $\nu\in\cM(\RR^d_0)$,
 denoted by \ $\nu_n \distrv \nu$ \ as \ $n\to\infty$, \ is defined by the condition \ $\nu_n(f) \to \nu(f)$ \ as \ $n\to\infty$
 \ for every bounded, continuous function \ $f : \RR^d_0 \to \RR_+$ \ with bounded support, where \ $\kappa(f) := \int_{\RR^d_0} f(\bx) \, \kappa(\dd\bx)$ \ for each \ $\kappa \in \cM(\RR^d_0)$.

\section{Regularly varying distributions}
\label{App_reg_var}

First, we recall the notions of slowly varying and regularly varying functions, respectively.

\begin{Def}
A measurable function \ $U: \RR_{++} \to \RR_{++}$ \ is called regularly varying at infinity with
 index \ $\rho \in \RR$ \ if for all \ $c \in \RR_{++}$,
 \[
   \lim_{x\to\infty} \frac{U(cx)}{U(x)} = c^\rho .
 \]
In case of \ $\rho = 0$, \ we call \ $U$ \ slowly varying at infinity.
\end{Def}

\begin{Def}
A random variable \ $Y$ \ is called regularly varying with index \ $\alpha \in \RR_{++}$ \ if \ $\PP(|Y| > x) \in \RR_{++}$ \ for all \ $x \in \RR_{++}$, \ the function \ $\RR_{++} \ni x \mapsto \PP(|Y| > x) \in \RR_{++}$ \ is regularly varying at infinity with index \ $-\alpha$, \ and a tail-balance condition holds:
 \begin{equation}\label{TB}
   \lim_{x\to\infty} \frac{\PP(Y>x)}{\PP(|Y|>x)} = p , \qquad
   \lim_{x\to\infty} \frac{\PP(Y\leq-x)}{\PP(|Y|>x)} = q ,
 \end{equation}
 where \ $p + q = 1$.
\end{Def}

\begin{Rem}\label{tail-balance}
In the tail-balance condition \eqref{TB}, the second convergence can be replaced by
 \begin{equation}\label{tail-balance+}
  \lim_{x\to\infty} \frac{\PP(Y<-x)}{\PP(|Y|>x)} = q ,
 \end{equation}
 see, e.g., Barczy et al. \cite[Remark C.3]{BarNedPap3_a}.
\proofend
\end{Rem}

\begin{Lem}\label{hom}
\renewcommand{\labelenumi}{{\rm(\roman{enumi})}}
 \begin{enumerate}
  \item
   A non-negative random variable \ $Y$ \ is regularly varying with index \ $\alpha \in \RR_{++}$ \ if and only if \ $\PP(Y > x) \in \RR_{++}$ \ for all \ $x \in \RR_{++}$, \ and the function \ $\RR_{++} \ni x \mapsto \PP(Y > x) \in \RR_{++}$ \ is regularly varying at infinity with index \ $-\alpha$.
  \item
   If \ $Y$ \ is a regularly varying random variable with index \ $\alpha \in \RR_{++}$, \ then for each \ $\beta \in \RR_{++}$, \ $|Y|^\beta$ \ is regularly varying with index \ $\alpha/\beta$.
 \end{enumerate}
\end{Lem}

The proof of the following lemma can be found, e.g., in Barczy et al. \cite[Lemma C.5]{BarNedPap3_a}.

\begin{Lem}\label{a_n}
If \ $Y$ \ is a regularly varying random variable with index \ $\alpha \in \RR_{++}$, \ then there exists a sequence \ $(a_n)_{n\in\NN}$ \ in \ $\RR_{++}$ \ such that \ $n \PP(|Y| > a_n) \to 1$ \ as \ $n \to \infty$.
\ If \ $(a_n)_{n\in\NN}$ \ is such a sequence, then \ $a_n \to \infty$ \ as \ $n \to \infty$.
\end{Lem}

\begin{Lem}\label{Lem_an_tan}
If \ $Y$ \ is a regularly varying random variable with index \ $\alpha \in \RR_{++}$ \ and \ $(a_n)_{n\in\NN}$ \ is a sequence in \ $\RR_{++}$ \ such that \ $n \PP(|Y| > a_n) \to 1$ \ as \ $n \to \infty$, \ then for each \ $c \in \RR$, \ the random variable \ $Y - c$ \ is regularly varying with index \ $\alpha$, \ and \ $n \PP(|Y - c| > a_n) \to 1$ \ as \ $n \to \infty$.
\end{Lem}

\noindent{\bf Proof.}
Let\ $c \in \RR$.
\ Then \ $Y - c$ \ is regularly varying with index \ $\alpha$, \ see, e.g., part (i) in Lemma C.3.1 in Buraczewski et al.\ \cite{BurDamMik}.
\ By Lemma \ref{a_n}, \ $a_n \to \infty$ \ as \ $n \to \infty$, \ hence for sufficiently large \ $n\in\NN$, \ we have
 \begin{align*}
  &n\PP(|Y - c| > a_n)
   = n \PP(|Y| > a_n) \frac{\PP(Y - c > a_n)+\PP(Y - c < -a_n)}{\PP(|Y| > a_n)} \\
  &= n \PP(|Y| > a_n)
     \biggl(\frac{\PP(Y>c+a_n)}{\PP(|Y|>c+a_n)} \frac{\PP(|Y|>a_n(1+c/a_n))}{\PP(|Y|>a_n)} \\
  &\phantom{= n \PP(|Y| > a_n)\biggl(}
            + \frac{\PP(Y<c-a_n)}{\PP(|Y|>a_n-c)} \frac{\PP(|Y|>a_n(1-c/a_n))}{\PP(|Y|>a_n)}\biggr) .
 \end{align*}
By the uniform convergence theorem for regularly varying functions (see, e.g., Bingham et al.\ \cite[Theorem 1.5.2]{BinGolTeu}) together with the fact that \ $1 + c/a_n \in [1/2, 2]$ \ and \ $1 - c/a_n \in [1/2, 2]$ \ for sufficiently large \ $n \in \NN$, \ we obtain
 \[
   \lim_{n\to\infty} \frac{\PP(|Y|>a_n(1+c/a_n))}{\PP(|Y|>a_n)} = 1 , \qquad
   \lim_{n\to\infty} \frac{\PP(|Y|>a_n(1-c/a_n))}{\PP(|Y|>a_n)} = 1 .
 \]
Hence, by the tail-balance condition \eqref{TB} and Remark \ref{tail-balance}, we conclude
 \[
   \lim_{n\to\infty} n\PP(|Y - c| > a_n)
   = \lim_{n\to\infty} \frac{\PP(Y>c+a_n)}{\PP(|Y|>c+a_n)}
     + \lim_{n\to\infty} \frac{\PP(Y<c-a_n)}{\PP(|Y|>a_n-c)}
   = p+q= 1 ,
 \]
 as desired.
\proofend

\begin{Lem}[Karamata's theorem for truncated moments]\label{truncated_moments}
Consider a non-negative regularly varying random variable \ $Y$ \ with index \ $\alpha \in \RR_{++}$.
\ Then
 \begin{align*}
  \lim_{x\to\infty}
   \frac{x^\beta\PP(Y>x)}{\EE(Y^\beta\bbone_{\{Y\leq x\}})}
  &= \frac{\beta-\alpha}{\alpha} \qquad \text{for \ $\beta \in [\alpha, \infty)$,} \\
  \lim_{x\to\infty}
   \frac{x^\beta\PP(Y>x)}{\EE(Y^\beta\bbone_{\{Y>x\}})}
  &= \frac{\alpha-\beta}{\alpha} \qquad \text{for \ $\beta \in (-\infty, \alpha)$.}
 \end{align*}
\end{Lem}

For Lemma \ref{truncated_moments}, see, e.g., Bingham et al.\ \cite[pages 26-27]{BinGolTeu} or Buraczewski et al.\ \cite[Appendix B.4]{BurDamMik}.

Next, based on Buraczewski et al.\ \cite[Appendix C]{BurDamMik}, we recall the definition
 and some properties of regularly varying random vectors.

\begin{Def}
A \ $d$-dimensional random vector \ $\bY$ \ and its distribution are called regularly
 varying with index \ $\alpha \in \RR_{++}$
 \ if there exists a probability measure \ $\psi$ \ on \ $\SSS^{d-1} := \{\bx \in \RR^d : \|\bx\| = 1\}$ \ such that for all \ $c \in \RR_{++}$,
 \[
   \frac{\PP\bigl(\|\bY\| > c x, \, \frac{\bY}{\|\bY\|} \in \cdot\bigr)}
        {\PP(\|\bY\| > x)}
   \distrw c^{-\alpha} \psi(\cdot) \qquad \text{as \ $x \to \infty$,}
 \]
 where \ $\distrw$ \ denotes the weak convergence of finite measures on \ $\SSS^{d-1}$.
\ The probability measure \ $\psi$ \ is called the spectral measure of \ $\bY$.
\end{Def}

The following equivalent characterization of multivariate regular variation can be derived, e.g., from Resnick \cite[page 69]{Res0},
 see Barczy et al. \cite[Proposition C.8]{BarNedPap3_a}.

\begin{Pro}\label{vague}
A $d$-dimensional random vector \ $\bY$ \ is regularly varying with some index \ $\alpha \in \RR_{++}$ \
 if and only if there exists a non-null locally finite measure \ $\mu$ \ on \ $\RR^d_0$ \ satisfying the limit relation
 \begin{equation}\label{vague_Kallenberg}
   \mu_x(\cdot)
   := \frac{\PP(x^{-1} \bY \in \cdot)}{\PP(\|\bY\| > x)}
   \distrv \mu(\cdot) \qquad \text{as \ $x \to \infty$,}
 \end{equation}
 where \ $\distrv$ \ denotes vague convergence of locally finite measures on \ $\RR^d_0$ \ (see Appendix \ref{App_vague}
  for the notion \ $\distrv$).
\ Further, \ $\mu$ \ satisfies the property \ $\mu(cB) = c^{-\alpha}\mu(B)$ \ for any \ $c\in\RR_{++}$ \ and
 \ $B\in\cB(\RR_0^d)$ \ (see, e.g., Theorems 1.14 and 1.15 and Remark 1.16 in Lindskog \cite{Lin}).
\end{Pro}
The measure \ $\mu$ \ in Proposition \ref{vague} is called the limit measure of \ $\bY$.

The next statement follows, e.g., from part (i) in Lemma C.3.1 in Buraczewski et al.\ \cite{BurDamMik}.

\begin{Lem}\label{Lem_shift}
If \ $\bY$ \ is a regularly varying \ $d$-dimensional random vector with index \ $\alpha \in \RR_{++}$, \ then for each \ $\bc \in \RR^d$, \ the random vector \ $\bY - \bc$ \ is regularly varying with index \ $\alpha$.
\end{Lem}

Recall that if \ $\bY$ \ is a regularly varying \ $d$-dimensional random vector with index \ $\alpha \in \RR_{++}$ \ and with limit measure \ $\mu$ \ given in \eqref{vague_Kallenberg}, and \ $f: \RR^d \to \RR$ \ is a continuous function with \ $f^{-1}(\{0\}) = \{\bzero\}$ \ and it is positively homogeneous of degree \ $\beta \in \RR_{++}$ \ (i.e., \ $f(c \bv) = c^\beta f(\bv)$ \ for every \ $c \in \RR_{++}$ \ and \ $\bv \in \RR^d)$, \ then \ $f(\bY)$ \ is regularly varying with index \ $\frac{\alpha}{\beta}$ \ and with limit measure \ $\mu(f^{-1}(\cdot))$, \ see, e.g., Buraczewski et al.\ \cite[page 282]{BurDamMik}.
Next we recall a result on the tail behaviour of \ $f(\bY)$ \ for appropriate positively homogeneous functions \ $f: \RR^d \to \RR$, \
 see Barczy et al. \cite[Proposition C.10]{BarNedPap3_a}.

\begin{Pro}\label{Pro_mapping}
Let \ $\bY$ \ be a regularly varying \ $d$-dimensional random vector with index \ $\alpha \in \RR_{++}$ \ and let \ $f: \RR^d \to \RR$ \ be a measurable function which is positively homogeneous of degree \ $\beta \in \RR_{++}$, \ continuous at \ $\bzero$ \ and \ $\mu(D_f) = 0$, \ where \ $\mu$ \ is the limit measure of \ $\bY$ \ given in \eqref{vague_Kallenberg} and \ $D_f$ \ denotes the set of discontinuities of \ $f$.
\ Then \ $\mu(\partial_{\RR_0^d}(f^{-1}((1, \infty)))) = 0$, \ where
\ $\partial_{\RR_0^d}(f^{-1}((1, \infty)))$ \ denotes the boundary of \ $f^{-1}((1, \infty))$ \ in \ $\RR_0^d$.
\ Consequently,
 \[
   \lim_{x\to\infty} \frac{\PP(f(\bY)>x)}{\PP(\|\bY\|^\beta>x)}
   = \mu(f^{-1}((1,\infty))) ,
 \]
 and \ $f(\bY)$ \ is regularly varying with tail index \ $\frac{\alpha}{\beta}$.
\end{Pro}

\section{Strongly stationary sequences}
\label{App_stat_proc}

In the next remark we check that for a strongly mixing, strongly stationary
 sequence its rate function given in \eqref{help_mixing_rate} can be represented in another form.

\begin{Rem}\label{Rem_strongly_mixing}
If \ $(Y_k)_{k\in\NN}$ \ is a strongly mixing, strongly stationary sequence with a rate function \ $(\alpha_h)_{h\in\NN}$ \ given in \eqref{help_mixing_rate}
 \ and \ $(\ldots,Y_{-2},Y_{-1},Y_0,Y_1,Y_2,\ldots)$ \ is a strongly stationary extension, then for each \ $h \in \NN$,
 \ we have \ $\alpha_h = \talpha_h$, \ and
 \[
   \talpha_h:= \sup_{k\in\NN} \sup_{A\in\cF^Y_{1,k},\; B\in\cF^Y_{k+h,\infty}} |\PP(A \cap B) - \PP(A) \PP(B)|
 \]
 with \ $\cF^Y_{i,j} := \sigma(Y_i, \ldots, Y_j)$ \ for \ $i, j \in \ZZ$ \ with \ $i \leq j$.
\ Indeed, for each \ $k, h \in \NN$, \ by the strong stationarity, we have
 \[
   \sup_{A\in\cF^Y_{1,k},\; B\in\cF^Y_{k+h,\infty}} |\PP(A \cap B) - \PP(A) \PP(B)|
   = \sup_{A\in\cF^Y_{1-k,0},\; B\in\cF^Y_{h,\infty}} |\PP(A \cap B) - \PP(A) \PP(B)|
   \leq \alpha_h ,
 \]
 hence \ $\talpha_h \leq \alpha_h$.
\ Moreover, for each \ $h \in \NN$, \ we have
 \[
   \alpha_h = \sup_{\ell\in\ZZ_+} \sup_{A\in\cF^Y_{-\ell,0},\; B\in\cF^Y_{h,\infty}} |\PP(A \cap B) - \PP(A) \PP(B)| ,
 \]
 where, for each \ $\ell \in \ZZ_+$ \ and \ $h \in \NN$, \ again by the strong stationarity, we have
 \[
   \sup_{A\in\cF^Y_{-\ell,0},\; B\in\cF^Y_{h,\infty}} |\PP(A \cap B) - \PP(A) \PP(B)|
   = \sup_{A\in\cF^Y_{1,\ell+1},\; B\in\cF^Y_{h+\ell+1,\infty}} |\PP(A \cap B) - \PP(A) \PP(B)|
   \leq \talpha_h ,
 \]
 hence \ $\alpha_h \leq \talpha_h$, \ and we conclude \ $\alpha_h = \talpha_h$, \ as claimed.
\proofend
\end{Rem}

\begin{Lem}\label{Lem_process_shift}
Let \ $(Y_k)_{k\in\ZZ_+}$ \ be a stochastic process, and let \ $c \in \RR$.
\renewcommand{\labelenumi}{{\rm(\roman{enumi})}}
 \begin{enumerate}
  \item
   If \ $(Y_k)_{k\in\ZZ_+}$ \ is strongly stationary, then \ $(Y_k - c)_{k\in\ZZ_+}$ \ is strongly stationary.
  \item
   If \ $(Y_k)_{k\in\ZZ_+}$ \ is jointly regularly varying with index \ $\alpha \in \RR_{++}$, \ then \ $(Y_k - c)_{k\in\ZZ_+}$ \ is jointly regularly varying with index \ $\alpha$.
  \item
   If \ $(Y_k)_{k\in\ZZ_+}$ \ is strongly stationary and strongly mixing with geometric rate function, then \ $(Y_k - c)_{k\in\ZZ_+}$ \ is strongly stationary and strongly mixing with the same geometric rate function.
 \end{enumerate}
\end{Lem}

\noindent{\bf Proof.}
(i). \ For each \ $k \in \NN$, \ let \ $g_k(x_1, \ldots, x_k) := (x_1 - c, \ldots, x_k - c)^\top$ \ for \ $(x_1, \ldots, x_k)^\top \in \RR^k$.
\ For each \ $k, m \in \NN$, \ using the strong stationarity of \ $(Y_k)_{k\in\ZZ_+}$, \ we have
 \[
   (Y_{m+1} - c, \ldots, Y_{m+k} - c)^\top
   = g_k(Y_{m+1}, \ldots, Y_{m+k})
   \distre g_k(Y_1, \ldots, Y_k)
   = (Y_1 - c, \ldots, Y_k - c)^\top ,
 \]
 hence \ $(Y_k - c)_{k\in\ZZ_+}$ \ is strongly stationary.

(ii). \ This is a consequence of Lemma \ref{Lem_shift}.

(iii). \ For each \ $k \in \NN$, \ we have \ $\sigma(Y_1 - c, \ldots, Y_k - c) = \sigma(g_k(Y_1, \ldots, Y_k)) \subset \sigma(Y_1, \ldots, Y_k)$, \ and, in a similar way, \ $\sigma(Y_1, \ldots, Y_k) \subset \sigma(Y_1 - c, \ldots, Y_k - c)$, \ and hence \ $\sigma(Y_1 - c, \ldots, Y_k - c) = \sigma(Y_1, \ldots, Y_k)$.
\ By the same reason, \ $\sigma(Y_j - c : j \geq k) = \sigma(Y_j : j \geq k)$ \ for each \ $k \in \NN$.
\ Using Remark \ref{Rem_strongly_mixing}, the strong mixing property of \ $(Y_k - c)_{k\in\ZZ_+}$ \ readily follows from the strong mixing property of \ $(Y_k)_{k\in\ZZ_+}$, \ and the rate functions of \ $(Y_k)_{k\in\ZZ_+}$ \ and \ $(Y_k - c)_{k\in\ZZ_+}$ \ coincide.
\proofend

\section{A representation of \ $(X_k)_{k\in\ZZ}$}
\label{App_representation}

Let $(\ldots,X_{-2},X_{-1},X_0,X_1,X_2,\ldots)$ be a strongly stationary extension of $(X_j)_{j\in\ZZ_+}$ given in Section \ref{application}.
The following representation of $(X_k)_{k\in\ZZ}$ can be found in Barczy et al.\ \cite[Lemma E.2]{BarBasKevPlaPap}.

\begin{Lem}\label{lem:reprX}
We have
 \begin{equation}\label{repr}
  (X_k)_{k\in\ZZ}
  \distre \biggl(\vare_k + \sum_{i=1}^\infty \theta_k^{(k-i)} \circ \cdots \circ \theta_{k-i+1}^{(k-i)} \circ \vare_{k-i}\biggr)_{k\in\ZZ} ,
 \end{equation}
 where \ $\{\vare_k : k \in \ZZ\}$ \ are independent random variables with the same distribution as \ $\vare$, \ and \ $\theta_k^{(\ell)}$, \ $k, \ell \in \ZZ$, are given by
 \[
   \theta_k^{(\ell)} \circ i
    := \begin{cases}
       \sum_{j=1}^i \xi_{k,j}^{(\ell)} , & \text{for \ $i \in \NN$,} \\
       0 , & \text{for \ $i = 0$,}
      \end{cases}
 \]
 where \ $\xi_{k,j}^{(\ell)}$, \ $j \in \NN$, \ $k, \ell \in \ZZ$, \ have the same distribution as \ $\xi$, \ and \ $\{\vare_k : k \in \ZZ\}$ \ and \ $\theta_k^{(\ell)}$, \ $k, \ell \in \ZZ$, \ are independent in the sense that the families \ $\{\vare_k : k \in \ZZ\}$ \ and \ $\{\xi_{k,j}^{(\ell)} : j \in \NN\}$, \ $k, \ell \in \ZZ$, \ occurring in \ $\theta_k^{(\ell)}$, \ $k, \ell \in \ZZ$, \ are independent families of independent random variables, and the series in the representation \eqref{repr} converge with probability one.
\end{Lem}

Next, we recall a useful representation of the random vectors \ $(X_0, X_1, \ldots, X_n)$, \ $n \in \NN$, \ from the proof of Theorem E.3
 in Barczy et al.\ \cite{BarBasKevPlaPap}.

\begin{Lem}\label{lem:reprXn}
For each \ $n \in \NN$, \ we have
 \[
   (X_0, X_1, \ldots, X_n)
   \distre
   (X_0, \kappa_1 + \theta_1^{(0)} \circ X_0, \ldots,
    \kappa_n + \theta_n^{(0)} \circ \cdots \circ \theta_1^{(0)} \circ X_0) ,
 \]
 where
 \[
   \kappa_k := \vare_k + \sum_{i=1}^{k-1} \theta_k^{(k-i)} \circ \cdots \circ \theta_{k-i+1}^{(k-i)} \circ \vare_{k-i} , \qquad k \in \NN .
 \]
Moreover, \ $\kappa_k$, \ $\theta_k^{(0)} \circ \cdots \circ \theta_1^{(0)} \circ j$ \ and \ $X_0$ \ are independent for any \ $k \in \NN$ \ and \ $j \in \ZZ_+$.
\end{Lem}

\section{Tail behaviour of \ $(X_k)_{k\in\ZZ_+}$}

Due to Basrak et al.\ \cite[Theorem 2.1.1]{BasKulPal}, we have the following tail
 behaviour of \ $(X_k)_{k\in\ZZ_+}$ \ given in Section \ref{application}.

\begin{Thm}\label{Xtail}
We have
 \[
   \lim_{x\to\infty} \frac{\pi((x, \infty))}{\PP(\vare > x)}
   = \sum_{i=0}^\infty m_\xi^{i\alpha}
   = \frac{1}{1-m_\xi^\alpha} ,
 \]
 where \ $\pi$ \ denotes the unique stationary distribution of the Markov chain \ $(X_k)_{k\in\ZZ_+}$, \ and consequently, \ $\pi$ \ is also regularly varying with index \ $\alpha$.
\end{Thm}

Note that in case of \ $\alpha = 1$ \ and \ $m_\vare = \infty$ \ Basrak et al.\
 \cite[Theorem 2.1.1]{BasKulPal} assume additionally that \ $\vare$ \ is consistently varying (or
 in other words intermediate varying), but, eventually, it follows from the fact that \ $\vare$ \ is regularly
 varying.

Let $(\ldots,X_{-2},X_{-1},X_0,X_1,X_2,\ldots)$  be a strongly stationary extension of $(X_k)_{k\in\ZZ_+}$.
Basrak et al. \cite[Lemma 3.1]{BasKulPal} described the so-called forward tail process of the
 strongly stationary process \ $(X_k)_{k\in\ZZ}$, \ and hence, due to Basrak and
 Segers \cite[Theorem 2.1]{BasSeg}, the strongly stationary process
 \ $(X_k)_{k\in\ZZ}$ \ is jointly regularly varying.
We summarize these results in the following theorem.

\begin{Thm}\label{Xtailprocess}
The finite dimensional conditional distributions of \ $(x^{-1} X_k)_{k\in\ZZ_+}$ \ with respect
 to the condition \ $X_0 > x$ \ converge weakly to the corresponding finite dimensional distributions of
 \ $(m_\xi^k Y)_{k\in\ZZ_+}$ \ as \ $x \to \infty$, \ where \ $Y$ \ is a random
 variable with Pareto distribution
 \ $\PP(Y \leq y) = (1 - y^{-\alpha}) \bbone_{[1,\infty)}(y)$, \ $y \in \RR$.
\ Consequently, the strongly stationary process \ $(X_k)_{k\in\ZZ}$ \ is jointly
 regularly varying with index \ $\alpha$, \ i.e., all its finite dimensional
 distributions are regularly varying with index \ $\alpha$.
\ The processes \ $(m_\xi^k Y)_{k\in\ZZ_+}$ \ and \ $(m_\xi^k)_{k\in\ZZ_+}$ \ are the so-called forward tail process and forward spectral tail process of \ $(X_k)_{k\in\ZZ}$, \ respectively.
\end{Thm}

\section{Mixing property of \ $(X_k)_{k\in\ZZ_+}$}
\label{mixing}

First, we present an auxiliary lemma stating that \ $(X_k)_{k\in\ZZ_+}$ \ given in Section \ref{application} is strongly mixing with geometric rate,
 see also Basrak et al.\ \cite[Remark 3.1]{BasKulPal}.

\begin{Lem}\label{lemma:strong_mixing}
The strongly stationary Markov chain \ $(X_k)_{k\in\ZZ_+}$ \ is aperiodic, strongly mixing with geometric rate, i.e., there exists a constant \ $q \in (0, 1)$ \ such that \ $\alpha_h = \OO(q^h)$ \ as \ $h \to \infty$.
\end{Lem}

\noindent{\bf Proof.}
We will apply part 1 of Theorem 2 in Jones \cite{Jon} in order to prove that \ $(X_k)_{k\in\ZZ_+}$ \ is strongly mixing.
For this, we need to check that \ $(X_k)_{k\in\ZZ_+}$ \ is aperiodic, \ $\psi$-irreducible and positive Harris recurrent,
 for the definitions, see Meyn and Tweedie \cite[pages 114, 84, 199 and 231]{MeyTwe}.
Since \ $m_\xi \in [0, 1)$, \ $\PP(\vare = 0) < 1$ \ and \ $\sum_{\ell=1}^\infty \log(\ell) \PP(\vare = \ell) < \infty$, \ there exists a unique stationary distribution \ $\pi$ \ of \ $(X_k)_{k\in\ZZ_+}$.
\ Indeed, one can apply Quine \cite[page 414]{Qui}, since the \ $1 \times 1$-matrix \ $m_\xi$ \ is irreducible and aperiodic
 (in the sense that there does not exist a positive integer \ $k \in \NN$ \ such that \ $m_\xi^{k+1} =  m_\xi$).
\ Since the state space \ $I \subset \ZZ_+$ \ of \ $(X_k)_{k\in\ZZ_+}$ \ is denumerable, the existence of a unique stationary distribution of \ $(X_k)_{k\in\ZZ_+}$ \ yields that there is exactly one positive (ergodic) communication class \ $D$ \ in \ $I$ \ which is the support of the unique stationary distribution \ $\pi$, \ see, e.g., Chung \cite[\S 7, Theorem 2]{chung:1960}.
Since the distribution of \ $X_0$ \ is \ $\pi$, \ then, by the definition of a communication class, \ $I = D$, \ hence \ $(X_k)_{k\in\ZZ_+}$ \ is irreducible and positive recurrent.
Since \ $I$ \ is denumerable, \ $(X_k)_{k\in\ZZ_+}$ \ is \ $\psi$-irreducible in the sense of Meyn and Tweedie \cite[page 84]{MeyTwe} with \ $\psi$ \ being the counting measure, and it is Harris recurrent in the sense of Meyn and Tweedie \cite[page 199]{MeyTwe}.
Next we check that \ $(X_k)_{k\in\ZZ_+}$ \ is aperiodic.
Let \ $i_{\min} := \inf\{\ell \in \ZZ_+ : \PP(\vare = \ell) > 0\}$.
\ Then \ $\PP(X_1 = i_{\min} \mid X_0 = j) > 0$ \ for all \ $j \in I$, \ since \ $m_\xi < 1$ \ yields \ $\PP(\xi = 0) > 0$.
\ Since \ $I$ \ consists of a single communication class, \ $i_{\min} \in I$.
\ Since $\PP(X_1 = i_{\min} \mid X_0 = i_{\min}) > 0$, \ the state \ $i_{\min}$, \ and hence,
 using again that there is a single communication class, the Markov chain \ $(X_i)_{i\in\ZZ_+}$ \ is aperiodic.

We will apply part 2 of Theorem 2 in Jones \cite{Jon} in order to prove that \ $(X_k)_{k\in\ZZ_+}$ \ is strongly mixing with geometric rate.
For this, we need to check that
 \begin{itemize}
  \item
   $(X_k)_{k\in\ZZ_+}$ \ is geometrically ergodic, i.e., there exists a function \ $M : I \to \RR_+$ \ and a constant \ $q \in (0, 1)$ \ such that
   \[
     \|\PP^n(i, \cdot) - \PP_{X_0}(\cdot)\| \leq M(i) q^n , \qquad i \in I , \quad n \in \NN ,
   \]
   where \ $\PP^n(i, A) := \PP(X_n \in A \mid X_0 = i)$, \ $n \in \NN$, \ $i \in I$, \ $A \subset I$, \ and \ $\|\cdot\|$ \ denotes the total variation norm of probability measures (i.e., \ $\|Q_1 - Q_2\| := \sup_{A\subset I} |Q_1(A) - Q_2(A)|$ \ for two probability measures \ $Q_1$ \ and \ $Q_2$),
  \item
   $\EE(M(X_0)) < \infty$.
 \end{itemize}
For geometric ergodicity, it is enough to check a so-called drift condition, namely, there exists a function \ $V : I \to [1, \infty)$, \ constants \ $d \in \RR_{++}$, \ $b \in \RR$, \ and a subset \ $C \subset I$ \ such that
 \begin{align}\label{drift_condition}
  \EE(V(X_1) \mid X_0 = i) - V(i) \leq - d V(i) + b \bbone_C(i) , \qquad i \in I ,
 \end{align}
 where \ $C$ \ is a so-called small set, see, e.g., Jones \cite[page 301, equation (5)]{Jon}.
Note also that condition \eqref{drift_condition} is nothing else but condition (15.28) in Meyn and Tweedie \cite{MeyTwe}.
We check that \eqref{drift_condition} holds with \ $V : I \to [1, \infty)$, \ $V(x) := 1 + x^p$, \ $x \in I$, \ with any \ $p \in (0, \alpha \land 1)$ \ and \ $C := \{0, 1, \ldots, K\} \cap I$ \ with some sufficiently large \ $K \in \NN$.
\ First, we verify that \ $\EE(V(X_1) \mid X_0 = i) \leq (1 - d) V(i)$ \ holds for all \ $i > K$, \ $i \in I$, \ with some \ $K \in \NN$.
\ Since \ $p \in (0, 1)$, \ we have \ $(x + y)^p \leq x^p + y^p$ \ for all \ $x, y \in \RR_+$ \ and the function \ $\RR_+ \ni x \mapsto x^p$ \ is concave, hence by Jensen's inequality,
 \begin{align*}
  \EE(X_1^p \mid X_0 = i)
  &= \EE\left(\left(\sum_{j=1}^i \xi_{1,j} + \vare_1\right)^p\right)
   \leq\EE\left(\left(\sum_{j=1}^i \xi_{1,j}\right)^p + \vare_1^p\right) \\
  &\leq \left(\EE\left(\sum_{j=1}^i \xi_{1,j}\right)\right)^p + \EE(\vare_1^p)
   = i^p m_\xi^p + \EE(\vare^p) , \qquad i \in I .
 \end{align*}
Hence \ $\EE(V(X_1) \mid X_0 = i) \leq 1 + m_\xi^p i^p + \EE(\vare^p)$, \ $i \in I$.
\ Consider an arbitrary constant \ $d \in (0, 1 - m_\xi^p)$.
\ Then \ $1 + m_\xi^p i^p + \EE(\vare^p) \leq (1 - d) V(i) = (1 - d) (1 + i^p)$ \ holds if and only if
 \[
   \EE(\vare^p) + d \leq (1 - d - m_\xi^p)i^p
   \qquad \Longleftrightarrow \qquad \frac{\EE(\vare^p)+d}{1-d-m_\xi^p} \leq i^p ,
 \]
 which holds for sufficiently large \ $i \in I$, \ e.g., one can choose \ $K$ \ to be \ $\Bigl\lfloor \left(\frac{\EE(\vare^p)+d}{1-d-m_\xi^p}\right)^{1/p}\Bigr\rfloor + 1$.
\ Next, we verify that \ $C$ \ is a small set (and consequently a petite set as well).
For all \ $i\in C$ \ and \ $B\subset \ZZ_+$, \ we have
 \[
   \PP(i,B) = \PP(X_1\in B \mid X_0=i)
            = \PP\left(\sum_{j=1}^i \xi_{1,j} + \vare_1 \in B\right)
            \geq \PP\left( \sum_{j=1}^K \xi_{1,j}=0, \vare_1\in B \right)
            =:\mu(B),
 \]
 where \ $\mu$ \ is a non-trivial finite measure on \ $(\ZZ_+,2^{\ZZ_+})$, \ since \ $\PP\left( \sum_{j=1}^K \xi_{1,j}=0 \right) = \PP(\xi_{1,1}=0)^K>0$ \
 following from the fact that \ $(X_k)_{k\in\ZZ_+}$ \ is subcritical.
Finally, the constant \ $b$ \ can be chosen as \ $\sup_{i\in C} |\EE_i(V(X_1)) + (1 - d)V(i)|$, \ which is finite due to the facts that \ $\EE(V(X_1) \mid X_0 = i) < \infty$, \ $i \in I$, \ and \ $V(i) < \infty$, \ $i \in I$.
\ Due to Remark 1 in Jones \cite{Jon}, \ $\EE(M(X_0)) < \infty$ \ holds, as well.
So we conclude that \ $(X_k)_{k\in\ZZ_+}$ \ is strongly mixing with geometric rate.
\proofend

For probability measures \ $\nu_1$ \ and \ $\nu_2$ \ on \ $(\ZZ_+,2^{\ZZ_+})$ \ and for a function
 \ $V:\ZZ_+\to [1,\infty)$, \ let
 \[
  \Vert \nu_1 - \nu_2\Vert_V:=\sup_{g:\ZZ_+\to\RR, \, \vert g\vert\leq V} \left\vert \int_{\ZZ_+} g(n)(\nu_1-\nu_2)(\dd n)\right\vert
                             = \sup_{g:\ZZ_+\to\RR, \, \vert g\vert\leq V} \left\vert \sum_{n=0}^\infty g(n)(\nu_1(\{n\}) - \nu_2(\{n\}) )  \right\vert.
 \]
Using Theorem 16.0.1 in Meyn and Tweedie \cite{MeyTwe} (equivalences of (ii) and (iv)) and the proof of Lemma \ref{lemma:strong_mixing},
 we have the following corollary.

\begin{Cor}\label{Cor_ergod}
The strongly stationary Markov chain \ $(X_k)_{k\in\ZZ_+}$ \ is \ $V$-uniformly ergodic, where \ $V:I\to[1,\infty)$, \ $V(x):=1+x^p$, \ $x\in I$, \
 with any \ $p\in(0,\alpha\wedge 1)$, \ and \ $I$ \ denotes the state space of \ $(X_k)_{k\in\ZZ_+}$, \ i.e.,
 there exist \ $\varrho\in(0,1)$ \ and \ $C>0$ \ such that
 \[
    \Vert \PP^n(i,\cdot) - \pi\Vert_V \leq C V(i)\varrho^n,\qquad n\in\ZZ_+, \;i\in I,
 \]
 where \ $\pi$ \ denotes the unique stationary distribution of \ $(X_k)_{k\in\ZZ_+}$.
\end{Cor}

\section*{Acknowledgements}
We would like to thank Thomas Mikosch for his suggestion to use the anti-clustering type condition \eqref{3.9} presented in Lemma \ref{Lemma2}, which will
 appear in his forthcoming book \cite{MikWin2} written jointly with Olivier Wintenberger.
We would like to thank the referee for the comments that helped us improve the paper.

\bibliographystyle{plain}
\bibliography{reg_var_aggr_3}

\end{document}